\pdfoutput=1
\documentclass[11pt]{amsart}
\usepackage[margin=1in]{geometry}
\usepackage{amsmath,amssymb,mathtools,mathrsfs,microtype,enumitem,graphicx,needspace}
\usepackage{tikz}
\usetikzlibrary{arrows.meta,fit}
\usepackage[colorlinks=true,allcolors=blue]{hyperref}
\numberwithin{equation}{section}
\newtheorem{theorem}{Theorem}[section]
\newtheorem{proposition}[theorem]{Proposition}
\newtheorem{lemma}[theorem]{Lemma}
\newtheorem{corollary}[theorem]{Corollary}
\newtheorem{remark}[theorem]{Remark}
\newcommand{\Holder}[1]{C^{#1,1/2}}
\newcommand{\T}{\mathbb T^3}
\newcommand{\R}{\mathbb R}
\newcommand{\dd}{\,\mathrm d}
\newcommand{\Id}{\operatorname{Id}}
\newcommand{\curl}{\operatorname{curl}}
\newcommand{\diver}{\operatorname{div}}
\newcommand{\supp}{\operatorname{supp}}
\newcommand{\tr}{\operatorname{tr}}
\newcommand{\norm}[1]{\left\lVert#1\right\rVert}
\newcommand{\av}[1]{\langle #1\rangle_s}
\newcommand{\PP}{\mathbb P}
\newcommand{\ip}[2]{\langle #1,#2\rangle}

\title[Singularity Formation for the Three-Dimensional Euler Equations]{Finite-Time Singularity Formation\\for the Three-Dimensional Euler Equations}
\hypersetup{pdftitle={Finite-Time Singularity Formation for the Three-Dimensional Euler Equations}}
\author{Levent Alp\"oge}
\author{Tristan Buckmaster}
\author{Matei P. Coiculescu}
\date{}
\begin{document}
\raggedbottom
\begin{abstract}
We construct smooth periodic initial data for the three-dimensional
incompressible Euler equations whose solution blows up in finite time.
The velocity remains bounded, while the vorticity magnitude at the origin
tends to infinity. Away from the origin, the velocity and its time derivative
have smooth spatial limits at the singular time. The construction uses
localized Craik--Criminale waves as building blocks of a multi-scale construction.
\end{abstract}
\maketitle
\setcounter{tocdepth}{1}
\tableofcontents
\section{Introduction}\label{cc:sec:problem}

The incompressible Euler equations are among the earliest partial
differential equations derived in mathematical physics. Euler presented them in
1755, less than a decade after d'Alembert derived the wave equation for
the vibrating string~\cite{DAlembert1747,Euler1757}.

The goal of this work is to construct smooth initial data for the
incompressible Euler equations on the three-dimensional torus whose
solution blows up in finite time. We represent $\mathbb T^3$ by the
cube $[-1/2,1/2]^3$ with periodic boundary conditions
(opposite faces identified). We write these equations as
\begin{equation}\label{cc:eq:euler}
 \partial_tu+(u\cdot\nabla)u+\nabla p=0,
 \qquad \diver u=0
 \quad\hbox{on }\mathbb T^3\times[0,T_*).
\end{equation}
Local existence and uniqueness for sufficiently regular data go back to
Lichtenstein and G\"unther in the 1920s~\cite{Lichtenstein1925,Gunther1926}.
Classical theory on $\mathbb T^3$ supplies a unique maximal smooth solution
from smooth divergence-free data~\cite{EbinMarsden,BedrossianVicol}.
By the Beale--Kato--Majda criterion~\cite{BKM,GGPP}, a smooth solution on
$[0,T)$ continues smoothly through $T<\infty$ if
$\int_0^T\|\curl u(t)\|_{L^\infty}\,\dd t<\infty$.

In two dimensions the scalar vorticity is transported, so its
$L^\infty$ norm is conserved, and smooth solutions exist for all
time~\cite{Wolibner1933}. Yudovich~\cite{Yudovich1963} proved global
existence and uniqueness for bounded, possibly nonsmooth initial
vorticity. In three dimensions, Fefferman's description of the
Navier--Stokes Millennium Problem stated that the global regularity
question was also open for the Euler equations on $\R^3$ and on
$\mathbb T^3$~\cite{Fefferman2006}.

Numerical searches for blow-up include the Taylor--Green vortex study
of Morf, Orszag and Frisch~\cite{MorfOrszagFrisch1980} and Kerr's
computations for antiparallel vortex tubes~\cite{Kerr1993}.
In a later computation, Hou and Li~\cite{HouLi2006} observed at most
double-exponential growth of the maximum vorticity past Kerr's predicted
singular time.
Luo and Hou~\cite{LuoHou2014} found numerical evidence for axisymmetric
blow-up at the wall of an axially periodic cylinder.
Childress, Gilbert and Valiant~\cite{ChildressGilbertValiant2016}
studied, analytically and numerically, eroding vortex dipoles in
axisymmetric flow without swirl, whose maximal vorticity grows like
$t^{4/3}$, and Childress and Gilbert~\cite{ChildressGilbert2018}
proposed, by formal asymptotic analysis, the hairpin geometry as a model
for finite-time blow-up.

On $\R^3$, Elgindi~\cite{ElgindiSingularity} constructed finite-time
blow-up from $C^{1,\alpha_0}$ initial data for some $\alpha_0>0$.
Elgindi, Ghoul and Masmoudi~\cite{ElgindiGhoulMasmoudi} obtained stable,
locally self-similar blow-up for finite-energy data in this regularity class.
Chen and Hou~\cite{ChenHouSmoothI,ChenHouSmoothII} prove nearly self-similar
blow-up from smooth axisymmetric data in a periodic cylinder with boundary,
using computer-assisted estimates. The theorem below concerns smooth
periodic data on $\mathbb T^3$, with bounded velocity and vorticity growth
concentrated at the origin.

\begin{theorem}\label{cc:target:unforced}
There is an odd divergence-free
$u_0\in C^\infty(\mathbb T^3;\R^3)$ whose maximal classical solution $u$
of \eqref{cc:eq:euler} is smooth and odd on
$\mathbb T^3\times[0,T_*)$ for a finite $T_*>0$, and satisfies
\begin{equation}\label{cc:uf:eq:main-breakdown}
 \int_0^{T_*}\|\nabla u(t)\|_{L^\infty}\dd t=\infty,
 \qquad
 \int_0^{T_*}\|\curl u(t)\|_{L^\infty}\dd t=\infty.
\end{equation}
Moreover,
\begin{equation}\label{cc:uf:eq:main-terminal}
 \sup_{0\le t<T_*}\|u(t)\|_{L^\infty}<\infty,
 \qquad \lim_{t\uparrow T_*}|\curl u(t,0)|=\infty.
\end{equation}
For every compact $K\subset\mathbb T^3\setminus\{0\}$ and every integer
$m\ge0$, both $u(t)$ and $\partial_tu(t)$ converge in $C^m(K)$ as
$t\uparrow T_*$.
\end{theorem}

\begin{remark}
The solution in Theorem~\ref{cc:target:unforced} may be chosen so that it
extends smoothly in space and time through $T_*$ away from the origin.
A proof is given in the auxiliary note~\cite{TerminalSmoothness}.
\end{remark}

We construct a sequence of smooth approximate solutions of the Euler
equations whose residuals tend to zero. The sequence converges in every
spatial $C^m$ norm, uniformly on compact time intervals before $T_*$, to
a solution of the Euler equations.

At each step in the sequence we add an oscillation at a higher frequency, localized
near the origin and initially small in velocity. The velocity gradient
of the preceding scales amplifies the new oscillation, whose gradient
in turn drives the growth of the next one. The growth intervals have
finite total length. This multiscale construction follows the approach of C\'ordoba and
Mart\'inez-Zoroa for the Euler equations with forcing and for forced
IPM~\cite{CMZForcedEuler,CMZForcedIPM}.
C\'ordoba, Mart\'inez-Zoroa and Zheng~\cite{CMZUnforcedEuler} used a cascade of
vorticity layers to obtain blow-up for the unforced Euler equations
on $\R^3$ from $C^{1,\alpha}$ data
smooth away from the origin. Related constructions give blow-up with forcing
for the inviscid Boussinesq equations, the variable-density Euler equations
and the hypodissipative Navier--Stokes equations~\cite{CLMZForcedBoussinesq,CLMZNonhomogeneousEuler,CMZHypodissipativeNS}.
Our companion papers proving singularity formation for forced IPM,
the forced Boussinesq equations, and the forced Euler
equations~\cite{AlpogeBuckmasterCoiculescuIPM,AlpogeBuckmasterBoussinesq,AlpogeBuckmasterForcedEuler}
use forces smooth in space and time up to the singular time. In the
present construction we use Craik--Criminale waves as building blocks
and choose the initial velocity at each step so that the limiting
solution is unforced.

Convex integration constructions of nonunique weak solutions of the Euler equations
also use oscillations at successive frequency scales. Below the Onsager
exponent, this method yields solutions with nonconstant kinetic
energy~\cite{IsettOnsager,BDLSVOnsager}. For the unforced Euler equations
in the plane, Coiculescu~\cite{CoiculescuEulerI,CoiculescuEulerII}
uses nested vortex rings at different scales. An instability propagates
from the inner rings outward, yielding distinct renormalized weak
solutions with the same initial vorticity. See also Choi~\cite{ChoiEulerNonuniqueness} and Hou, Qin, Song and
Zhang~\cite{HouQinSongZhangEulerNonuniqueness}.
For the three-dimensional Navier--Stokes equations, Jia and
\v{S}ver\'ak~\cite{JiaSverak2015} gave a conditional construction of
nonunique Leray--Hopf solutions based
on spectral assumptions about forward self-similar flows.
Nonuniqueness was shown for
finite-energy weak solutions by Buckmaster and Vicol~\cite{BuckmasterVicolNS},
and in the Leray--Hopf class with a common force by Albritton, Bru\'e and
Colombo~\cite{AlbrittonBrueColombo}. In a recent preprint, Buckmaster,
Coiculescu and Li~\cite{BuckmasterCoiculescuLiNS} obtain nonuniqueness
of unforced Leray--Hopf solutions.

\subsection{Waves and successive amplification}\label{cc:sec:intro-waves}
We first describe the wave construction of Craik and
Criminale~\cite{CraikCriminale}. Let $D(t)$ be a smooth matrix with
$\operatorname{tr}D=0$ and $D'+D^2$ symmetric. These conditions ensure that
$U(t,x)=D(t)x$, linear in space, solves the Euler equations on $\R^3$, with pressure
$p_D(t,x)=-x\cdot(D'+D^2)x/2$. We add to this linear background a velocity
\[
 w(t,x)=\mathbf a(t)F(k(t)\cdot x),\qquad \mathbf a\cdot k=0,
\]
where $k(t)\ne0$ is the wave vector, $\mathbf a(t)$ is the amplitude
vector, and $F$ is a smooth odd $1$-periodic profile of the scalar phase $s$.
For example, one may take $F(s)=\sin(2\pi s)$. The vector $\mathbf a$
records both size and direction: the velocity amplitude is $|\mathbf a|\|F\|_\infty$, and
$\mathbf a/|\mathbf a|$ is its unit direction when $\mathbf a\ne0$.
The unit wave vector is $k/|k|$. Since $F$ has period one, $w$ is
periodic in the direction of $k$ with period $1/|k|$, so $|k|$ is the
spatial frequency of the wave.

The sum $u=Dx+w$ is an exact solution of the Euler equations if $k$ and $\mathbf a$ solve
\begin{equation}\label{cc:intro:wave-ode}
 k'=-D^Tk,\qquad
 \mathbf a'=-D\mathbf a+p_{\rm w} k,\qquad
 p_{\rm w}=\frac{2k\cdot D\mathbf a}{|k|^2}.
\end{equation}
Its pressure is $p=p_D-p_{\rm w} P(k\cdot x)$, where $P'=F$.
Indeed, $\mathbf a\cdot k=0$ gives $(w\cdot\nabla)w=0$, the first
ordinary differential equation transports the phase by $Dx$, and the
second cancels the remaining acceleration with the gradient of the
wave pressure $-p_{\rm w}P(k\cdot x)$.
The displayed value of $p_{\rm w}$ preserves $\mathbf a\cdot k=0$.
We call \eqref{cc:intro:wave-ode} the Craik--Criminale equations for the
wave vector and the amplitude vector. Normalizing $F'(0)=1$ gives
$\nabla w(t,0)=\mathbf a\otimes k$.
Related short-wave instabilities of planar vortices were studied by
Pierrehumbert~\cite{Pierrehumbert1986} and Bayly~\cite{Bayly1986}.

To describe amplification of a new wave, write $k_p$ and $\mathbf a_p$
for the nonzero wave vector and amplitude vector of the preceding wave,
and set
\[
 E_1=\frac{k_p}{|k_p|},\qquad
 E_2=\frac{\mathbf a_p}{|\mathbf a_p|},\qquad E_3=E_1\times E_2,
 \qquad \sigma=|\mathbf a_p|\,|k_p|.
\]
$E_1$ is the unit wave vector and $E_2$ the direction vector of the
preceding wave.
The preceding wave contributes the velocity gradient
$\mathbf a_p\otimes k_p=\sigma E_2\otimes E_1$ at the origin.
Here $(E_2\otimes E_1)v=(E_1\cdot v)E_2$: the corresponding linear
velocity is $\sigma(E_1\cdot x)E_2$, a shear in direction $E_2$
varying in direction $E_1$. The new wave therefore sees
$D=B+\sigma E_2\otimes E_1$, where $B$ is the older background gradient.
To model the wave equations at the origin, take $B$ with just two
nonzero entries in this moving frame, both constant,
\[
 E_1\cdot BE_2=\nu>0,\qquad E_3\cdot BE_2=\eta>0,
 \qquad \varepsilon=\eta/\nu\ll1.
\]
Start the new wave at time $t=0$ with $k=\lambda E_3$, where $\lambda>0$
is its initial frequency, and write $a_i=\mathbf a\cdot E_i$. The
preceding wave evolves by \eqref{cc:intro:wave-ode} with $D=B$.
Differentiating its normalized vectors gives, in this case,
$E_1'=-\nu E_2$, $E_2'=\nu E_1-\eta E_3$ and $E_3'=\eta E_2$.
Combining these frame equations with \eqref{cc:intro:wave-ode} for the
new wave gives, at $t=0$,
\[
 \dot a_1=-2\nu a_2,\qquad
 \dot a_2=-(\sigma-\nu)a_1.
\]
Freezing these coefficients gives a growing mode with rate
$\sqrt{2\nu(\sigma-\nu)}$ when $\sigma>\nu$. We choose the initial
direction of the amplitude along $E_2-\sqrt{2\nu/\sigma}\,E_1$, the growing direction
to leading order for $\sigma\gg\nu$.

\begin{figure}[tbp]
\centering
\includegraphics[width=\textwidth]{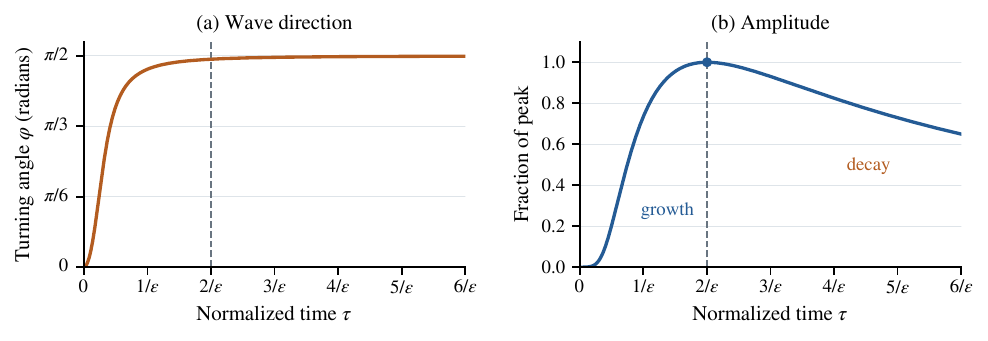}
\caption{Turning and growth in the case of two nonzero entries, in the
limit of a strong shear.
Left: the angle $\varphi$ of the wave vector from $E_3$ toward $E_1$.
Right: amplitude as a fraction of its peak. The dashed line marks
the amplitude peak, by which time $\varphi$ is close to $\pi/2$.}
\label{cc:fig:wave-turning}
\end{figure}

The coefficients evolve as the wave vector turns. Let $\varphi$
be the angle of $k$ from $E_3$ toward $E_1$:
\[
 \tan\varphi=\frac{k\cdot E_1}{k\cdot E_3},
\]
so $\varphi=0$ at $t=0$. The coupling $\eta$ tilts $k$ toward $-E_2$,
and the shear carries this tilt into $E_1$, so $\varphi$ increases;
while $\eta t$ is small, $\tan\varphi$ grows like $(\sigma+\nu)\eta t^2$.
In the strong-shear limit $\sigma/\nu\to\infty$, the wave pressure
multiplies the driving coefficient $2\nu$ in the equation for $a_1$ by
$\cos^2\varphi\,(1-\varepsilon\tan\varphi)$
(Lemma~\ref{cc:lem:strong-shear}). This factor changes sign
at $\tan\varphi=1/\varepsilon$, and in the strong-shear limit the
amplitude $|\mathbf a|$ reaches its peak after this sign change, as shown
in Figure~\ref{cc:fig:wave-turning}.
For the perturbed wave with this initial direction,
Proposition~\ref{cc:prop:finite-return} gives a \emph{gain} of at least
a fixed multiple of $\exp(\varepsilon^{-1/2})$ on its prescribed interval.
The gain is the ratio of the maximum of $|\mathbf a|$ on that interval
to its initial value.

Size alone does not determine the successor wave. The new wave must
supply both a large shear and a suitable frame for the next step of
the iteration. In its frame (the subscript $+$ denotes the frame of
the newly added wave),
\[
 E_{1,+}=k/|k|,\qquad E_{2,+}=\mathbf a/|\mathbf a|,\qquad
 E_{3,+}=E_{1,+}\times E_{2,+},
\]
its shear is $\sigma_+E_{2,+}\otimes E_{1,+}$, with
$\sigma_+=|\mathbf a|\,|k|$. The older gradient $D$, expressed in this
new frame, must again have positive couplings
\[
 \nu_+=E_{1,+}\cdot D E_{2,+},\qquad
 \eta_+=E_{3,+}\cdot D E_{2,+}.
\]
These play the roles of $\nu$ and $\eta$ in the next amplification
step. We prescribe a smaller ratio $\eta_+/\nu_+=\varepsilon_+<\varepsilon$
to give the next wave enough amplification to keep the errors small.
The turning of the wave vector supplies this orientation. In the
strong-shear limit, $\nu_+$ starts at $\varepsilon\nu$, and $\eta_+$
stays within a factor $1+O(\varepsilon)$ of $\nu\cos\varphi$; at late
times this results from a cancellation between terms much larger than
$\nu\cos\varphi$. For the perturbed wave,
Lemma~\ref{cc:lem:new-wave-couplings} shows that at late times $\nu_+$
is close to $\nu$ and $\eta_+/\nu_+$ is comparable to $\cos\varphi$.
Under the hypotheses of
Proposition~\ref{cc:prop:finite-return}, the ratio reaches the prescribed
value after the amplitude peaks in the strong-shear limit, while the new
shear remains large. We start the successor when the ratio reaches that threshold.
Thus the Craik--Criminale equations produce both amplification and the
orientation needed to repeat the construction.

\Needspace{10\baselineskip}
The periodic construction starts from the stationary shear flow
\[
 U_0(x)=\sin(2\pi x_2)(e_1+\rho e_3),
\]
where $e_i$ are the coordinate vectors and $\rho>0$ is small. It solves
the Euler equations with zero pressure.
In the initial coordinate frame its central gradient has
$\nu=2\pi$ and $\eta=2\pi\rho$, so that $\varepsilon=\rho$. Its explicit flow supplies the first wave;
subsequent waves use the successive gradients they create.

Craik--Criminale waves on $\R^3$ provide a local model for amplification.
For the periodic construction, we localize each wave by a smooth cutoff
supported in the interior of a fundamental cell in Lagrangian
coordinates. We extend the localized field periodically and transport its
phase and cutoff by the preceding periodic velocity, as in geometric optics
\cite{LifschitzHameiri1991,FriedlanderVishik1991,Lax1957,ChoquetBruhat1969}.
Since the cutoff
vanishes near the cell boundary, the wave vector need not lie
in $\mathbb Z^3$. Localization and the spatial variation of the background
introduce errors, which the corrections below make small at each step
and remove in the limit.

\subsection{Controlling corrections in an amplifying flow}
To correct a velocity in an amplifying flow, we use the freedom to choose
its initial value. To avoid exciting growing modes, we prescribe the final
particle positions and choose the initial velocity to reach them.

Let $u,p$ solve the Euler equations on $[0,T]$. The flow $X(t,a)$ of
$u$ gives the position at time $t$ of the particle that starts at $a$.
It satisfies
\[
 \partial_tX(t,a)=u(t,X(t,a)),\qquad X(0,a)=a.
\]
For a divergence-free velocity correction $b$, let $X_\delta$ be the
flow of $u+\delta b$ with the same initial positions, so that
\[
 \partial_tX_\delta(t,a)
       =(u+\delta b)(t,X_\delta(t,a)),\qquad X_\delta(0,a)=a.
\]
The first-order change $\xi$ in the particle trajectories is defined by
\[
 X_\delta(t,a)=X(t,a)+\delta\xi(t,a)+O(\delta^2).
\]
Differentiating the flow equation at $\delta=0$ gives
\begin{equation}\label{cc:intro:first-variation}
 \partial_t\xi(t,a)=\nabla u(t,X(t,a))\xi(t,a)+b(t,X(t,a)),
 \qquad \xi(0,a)=0.
\end{equation}
The initial condition follows from the common starting positions.
Since also $X(0,a)=a$, evaluating \eqref{cc:intro:first-variation} at
$t=0$ gives $\partial_t\xi(0,a)=b(0,a)$. This identity holds for every
$b$ and imposes no condition. We impose instead the endpoint condition
$\xi(T,a)=0$. It states that $X_\delta(T,a)=X(T,a)+O(\delta^2)$: the
corrected flow carries each particle to its original final position,
up to an error of order $\delta^2$. When $b$ solves the linearized
Euler equations \eqref{cc:intro:linearized-momentum} below, it is
determined by $b(0)$ and the source. Under the bounds $H\le\Id$ and
$T\le1$ below, the endpoint condition selects $b(0)$ uniquely.

For a prescribed momentum source $f$, the inhomogeneous problem is to
find $b,q$ satisfying the linearized Euler equations
\begin{equation}\label{cc:intro:linearized-momentum}
 \partial_t b+(u\cdot\nabla)b+(b\cdot\nabla)u+\nabla q=f,
 \qquad \diver b=0.
\end{equation}
Here $q$ is the pressure variation, and the displacement induced by $b$
must satisfy $\xi(0)=\xi(T)=0$. Differentiating
\eqref{cc:intro:first-variation} in time and using
\eqref{cc:intro:linearized-momentum} gives
\begin{equation}\label{cc:intro:jacobi}
 \begin{gathered}
 \xi_{tt}+H\xi+(\nabla q)\circ X=f\circ X,
 \qquad H=(D_x^2p)\circ X,\\
 \diver_x(\xi\circ X^{-1})=0,\qquad \xi(0)=\xi(T)=0.
 \end{gathered}
\end{equation}
The homogeneous equation is the Jacobi equation along the geodesic
associated with this solution in Arnold's formulation~\cite{Arnold1966}.
The constraint $\diver_x(\xi\circ X^{-1})=0$ expresses volume
preservation to first order.

For a solution of \eqref{cc:intro:jacobi}, testing with $-\xi$ and
integrating by parts gives the identity
\begin{equation}\label{cc:intro:Q-identity}
 \underbrace{\int_0^T\!\int_{\mathbb T^3}
   \bigl(|\xi_t|^2-\xi\cdot H\xi\bigr)\,\dd a\,\dd t}_{\mathcal Q(\xi)}
 =-\int_0^T\!\int_{\mathbb T^3}(f\circ X)\cdot\xi\,\dd a\,\dd t.
\end{equation}
Volume preservation and the divergence constraint remove the pressure
variation. The quadratic form $\mathcal Q$ is the fixed-endpoint second variation of the action.
If $H\le\Id$ and $T\le1$, then
${|\xi(t,a)|^2\le(T-t)\int_t^T|\xi_\tau(\tau,a)|^2\,\dd\tau}$ gives
\begin{equation}\label{cc:intro:coercivity}
 \mathcal Q(\xi)\ge\Bigl(1-\frac{T^2}{2}\Bigr)\|\xi_t\|_{L^2_{t,a}}^2
 \ge\frac12\|\xi_t\|_{L^2_{t,a}}^2 .
\end{equation}
The term $-\xi\cdot H\xi$ lowers the integrand where $\xi\cdot H\xi>0$
and raises it where $\xi\cdot H\xi<0$. Thus an upper bound for $H$ gives
the lower bound in \eqref{cc:intro:coercivity}.
For $f=0$, \eqref{cc:intro:Q-identity} gives $\mathcal Q(\xi)=0$. Then
\eqref{cc:intro:coercivity} forces $\xi_t=0$, and $\xi(0)=0$ gives $\xi=0$.
Equation~\eqref{cc:intro:first-variation} gives $b=0$, proving that the linearized
endpoint map $b(0)\mapsto\xi(T)$ has trivial kernel.
Upper bounds on the pressure Hessian also enter Brenier's fixed-endpoint
action-minimization criterion~\cite{Brenier1989,Brenier1999}.

For an approximate parent $U$ with residual $R_U$,
Subsection~\ref{cc:uf:sec:parent-step} uses
$H^\#=(D_x^2p_U-\nabla R_U)\circ X$ in place of $H$, under the same
bounds $H^\#\le\Id$ and $T\le1$. Subsection~\ref{cc:uf:sec:mean} uses
the resulting coercivity to construct a bounded inverse for the forced
Jacobi problem, and
Appendix~\ref{cc:app:mean} establishes the quantitative H\"older bounds
of Theorem~\ref{cc:uf:mean:inverse}. The inverse exists on every
interval $[0,T]$ with $0<T\le1$ on which $H^\#\le\Id$, with a bound that
depends on the parent but requires no further smallness of $T$. Each
correction therefore solves a linear two-point boundary value problem
along the flow of the parent, and no short-time argument is used.

In the affine wave of Subsection~\ref{cc:sec:intro-waves}, the Hessian of
$-p_{\rm w}P(k\cdot x)$ is $-p_{\rm w}F'(k\cdot x)k\otimes k$.
We seek an odd periodic profile
$F_h$ satisfying
\begin{equation}\label{cc:intro:profile-slopes}
 F_h'(0)=1,\qquad -Ch\le F_h'(s)\le1 .
\end{equation}
The first condition preserves the wave gradient $\mathbf a\otimes k$
at the origin, and hence the amplification described in Subsection~\ref{cc:sec:intro-waves}.
The lower bound on $F_h'$ controls the largest eigenvalue of
$D_x^2[-p_{\rm w}P_h(k\cdot x)]$, where $P_h'=F_h$. On an interval
with $p_{\rm w}>0$, its positive part satisfies
\begin{equation}\label{cc:intro:profile-curvature}
 p_{\rm w}|k|^2\bigl(-F_h'(k\cdot x)\bigr)_+
 \le Ch\,p_{\rm w}|k|^2,
\end{equation}
where $(r)_+=\max\{r,0\}$. For fixed wave coefficients, taking $h$ small
makes this positive part small at every point. In the construction $h$
is a negative power of the frequency of the wave. As $h$ shrinks, the
wave produces thin layers of strong shear while its velocity amplitude
is $O(h|\mathbf a|)$. Its gradient at the origin remains
$\mathbf a\otimes k$, preserving the local shear that amplifies the
next wave. The pressure Hessian added at each step is bounded above
by a small multiple of the identity. Section~\ref{cc:uf:sec:tower}
sums these bounds and controls $\|\nabla R_U\|_\infty$ to obtain
$H^\#\le\Id$ as a quadratic form.

\begin{figure}[!hbp]
\centering
\includegraphics[width=.75\textwidth]{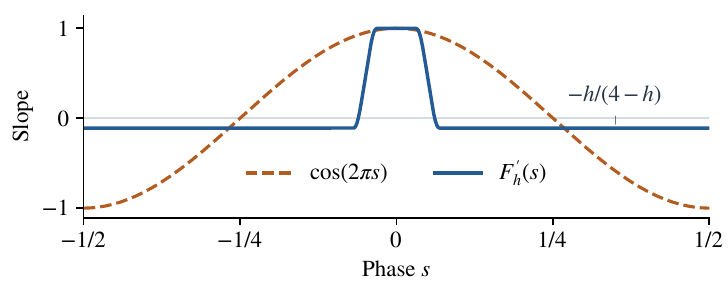}
\caption{The sine slope divided by its value at the origin (dashed)
and the slope $F_h'$ of the modified profile (solid). Both curves equal
one at the origin.}
\label{cc:fig:pressure-profile}
\end{figure}

For $F(s)=\sin(2\pi s)$, the ratio $F'(1/2)/F'(0)$ is $-1$.
The sine profile with slope one at the origin therefore violates the
lower bound in \eqref{cc:intro:profile-slopes} for small $h$.
Subsection~\ref{cc:sec:profile} constructs $F_h$ with a narrow positive
bump in $F_h'$ and a shallow negative part, as shown in
Figure~\ref{cc:fig:pressure-profile}.

\subsection{Correcting a localized wave}
We smoothly truncate the velocity of each step in frequency at the end of the step,
to simplify the bookkeeping of derivative estimates, as in the mollification
step of convex integration~\cite[Section~2.4]{BDLSVOnsager}.
This introduces a small error, carried to the next step;
Subsection~\ref{cc:uf:trunc:sec:truncation} gives the estimates.
Write $U$ for the truncated velocity of the preceding step and $R_U$ for
its residual in the Euler equations.

We add the localized wave described in Subsection~\ref{cc:sec:intro-waves},
using $F_h$ and a divergence correction. Let $X$ be the flow of $U$ from
time zero, and fix a nonzero vector $k_*$. The phase $\theta$ of the
localized wave is defined by
\[
 \theta(t,X(t,a))=k_*\cdot a
\]
for Lagrangian coordinates $a$ in the fundamental cell. It is the plane
phase $k_*\cdot x$ at time zero, transported by $U$. Its gradient
$\kappa=\nabla_x\theta$ is the wave vector of the localized wave; for a
linear flow, $\theta=\kappa\cdot x$, as in
Subsection~\ref{cc:sec:intro-waves}. For a wave after the first, $k_*$
is chosen so that $\kappa=\lambda E_3$ at the origin at its activation
time, as for the new wave of Subsection~\ref{cc:sec:intro-waves}.
The leading velocity is $A_1(t,x,\theta(t,x))$. Writing the profile as
$A_1(t,x,s)$ with $s$ independent lets us distinguish derivatives of
the coefficients from the large term $\kappa\partial_s$ that arises
when differentiating through the phase.
A wave after the first is multiplied
by a cutoff $\widetilde\chi(t)$ in time that vanishes before its
amplification interval; for the first wave $\widetilde\chi=1$.
The derivative of the cutoff leaves an error of the size of the wave at
the start of that interval, also carried to the next step.

The velocity added to $U$ is a finite expansion
\[
 \sum_{p=1}^PZ_p+C_{P+1},\qquad Z_p=V_p+C_p+B_p,
\]
where $V_p(t,x,s)$ is a profile perpendicular to $\kappa$ with zero
phase average, $C_{p+1}$ is the divergence correction of $V_p$, and
$B_p(t,x)$ is a mean, independent of the phase; the profiles are
evaluated at $s=\theta$. The leading term is the wave,
$Z_1=V_1=\widetilde\chi A_1$, with $C_1=0$ and $B_1=0$. We call $p$ the
\emph{level} of $Z_p$ and of each of $V_p$, $C_p$ and $B_p$: up to a
fixed power of a number that bounds the parent, the bound for $Z_p$ is
that of the wave times $\lambda^{-2(p-1)/3}$, where $\lambda$ is the
frequency of the wave. The corrections are computed one level at a
time up to the depth $P$, which is fixed in
Section~\ref{cc:uf:sec:tower}.
The error left by the wave and its divergence correction $C_2$ is a
profile $H_2(t,x,s)$ of level two.
We split $H_2$ into its phase average $\av{H_2}$ and its
oscillatory part $H_2-\av{H_2}$, as in the expansions of
Cheverry~\cite[Section~4.1]{Cheverry2006}.

The mean correction $B_2$ cancels the phase average of $H_2$ together
with the parent residual. With $D_t=\partial_t+U\cdot\nabla_x$, the
linearized Euler equation \eqref{cc:intro:linearized-momentum} gives
$D_t B_2+(\nabla U)B_2+\nabla q=-\av{H_2}-R_U$ with
$\diver B_2=0$. The Jacobi inverse selects $B_2(0)$ so that the
correction leaves the truncated parent's final particle positions fixed
to first order.

With $B_2$ determined, the oscillatory correction $V_2$ and its pressure
$\pi_2$ solve
\begin{equation}\label{cc:intro:transverse-correction}
 \begin{split}
 D_tV_2+(\nabla U)V_2+\kappa\partial_s\pi_2
   &=-(H_2-\av{H_2})-(B_2\cdot\kappa)\partial_s(\widetilde\chi A_1),\\
 \kappa\cdot V_2&=0,\qquad \av{V_2}=0.
 \end{split}
\end{equation}
The mean $B_2$ transports the leading wave's fast phase. This term has
zero phase average and, like the oscillatory part $H_2-\av{H_2}$,
belongs to level two, so we solve for the mean $B_2$ before the
oscillatory correction $V_2$. At fixed phase along each trajectory,
\eqref{cc:intro:transverse-correction} is a forced Craik--Criminale
amplitude equation. Every wave uses
zero initial data for this correction. The curl construction of
Subsection~\ref{cc:sec:material} gives the divergence correction $C_3$
of $V_2$, so that $V_2+C_3$ is exactly divergence free. The term $C_3$
belongs to level three, where we solve for the mean and oscillatory
corrections again.

\begin{figure}[tbp]
\centering
\begin{tikzpicture}[
  font=\small,
  block/.style={draw=black!65, rounded corners=2pt, align=center,
                inner xsep=5pt, inner ysep=4pt},
  correction/.style={block, text width=37mm, minimum height=8mm},
  line/.style={-{Stealth[length=1.7mm]}, semithick, draw=black!75}
]
\node[block] (parent) at (0,0) {Parent with residual $R_U$};
\node[block] (wave) at (0,-1.2)
  {Add localized wave and divergence-free corrector (\S\ref{cc:sec:material})\\
   \emph{Amplify vorticity at the origin}};
\node[inner sep=0pt] (linear-title) at (0,-2.55)
  {Linear corrections (\S\ref{cc:uf:sec:hierarchy})};
\node[correction] (mean) at (-4.65,-3.20)
  {Mean (\S\ref{cc:uf:sec:mean})};
\node[correction] (transverse) at (0,-3.20)
  {Oscillatory};
\node[correction] (curl) at (4.65,-3.20)
  {Divergence-free};
\node[draw=black!40, dashed, rounded corners=3pt,
      fit=(linear-title)(mean)(transverse)(curl), inner sep=4pt] (linear) {};
\node[block] (completion) at (0,-4.50)
  {Truncate in frequency (\S\ref{cc:uf:trunc:sec:truncation})};
\node[block] (exact) at (0,-5.30)
  {Next parent (\S\ref{cc:uf:trunc:sec:truncation})};
\draw[line] (parent) -- (wave);
\draw[line] (wave) -- (linear);
\draw[line] (mean) -- (transverse);
\draw[line] (transverse) -- (curl);
\draw[line] (linear) -- node[right, xshift=5pt]
  {Small remaining error} (completion);
\draw[line] (completion) -- (exact);
\end{tikzpicture}
\caption{One step of the construction.}
\label{cc:fig:correction-flow}
\end{figure}
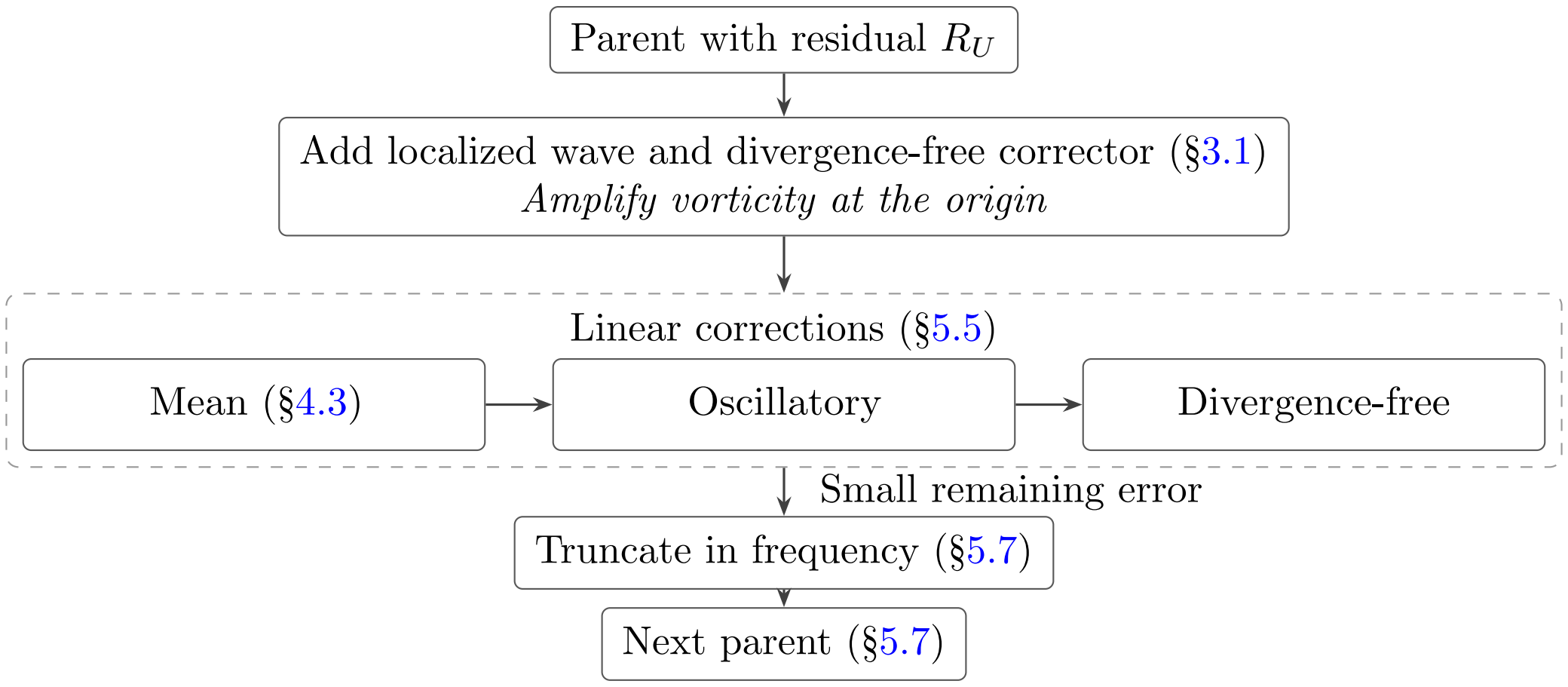

We stop the expansion at the depth $P$, keeping the divergence
correction $C_{P+1}$ of $V_P$. The remainder then has a small factor
$\lambda^{-2P/3}$, once $\lambda$ is larger
than a fixed power of a number that bounds the parent
(Theorem~\ref{cc:uf:hier:thm:depth}). The truncated parent, localized wave
and linear corrections give an approximate solution $V$ of the Euler
equations. We truncate $V$ in frequency. Three choices make the three
errors of the step smaller than a fixed negative power of the frequency
of the next wave. The depth $P$ controls the remainder of the
expansion. The frequency cutoff controls the truncation error. The
ratio $\varepsilon$ of Subsection~\ref{cc:sec:intro-waves} sets the
gain, which controls the error of the cutoff in time. This cutoff acts
while the wave's amplitude is at most a fixed multiple of
$\exp(-\varepsilon^{-1/2})$ times its peak. These three errors make up
the residual of the truncated velocity, which is the parent of the next step
(Proposition~\ref{cc:uf:step:thm:step}).
The pressure and residual-gradient estimates give coercivity for
the mean equation of the next step. The leading wave and oscillatory
corrections stay in the transported cutoff region. The mean corrections
and the truncation difference may extend over the torus. Figure~\ref{cc:fig:correction-flow} summarizes the step.

\subsection{Comparison with blow-up constructions with forcing}
The construction of blow-up for the forced Boussinesq equations by
Alp\"oge and Buckmaster \cite[Section~1.2]{AlpogeBuckmasterBoussinesq}
also uses coupled wave
amplitudes to amplify a gradient and prepare the next wave. There the
temperature gradient grows together with a vorticity amplitude. The
accompanying shear would obstruct the next growth stage, so an external
force rotates the temperature gradient and wave direction together.
This preserves the coupling that drives the temperature amplitude while
changing the horizontal component of the wave vector, and hence the sign
of vorticity production. A scalar shooting parameter selects a rotation
pulse that returns the vorticity amplitude and that horizontal component
to zero, retaining the temperature gain. Keeping the wave vector vertical
then holds both amplitudes fixed while the next wave grows.

The construction of blow-up for the forced Euler equations by
Alp\"oge and Buckmaster \cite[Section~1.2]{AlpogeBuckmasterForcedEuler}
builds on this mechanism for axisymmetric flow with swirl. Circulation replaces temperature and
couples to azimuthal vorticity. The angular control returns the principal
azimuthal-vorticity amplitude at the material center to zero while
retaining the circulation-gradient gain. The radial dependence of the
elliptic operator and of this coupling requires amplitudes that vary
across the transported support, with transport by the preceding velocity,
including its corrections.

In the present unforced construction, the amplified shear itself drives the successor. Its orientation
is obtained from the Craik--Criminale wave evolution; the pressure coefficient
is determined by incompressibility and is not a rotation control.
The next wave is activated at the start of its amplification interval,
with its amplitude in the growing direction (Subsection~\ref{cc:uf:sec:later-wave}).

Localization presents a further difference. In the construction for
the Boussinesq equations the profiles are linear near zero, so each layer's velocity and temperature are exactly linear
on a region containing the later layers. The streamfunction corrections vanish there, preserving
this region. In both forced constructions, the remaining errors are
absorbed by a force smooth through the singular time. To solve the unforced Euler equations,
we must remove those errors. The pressure and velocity
corrections need not preserve a neighborhood on which the velocity is linear;
the next wave uses the central gradient and estimates on its spatial variation.

The organization of this paper is as follows:
Section~\ref{cc:sec:waves} proves amplification and its stability,
Sections~\ref{cc:sec:localized}--\ref{cc:sec:exact} localize the waves
and construct the linear corrections,
Sections~\ref{cc:sec:iteration} and~\ref{cc:uf:sec:tower} start and
iterate the construction, and Section~\ref{cc:sec:limit} passes to the limit.

The public repository
\url{https://github.com/mcoiculescu/EulerBlowupPublic} contains a Lean
formalization of the proof and the latest version of this manuscript.

\subsection*{Conventions}
A time supremum is indicated by $C_t$ or by the time interval. The
Leray projection is
$\PP=I-\nabla\Delta^{-1}\diver$, where $\Delta^{-1}$ acts on zero-mean
scalars. The periodic velocities
constructed below are odd and hence have zero spatial mean.
Pressure is normalized to have zero spatial mean. Spatial differences use the periodic
metric, and flows are lifted so that their displacements are periodic.
We use $(\nabla U)_{ij}=\partial_jU_i$, and for a matrix field $A$,
$(\diver A)_i=\sum_j\partial_jA_{ij}$. A matrix acts on vectors by
$(Av)_i=\sum_jA_{ij}v_j$. For vectors $a,b$ we write
$a\otimes b=ab^T$, so that $(a\otimes b)_{ij}=a_ib_j$ and
$(a\otimes b)v=(b\cdot v)a$. For an orthonormal basis $(E_1,E_2,E_3)$,
the entries of $A$ in this basis are $A_{ij}=E_i\cdot AE_j$; then
$AE_j=\sum_iA_{ij}E_i$, $A^TE_i=\sum_jA_{ij}E_j$ and
$Av=\sum_{i,j}A_{ij}(E_j\cdot v)E_i$.

\subsection{Choice of parameters}\label{cc:sec:parameters}
The frequencies and orders are chosen so that the truncation and correction
errors are summable. The gains also make each new wave small at its
activation time. We fix the scales here; Appendix~\ref{cc:app:parameters}
proves Lemma~\ref{cc:uf:tower:lem:slow-orders}.

Fix the smooth profiles used in the construction.
The integer $b\ge1$ is the fixed integer specified in
Subsection~\ref{cc:uf:sec:approximate}.
For integers $s\ge1$ put $K(s)=320b(s+1)$.
The exponent $K$ separates successive frequencies.

The estimates of Sections~\ref{cc:sec:localized}, \ref{cc:sec:exact},
\ref{cc:sec:iteration} and~\ref{cc:uf:sec:tower} and of
Lemma~\ref{cc:uf:trunc:lem:cutoff} contain constants that depend only
on an order of differentiation, the depth of the expansion and the
fixed profiles. We call them the fixed-order constants. Choose a
nondecreasing function $\mathcal B$ on the positive integers such that
$\mathcal B(s)$ bounds the fixed-order constants at all orders and
depths at most $7K(s)$, and
\begin{equation}\label{cc:par:floor}
 \log\mathcal B(s)\ge4K(s+1)^3\qquad(s\ge1).
\end{equation}

Let $N\ge\mathcal B(1)$ be the initial frequency parameter. The integer
$q_j$ specifies the spatial derivative order through which step $j$
makes only a small change to the velocity on earlier time intervals.
We define the frequencies and $q_j$ together, starting from $q_1=1$.
Put
\begin{equation}\label{cc:uf:tower:eq:schedule}
 \begin{gathered}
 K_j=K(q_j),\qquad
 \lambda_0=N,\qquad \lambda_j=\lambda_{j-1}^{K_j}\quad(j\ge1),\qquad
 \Lambda_j=\lambda_j^5\quad(j\ge0),\\
 q_{j+1}=\max\{1\le s\le q_j+1:\ \mathcal B(s)\le\lambda_{j-1}\}\quad(j\ge1),\\
 L_j=4\log\lambda_{j+1},\qquad
 \varepsilon_{j-1}=L_j^{-2}
 \quad(j\ge2).
 \end{gathered}
\end{equation}
Here $\lambda_j$, $j\ge1$, is the frequency of wave $j$, and $L_j$ is the
prescribed logarithmic gain. For $j\ge1$, $\Lambda_j$ is the Fourier cutoff
applied to the velocity constructed at step $j$, and $\Lambda_0$ bounds
the frequencies of the shear from which the construction starts.

The orders $q_j$ must tend to infinity for the limit to be smooth;
however fast $\mathcal B$ grows, it only slows their growth. They must
also grow slowly: the amplification estimates of
Section~\ref{cc:uf:sec:tower} lose fixed powers of the logarithmic gain
$L_{j+1}$, which grows with the orders through the exponents $K_j$, and
these losses must stay below small powers of $\lambda_{j-1}$. The lower
bound \eqref{cc:par:floor} is imposed for this purpose.

\begin{lemma}\label{cc:uf:tower:lem:slow-orders}
The integers $q_j$ satisfy $q_j\le q_{j+1}\le q_j+1$ and tend to
infinity, and
\begin{equation}\label{cc:uf:tower:eq:uniform}
 L_{j+1}\le(\log\lambda_{j-1})^2\qquad(j\ge1).
\end{equation}
\end{lemma}

\subsection{Provenance}\label{cc:sec:provenance}

Before OpenAI released its paper on the unforced Euler equations, the first
and second authors had completed a second proof of blow-up for the Euler equations
with smooth forcing. Unlike their public
proof~\cite{AlpogeBuckmasterForcedEuler}, which uses axisymmetric flow
with swirl, this unreleased proof uses the same Craik--Criminale
amplification mechanism as OpenAI's paper, including the same reduced
equations for the amplitude of a wave, the system
\eqref{cc:eq:model-system} below. In our
unpublished work on the hypodissipative Navier--Stokes equations, much of the argument
for removing the force was developed, with human prompts guiding that
work. This is consistent with the second author's statement of
September~7 that the hypodissipative work pointed toward a blow-up
construction for the unforced Euler equations~\cite{BuckmasterSeptemberStatement}.

The second author's account of his communications with OpenAI about
these results is given in~\cite{BuckmasterSeptemberStatement}.

\subsection{AI statement}\label{cc:sec:ai}

We used Claude and Codex throughout this work. The first proof without
force was written by Anthropic models from the drafts on the Euler equations
with forcing and on the hypodissipative Navier--Stokes equations described in
Subsection~\ref{cc:sec:provenance}. Like the first writeups of our
proofs for the Boussinesq and Euler equations, it was long and hard to read. We then
rewrote it with Claude and Codex under our direction, iterating with
the models on the structure of the argument, isolating the key ideas,
and simplifying it over successive revisions. In the course of this we
replaced parts of the argument by ideas of our own, which the models
implemented. Together with our own editing, this led to the present
text.

\section{The wave and its amplification}\label{cc:sec:waves}
\subsection{Exact waves over a linear background}\label{cc:sec:affine}
We first work on $\R^3$, where linear velocity fields give local models for the
periodic construction. Let $D(t)$ be smooth and trace free, with
$D'+D^2$ symmetric. Then $Dx$ solves the Euler equations with pressure
$-x\cdot(D'+D^2)x/2$. Choose a smooth odd $1$-periodic profile
$F$ with $F'(0)=1$, and let $P'=F$, $\av P=0$.
Here $\av H=\int_0^1H(s)\,\mathrm ds$ is the phase average.

\begin{lemma}\label{cc:lem:wave}
Let the wave vector $k(t)\ne0$ and the amplitude vector $\mathbf a(t)$,
functions of time alone, satisfy
\begin{equation}\label{cc:eq:wave-ode}
 k'=-D^Tk,\qquad
 \mathbf a'=-D\mathbf a+p_{\rm w} k,\qquad
 p_{\rm w}=\frac{2k\cdot D\mathbf a}{|k|^2},\qquad
 k\cdot\mathbf a=0.
\end{equation}
Then $Dx+\mathbf aF(k\cdot x)$ is divergence free and solves the Euler equations on $\R^3$ with
pressure $-x\cdot(D'+D^2)x/2-p_{\rm w} P(k\cdot x)$.
If $k(t_0)\ne0$, $\mathbf a(t_0)\ne0$ and
$k(t_0)\cdot\mathbf a(t_0)=0$, then both vectors remain nonzero and
perpendicular for every $t\ge t_0$ in the interval of existence of $D$.
\end{lemma}
\begin{proof}
Write $w=\mathbf aF(k\cdot x)$ and $u=Dx+w$. Since $k'=-D^Tk$,
\[
 (\partial_t+Dx\cdot\nabla)(k\cdot x)=(k'+D^Tk)\cdot x=0,
 \qquad (\partial_t+Dx\cdot\nabla)w=\mathbf a'F.
\]
Also $(w\cdot\nabla)Dx=D\mathbf aF$, and perpendicularity cancels
the self-advection,
\[
 (w\cdot\nabla)w=(k\cdot\mathbf a)FF'\mathbf a=0.
\]
Consequently
\[
 \partial_tu+(u\cdot\nabla)u=(D'+D^2)x+(\mathbf a'+D\mathbf a)F.
\]
The first term is minus the gradient of the background pressure, and
the second equals $p_{\rm w}kF$, which is minus the gradient of the
wave pressure $-p_{\rm w}P(k\cdot x)$. Moreover
$\diver u=\tr D+(k\cdot\mathbf a)F'=0$.
The choice of $p_{\rm w}$ gives
$(k\cdot\mathbf a)'=-2k\cdot D\mathbf a+p_{\rm w}|k|^2=0$.
For a given $D$, the equation for $k$ is linear and has an invertible
fundamental matrix. Hence $k(t_0)\ne0$ implies $k(t)\ne0$. With this
$k$, the equation
$\mathbf a'=-(\Id-2\frac{k}{|k|}\frac{k^T}{|k|})D\mathbf a$
is also linear with an invertible fundamental matrix, so
$\mathbf a(t_0)\ne0$ implies $\mathbf a(t)\ne0$.
\end{proof}

For a single profile, this is the inviscid wave construction of
Craik and Criminale~\cite{CraikCriminale}.

\subsubsection{The preceding frame}\label{cc:sec:frame}
A wave evolving over a background matrix $B$ becomes part of the
background for the next wave. At the origin, its contribution to the
velocity gradient is $\mathbf a\otimes k$, and its unit wave vector and
velocity direction provide a moving orthonormal frame. We express the
older gradient $B$ in this frame to identify the couplings that drive the
next wave. The following calculation uses only $\tr B=0$; it requires no
symmetry of $B'+B^2$.

\begin{lemma}\label{cc:lem:frame}
Let $B(t)$ be a smooth trace-free matrix, and let $k\ne0$ and
$\mathbf a\ne0$ satisfy \eqref{cc:eq:wave-ode} with $D$ replaced by $B$.
Put $E_1=\frac{k}{|k|}$, $E_2=\frac{\mathbf a}{|\mathbf a|}$ and $E_3=E_1\times E_2$,
write $B_{ij}=E_i\cdot BE_j$ $(i,j=1,2,3)$ for the entries of $B$ in this
orthonormal frame, and let $\sigma=|k|\,|\mathbf a|$, $\nu=B_{12}$ and
$\eta=B_{32}$. Then
\begin{equation}\label{cc:eq:frame}
 E_1'=-\nu E_2-B_{13}E_3,\quad
 E_2'=\nu E_1-\eta E_3,\quad
 E_3'=B_{13}E_1+\eta E_2,
 \qquad(\log\sigma)'=B_{33}.
\end{equation}
\end{lemma}
\begin{proof}
For an arbitrary nonzero vector $v(t)$ depending on time, put
$\widehat v=v/|v|$. In this proof primes denote derivatives in $t$.
Then
\[
 \widehat v'=(\Id-\widehat v\widehat v^T)\frac{v'}{|v|},\qquad
 (\log|v|)'=\widehat v\cdot\frac{v'}{|v|}.
\]
Since $B^TE_1=\sum_jB_{1j}E_j$ and $k'=-B^Tk$, this gives
\[
 E_1'=-(\Id-E_1E_1^T)B^TE_1=-B_{12}E_2-B_{13}E_3,\qquad
 (\log|k|)'=-B_{11}.
\]
Similarly $BE_2=\sum_iB_{i2}E_i$ and
$p_{\rm w}|k|=2|\mathbf a|B_{12}$, so
\[
 E_2'=(\Id-E_2E_2^T)(-BE_2+2B_{12}E_1)=B_{12}E_1-B_{32}E_3,
 \qquad (\log|\mathbf a|)'=-B_{22}.
\]
Using $E_3\times E_2=-E_1$ and $E_1\times E_3=-E_2$ gives
\[
 E_3'=E_1'\times E_2+E_1\times E_2'=B_{13}E_1+B_{32}E_2.
\]
Finally,
$(\log\sigma)'=-B_{11}-B_{22}=B_{33}$ because $\tr B=0$.
\end{proof}

\subsection{Amplification on a prescribed interval}\label{cc:sec:amplification}

Of the two entries $\nu=E_1\cdot BE_2$ and $\eta=E_3\cdot BE_2$ of the
older gradient introduced in Lemma~\ref{cc:lem:frame}, the coupling
$\nu$ helps the preceding shear amplify a new wave, while $\eta$
changes its direction and eventually ends the growth. A small ratio
$\eta/\nu$ therefore allows a large gain before decay. In the model
case of Subsection~\ref{cc:sec:intro-waves} and the limit of an
infinitely strong preceding shear, the amplitude of a new wave solves
a linear system of two equations with explicit coefficients
(Lemma~\ref{cc:lem:strong-shear}). Up to errors measured by one number
$\delta$, the amplitude of a general new wave solves the same system
(Lemma~\ref{cc:uf:exact:lem:additive}). Its solution with growing
initial data peaks after a large gain (Lemma~\ref{cc:lem:model}), and
the wave follows this solution with relative errors that grow only
polynomially in the length of the interval
(Lemma~\ref{cc:lem:finite-comparison}). The couplings of the new wave
are estimated in Lemma~\ref{cc:lem:new-wave-couplings}, and
Proposition~\ref{cc:prop:finite-return} locates the ratio prescribed
for its successor after the amplitude peaks in the strong-shear limit.

\subsubsection{The new wave in the preceding frame}\label{cc:sec:comparison}
Let $B(t)$ be a smooth trace-free matrix, and let the preceding wave
$(k_p,\mathbf a_p)$, with $k_p\ne0$ and $\mathbf a_p\ne0$, satisfy
\eqref{cc:eq:wave-ode} with $D$ replaced by $B$. Let $(E_1,E_2,E_3)$ be
its frame, $B_{ij}=E_i\cdot BE_j$ the entries of $B$ in this frame, and
$\sigma_p=|k_p|\,|\mathbf a_p|$ its strength, as in
Lemma~\ref{cc:lem:frame}. Thus \eqref{cc:eq:frame} holds with the
functions $B_{12}$, $B_{32}$ and $\sigma_p$ in the places of
$\nu$, $\eta$ and $\sigma$. The new wave evolves over
\[
 M=B+\mathbf a_p\otimes k_p+E_c=B+\sigma_p E_2\otimes E_1+E_c,
\]
the sum of the older gradient, the gradient of the preceding wave at
the origin, and a continuous matrix $E_c(t)$ of errors. Its wave vector
$k$ solves the first equation of \eqref{cc:eq:wave-ode} with $D=M$ and
equals $\lambda E_3$, with $\lambda>0$, at the activation time $t_0$.
For this $k$, an \emph{amplitude vector} is a solution $\mathbf a$ of
the second equation in \eqref{cc:eq:wave-ode} with $D=M$ satisfying
$\mathbf a(t)\cdot k(t)=0$. Only these ordinary differential equations
enter this subsection: the matrices $B$ and $E_c$ are prescribed, and no
equation is assumed for them. Subsection~\ref{cc:uf:sec:later-wave}
applies the results below with $M$ the velocity gradient of the parent
at the origin, which is a fixed point of the flow.

Freeze the strength and the two couplings at the activation time,
\[
 \sigma=\sigma_p(t_0),\qquad \nu=B_{12}(t_0)>0,\qquad
 \eta=B_{32}(t_0)>0,
\]
and put
\[
 r=\sqrt{2\sigma\nu},\qquad S=\frac\sigma r=\sqrt{\frac{\sigma}{2\nu}},\qquad
 \varepsilon=\frac\eta\nu,\qquad \tau=r(t-t_0).
\]
Since $\nu/\sigma=1/(2S^2)$, the growth rate
$\sqrt{2\nu(\sigma-\nu)}$ found in Subsection~\ref{cc:sec:intro-waves}
equals $r(1-1/(2S^2))^{1/2}$; thus $r$ is this rate up to a factor
close to one when $S$ is large. The variable
$\tau$ measures time in units of the growth time $r^{-1}$, and $S$ is
the ratio of the preceding shear rate to the growth rate. We work on the
interval $0\le\tau\le T$, or $t_0\le t\le t_0+T/r$, with $T\ge2$.
In \eqref{cc:eq:wave-ode} and \eqref{cc:eq:frame} primes denote
derivatives in $t$. From now on in this subsection they denote
derivatives in $\tau$, so that $d/d\tau=r^{-1}d/dt$. All suprema are
over $0\le\tau\le T$, and $b_{ij}=B_{ij}/r$.

Define the scaled components $\zeta=(\zeta_1,\zeta_2,\zeta_3)$ and
$\alpha=(\alpha_1,\alpha_2,\alpha_3)$ of the wave vector and an amplitude
vector by
\begin{equation}\label{cc:eq:scaled-pair}
 k=\lambda\left(\zeta_1E_1+\frac{\zeta_2}{S}E_2+\zeta_3E_3\right),
 \qquad
 \mathbf a=\frac{\alpha_1}{S}E_1+\alpha_2E_2+\frac{\alpha_3}{S}E_3.
\end{equation}
Then $\zeta(0)=(0,0,1)$ and
\[
 |k|^2=\lambda^2m,\qquad
 m=\zeta_1^2+\zeta_2^2/S^2+\zeta_3^2.
\]
Transversality is $\zeta\cdot\alpha=0$, so where $\zeta_3\ne0$,
\[
 \alpha_3=-\frac{\zeta_1\alpha_1+\zeta_2\alpha_2}{\zeta_3}.
\]
The two components $\alpha_1,\alpha_2$ therefore determine the amplitude
vector. The factor $S^{-1}$ in the $E_2$ coordinate of $k$ makes the
contribution of $B$ to the coefficient of $\zeta_3$ in the equation for
$\zeta_2$ equal to $-B_{32}/\nu$, which equals $-\varepsilon$ at
$\tau=0$. In the growing
mode of Subsection~\ref{cc:sec:intro-waves}, the components of the
amplitude along $E_1$ and $E_2$ have opposite signs, and the first is
$S^{-1}(1-\nu/\sigma)^{-1/2}$ times the second in absolute value.
With the scaling in \eqref{cc:eq:scaled-pair}, this mode is
$\alpha_1=-\alpha_2$ to leading order. In
Subsection~\ref{cc:sec:scalar} the components $(\alpha_1,\alpha_2)$ are
divided further by the amplitude in the strong-shear limit, which gives
the normalized coordinates \eqref{cc:eq:normalized}.

In the model case of Subsection~\ref{cc:sec:intro-waves}, the entries
of $B$ in the frame are the constants $B_{12}=\nu$ and $B_{32}=\eta$,
all other entries vanish, and $E_c=0$. Put
\begin{equation}\label{cc:eq:reference-angle}
 \tan\varphi=\tfrac12\varepsilon\tau^2,\qquad
 0\le\varphi<\tfrac\pi2,\qquad
 \zeta_0=(\tan\varphi,-\varepsilon\tau,1).
\end{equation}
The next lemma derives the equations for the scaled components in the
model case and their limit as the preceding shear becomes infinitely
strong. In the limit the wave vector is $\lambda(\tan\varphi E_1+E_3)$,
so $\varphi$ is its angle from $E_3$ toward $E_1$, as in
Subsection~\ref{cc:sec:intro-waves}.

\begin{lemma}\label{cc:lem:strong-shear}
Suppose that on $0\le\tau\le T$ the entries of $B$ in the frame are
$B_{12}=\nu$, $B_{32}=\eta$ and $B_{ij}=0$ otherwise, and that $E_c=0$.
Then $\sigma_p=\sigma$, and the scaled components of the wave vector
satisfy
\[
 \zeta_1'=-\Bigl(1+\frac1{2S^2}\Bigr)\zeta_2,\qquad
 \zeta_2'=-\varepsilon\zeta_3,\qquad
 \zeta_3'=\frac{\varepsilon}{2S^2}\zeta_2 .
\]
Where $\zeta_3\ne0$, the components $\alpha_1,\alpha_2$ of every
amplitude vector satisfy
\[
 \begin{aligned}
 \alpha_1'&=-\alpha_2+\zeta_1\Pi,\\
 \alpha_2'&=-\Bigl(1-\frac1{2S^2}\Bigr)\alpha_1
   +\frac{\varepsilon}{2S^2}\,\frac{\zeta_1\alpha_1+\zeta_2\alpha_2}{\zeta_3}
   +\frac{\zeta_2}{S^2}\,\Pi,
 \end{aligned}
 \qquad
 \Pi=\frac{2\zeta_2\alpha_1+(\zeta_1+\varepsilon\zeta_3)\alpha_2}m .
\]
As $S\to\infty$ with $\varepsilon$ and $T$ fixed, $\zeta\to\zeta_0$
uniformly on $[0,T]$, and the coefficients of this system for
$(\alpha_1,\alpha_2)$ converge uniformly on $[0,T]$ to those of
\begin{equation}\label{cc:eq:model-system}
 \alpha_1'=-\cos^2\varphi(1-\varepsilon\tan\varphi)\alpha_2
                 -\varepsilon\tau\sin2\varphi\,\alpha_1,
 \qquad \alpha_2'=-\alpha_1 .
\end{equation}
\end{lemma}
\begin{proof}
In this proof write $k_i=k\cdot E_i$ and $a_i=\mathbf a\cdot E_i$, and
let a dot denote the derivative in $t$. The frame equations
\eqref{cc:eq:frame} read $\dot E_1=-\nu E_2$,
$\dot E_2=\nu E_1-\eta E_3$ and $\dot E_3=\eta E_2$, and their last
equation gives $d(\log\sigma_p)/dt=B_{33}=0$, so $\sigma_p=\sigma$. Thus
$M=\nu E_1\otimes E_2+\eta E_3\otimes E_2+\sigma E_2\otimes E_1$, and
$ME_1=\sigma E_2$, $ME_2=\nu E_1+\eta E_3$, $ME_3=0$. Since
$\dot k=-M^Tk$, we have $\dot k_i=-k\cdot ME_i+k\cdot\dot E_i$, hence
\[
 \dot k_1=-(\sigma+\nu)k_2,\qquad \dot k_2=-2\eta k_3,\qquad
 \dot k_3=\eta k_2 .
\]
Since $\dot{\mathbf a}=-M\mathbf a+p_{\rm w}k$ and
$M\mathbf a=\nu a_2E_1+\sigma a_1E_2+\eta a_2E_3$, also
$\dot a_i=-E_i\cdot M\mathbf a+p_{\rm w}k_i+\mathbf a\cdot\dot E_i$, hence
\[
 \dot a_1=-2\nu a_2+p_{\rm w}k_1,\qquad
 \dot a_2=-(\sigma-\nu)a_1-\eta a_3+p_{\rm w}k_2,\qquad
 p_{\rm w}=\frac{2(\nu k_1a_2+\sigma k_2a_1+\eta k_3a_2)}{|k|^2}.
\]
Substitute \eqref{cc:eq:scaled-pair}, divide by $r$ to pass to the time
$\tau$, and use
\[
 rS=\sigma,\qquad \frac{\nu}{\sigma}=\frac1{2S^2},\qquad
 \frac{\nu S}r=\frac12,\qquad \frac{\eta S}r=\frac\varepsilon2,\qquad
 \frac\eta\sigma=\frac{\varepsilon}{2S^2}.
\]
The equations for $k_i$ give those for $\zeta$, and
$\lambda Sp_{\rm w}/r=\Pi$. The equations for $a_1$ and $a_2$ then give
those for $\alpha_1=Sa_1$ and $\alpha_2=a_2$, after substituting
$\alpha_3=-(\zeta_1\alpha_1+\zeta_2\alpha_2)/\zeta_3$.

As $S\to\infty$, the coefficients of the system for $\zeta$ converge to
those of $\zeta_1'=-\zeta_2$, $\zeta_2'=-\varepsilon\zeta_3$,
$\zeta_3'=0$, whose solution with $\zeta(0)=(0,0,1)$ is $\zeta_0$. Both
systems are linear with bounded coefficients on $[0,T]$, so
$\zeta\to\zeta_0$ uniformly there by Gronwall's inequality. In
particular $\zeta_3\to1$ and $m\to1+\tan^2\varphi=\cos^{-2}\varphi$.
In the equation for $\alpha_1$ the coefficients of $\alpha_1$ and
$\alpha_2$ therefore converge to
\[
 \frac{2\zeta_1\zeta_2}m\to-2\varepsilon\tau\tan\varphi\cos^2\varphi
   =-\varepsilon\tau\sin2\varphi,\qquad
 -1+\frac{\zeta_1(\zeta_1+\varepsilon\zeta_3)}m
 \to-\cos^2\varphi(1-\varepsilon\tan\varphi),
\]
since $(\tan\varphi+\varepsilon)\tan\varphi\cos^2\varphi
=\sin^2\varphi+\varepsilon\sin\varphi\cos\varphi$. In the equation for
$\alpha_2$ every term except $-\alpha_1$ has a factor $S^{-2}$ and a
bounded coefficient, so the two coefficients converge to $-1$ and $0$.
\end{proof}

We call \eqref{cc:eq:model-system} the \emph{reference system}. The
vector $\zeta_0$ solves the limiting system $\zeta_1'=-\zeta_2$,
$\zeta_2'=-\varepsilon\zeta_3$, $\zeta_3'=0$ for the wave vector.

The errors of the comparison are measured by
\begin{equation}\label{cc:uf:exact:eq:additive-error}
 \delta=\left\|\frac{B_{12}}{\nu}-1\right\|_\infty
       +\left\|\frac{B_{32}}{\eta}-1\right\|_\infty
       +\frac{\|B\|_\infty}{r}+\frac{\|E_c\|_\infty}{\nu}.
\end{equation}
The first two terms are the relative variations of the two couplings on
the interval. The third compares the older gradient with the growth
rate; it is at least $\nu/r=1/(2S)$, as shown below, so it is small
only if the preceding shear is strong. The last compares the error
$E_c$ with the coupling $\nu$ rather than with $r$. The scaled
components multiply the $E_2$ component of $k$ and the $E_1$ and $E_3$
components of $\mathbf a$ by $S$, so in the time $\tau$ an entry of
$E_c$ enters their equations with a factor up to $S/r=1/(2\nu)$; see
the proof of Lemma~\ref{cc:uf:exact:lem:additive}.

Fix an absolute constant $C_0$, sufficiently large for the finitely
many estimates of this subsection, and assume
\begin{equation}\label{cc:eq:finite-smallness}
 C_0\varepsilon\le1,\qquad C_0T^{C_0}\delta\le\varepsilon .
\end{equation}
All other constants and exponents $C,c$ in this subsection are
absolute and do not depend on $C_0$; they may change from one estimate
to the next. The errors below have the form $CT^C\delta$, and by
\eqref{cc:eq:finite-smallness} such an error is at most $\varepsilon$
times any prescribed negative power of $T$ once $C_0$ is large.

Condition \eqref{cc:eq:finite-smallness} forces the preceding shear to
be strong. Indeed, $|B|\ge|E_1\cdot BE_2|$, and at $\tau=0$ this entry
is $B_{12}(t_0)=\nu$. Since $\nu/r=\sqrt{\nu/(2\sigma)}=1/(2S)$, the
third term of \eqref{cc:uf:exact:eq:additive-error} gives
$\delta\ge\nu/r=1/(2S)$. On the other hand,
\eqref{cc:eq:finite-smallness} gives
$\delta\le\varepsilon C_0^{-1}T^{-C_0}\le C_0^{-2}T^{-C_0}$. For
$C_0\ge2$ and $T\ge2$ it follows that
\begin{equation}\label{cc:eq:S-bounds}
 S\ge\frac1{2\delta}\ge\frac{C_0^2T^{C_0}}2\ge T^2,\qquad
 S^{-2}\le4\delta^2\le\delta,
\end{equation}
where the last inequality uses $\delta\le1/4$.

In the applications in Subsection~\ref{cc:uf:sec:later-wave} and
Section~\ref{cc:uf:sec:tower}, the parameter choice
\eqref{cc:uf:tower:eq:schedule} prescribes $\varepsilon$, and
$\varepsilon^{-1}$ and $T$ are bounded by powers of the logarithm of
the frequency of the preceding wave, while $\delta$ is at most a
negative power of that frequency times a power of its logarithm.
Item~2 of the proof of Theorem~\ref{cc:uf:tower:thm:finite} verifies
\eqref{cc:eq:finite-smallness} in this way.

For a general background and $E_c=0$, each coefficient in the proof of
the next lemma is a polynomial in $\zeta$, $1/\zeta_3$, $1/m$ and the quantities
\begin{equation}\label{cc:eq:coefficient-list}
 b_{11},\ b_{22},\ b_{33},\ b_{13},\ b_{31},\
 \frac{b_{21}}S,\ \frac{b_{12}}S,\ \frac{b_{23}}S,\ \frac{b_{32}}S,\
 2Sb_{12},\ 2Sb_{32},\ \frac{\sigma_p}\sigma,\ S^{-2}.
\end{equation}
The \emph{reference value} of such an expression is its value at
$\zeta=\zeta_0$, $2Sb_{12}=1$, $2Sb_{32}=\varepsilon$ and
$\sigma_p/\sigma=1$, with all other quantities in
\eqref{cc:eq:coefficient-list} equal to zero. By
Lemma~\ref{cc:lem:strong-shear}, it is the limit of the expression in
the model case as $S\to\infty$.

\begin{lemma}\label{cc:uf:exact:lem:additive}
Assume \eqref{cc:eq:finite-smallness}. The scaled components of the wave
vector satisfy
\begin{equation}\label{cc:eq:wave-vector-comparison}
 \sup|\zeta-\zeta_0|\le CT^C\delta\le\tfrac12,\qquad
 \sup\left|\frac{|k|^2\cos^2\varphi}{\lambda^2}-1\right|
       \le CT^C\delta\le\tfrac12.
\end{equation}
In particular
\begin{equation}\label{cc:eq:zeta-bounds}
 \tfrac12\le\zeta_3\le\tfrac32,\qquad |\zeta|\le3T^2 .
\end{equation}
The components $\alpha_1,\alpha_2$ of every amplitude vector satisfy
\[
 \frac{d}{d\tau}\binom{\alpha_1}{\alpha_2}
 =\left\{-\begin{pmatrix}
 \varepsilon\tau\sin2\varphi&\cos^2\varphi(1-\varepsilon\tan\varphi)\\
 1&0
 \end{pmatrix}+\mathcal E(\tau)\right\}\binom{\alpha_1}{\alpha_2},
 \qquad \sup|\mathcal E|\le CT^C\delta,
\]
where $\mathcal E$ is the same for all amplitude vectors.
\end{lemma}
\begin{proof}
The equations \eqref{cc:eq:wave-ode} and \eqref{cc:eq:frame} are in the
time $t$, and every coefficient in them is a rate of change in $t$.
In this proof primes denote derivatives in $\tau$. Since
$d/d\tau=r^{-1}d/dt$, each coefficient is divided by $r$,
which turns $B_{ij}$ into $b_{ij}$ and $\sigma_p$ into
$S\sigma_p/\sigma$. By \eqref{cc:uf:exact:eq:additive-error} and
$r/S=2\nu$,
\[
 |b_{ij}|\le\delta,\qquad
 |2Sb_{12}-1|=\left|\frac{B_{12}}{\nu}-1\right|\le\delta,\qquad
 |2Sb_{32}-\varepsilon|
 =\varepsilon\left|\frac{B_{32}}{\eta}-1\right|\le\delta .
\]
In the time $\tau$ the last equation of \eqref{cc:eq:frame} reads
$(\log\sigma_p)'=b_{33}$, and $\sigma_p=\sigma$ at $\tau=0$. Hence
\[
 \frac{\sigma_p(\tau)}{\sigma}
 =\exp\left(\int_0^\tau b_{33}(\tau')\,d\tau'\right),
\]
and $\sup|\sigma_p/\sigma-1|\le2T\delta$, because the integral is at
most $T\delta\le1$ in absolute value.

When $E_c=0$, the equation for $k$ and
\eqref{cc:eq:frame} give
\[
 \zeta_2'=-b_{22}\zeta_2-2Sb_{32}\zeta_3,\qquad
 \zeta_3'=-b_{33}\zeta_3-(b_{23}-b_{32})\zeta_2/S,
\]
\[
 \zeta_1'=-b_{11}\zeta_1
  -(\sigma_p/\sigma+b_{21}/S+b_{12}/S)\zeta_2
  -(b_{31}+b_{13})\zeta_3.
\]
In the model case these are the equations of
Lemma~\ref{cc:lem:strong-shear}. When $E_c\ne0$, the term $-E_c^Tk$ is added to the equation for $k$,
and the frame still evolves over $B$. In the time $\tau$, the
additional coefficient matrix for $\zeta$ is
\[
 -r^{-1}\operatorname{diag}(1,S,1)\,E_c^T
                      \operatorname{diag}(1,S^{-1},1),
\]
where $E_c^T$ is the transpose of $E_c$ written in the frame $(E_i)$
of the preceding wave, the matrix with entries $E_j\cdot E_cE_i$.
Since $S\ge1$, its norm is at
most $(S/r)|E_c|=|E_c|/(2\nu)\le\delta$.

Write $\zeta'=(A_0+A_1)\zeta$, where $A_0$ is the matrix of the
limiting system $\zeta_1'=-\zeta_2$, $\zeta_2'=-\varepsilon\zeta_3$,
$\zeta_3'=0$, and $A_1$ collects the remaining coefficients. The bounds
above give $\sup|A_1|\le CT\delta$. The vector $\zeta_0$ solves
$\zeta_0'=A_0\zeta_0$ with the same initial value as $\zeta$, and
$|\zeta_0|\le2T^2$. Since $A_0^3=0$, its propagator is
$\exp(uA_0)=\Id+uA_0+u^2A_0^2/2$, of norm at most $2T^2$ for
$0\le u\le T$. Therefore
\[
 \zeta(\tau)-\zeta_0(\tau)
 =\int_0^\tau e^{(\tau-\tau')A_0}A_1(\tau')\zeta(\tau')\,d\tau',
\]
so
\[
 \sup|\zeta-\zeta_0|
 \le CT^C\delta\,\bigl(1+\sup|\zeta-\zeta_0|\bigr).
\]
Condition \eqref{cc:eq:finite-smallness} gives the first bound of
\eqref{cc:eq:wave-vector-comparison}, hence \eqref{cc:eq:zeta-bounds}
and $m\ge1/4$. Since $m=\zeta_1^2+S^{-2}\zeta_2^2+\zeta_3^2$ and
$S^{-2}\le\delta$ by \eqref{cc:eq:S-bounds}, it also
gives $|m-\cos^{-2}\varphi|\le CT^C\delta$. Multiplication by
$\cos^2\varphi\le1$ proves the second bound.

To express the pressure coefficient in the amplitude equation for
$E_c=0$, put
\begin{align*}
 F_1={}&\zeta_1(b_{11}-b_{13}\zeta_1/\zeta_3)
       +\zeta_2(\sigma_p/\sigma+b_{21}/S-(b_{23}/S)\zeta_1/\zeta_3)\\
      &+\zeta_3(b_{31}-b_{33}\zeta_1/\zeta_3),\\
 F_2={}&\zeta_1(Sb_{12}-b_{13}\zeta_2/\zeta_3)
       +\zeta_2(b_{22}-(b_{23}/S)\zeta_2/\zeta_3)\\
      &+\zeta_3(Sb_{32}-b_{33}\zeta_2/\zeta_3).
\end{align*}
Then the pressure coefficient of \eqref{cc:eq:wave-ode} satisfies
\[
 \frac{\lambda S}{r}p_{\rm w}
   =\frac{2(F_1\alpha_1+F_2\alpha_2)}m.
\]
Substituting for $\alpha_3$ in the first two amplitude equations gives
a system for $(\alpha_1,\alpha_2)^T$ whose coefficient matrix $N$ has entries
\begin{align*}
 N_{11}&=-b_{11}+2b_{13}\zeta_1/\zeta_3+\boxed{2\zeta_1F_1/m},\\
 N_{12}&=\boxed{-2Sb_{12}}+2b_{13}\zeta_2/\zeta_3+\boxed{2\zeta_1F_2/m},\\
 N_{21}&=\boxed{-\sigma_p/\sigma}-(b_{21}-b_{12})/S
       +((b_{23}+b_{32})/S)\zeta_1/\zeta_3+2\zeta_2F_1/(S^2m),\\
 N_{22}&=-b_{22}+((b_{23}+b_{32})/S)\zeta_2/\zeta_3
                                      +2\zeta_2F_2/(S^2m).
\end{align*}
Only the boxed terms have nonzero reference values. Every other term
has one of $b_{11}$, $b_{22}$, $b_{13}$, the quotients $b_{ij}/S$ or
$S^{-2}$ as a factor, so its reference value is zero. In $F_1$ only
the term $\zeta_2\sigma_p/\sigma$, and in $F_2$ only $\zeta_1Sb_{12}$
and $\zeta_3Sb_{32}$, have nonzero reference values. Hence the
reference values of $F_1$, $F_2$ and $m$ are $-\varepsilon\tau$,
$(\tan\varphi+\varepsilon)/2$ and $\cos^{-2}\varphi$. Thus the
reference entries of $N$ are
\[
 -\varepsilon\tau\sin2\varphi,\qquad
 -1+\tan\varphi(\tan\varphi+\varepsilon)\cos^2\varphi
     =-\cos^2\varphi(1-\varepsilon\tan\varphi),\qquad -1,\qquad 0,
\]
which are the coefficients of \eqref{cc:eq:model-system}.

Each entry of $N$ is a polynomial of bounded degree in $\zeta_1$,
$\zeta_2$, $\zeta_3$, $1/\zeta_3$, $1/m$ and the quantities
\eqref{cc:eq:coefficient-list}. These are at most $CT^2$.
Each differs from its reference value by at most $CT^C\delta$: the
coordinates by \eqref{cc:eq:wave-vector-comparison}, the reciprocals
by a multiple of $|m-\cos^{-2}\varphi|+|\zeta_3-1|$, the
coefficients by the bounds at the beginning of the proof, and $S^{-2}$
by \eqref{cc:eq:S-bounds}. Hence each
entry of $N$ differs from its reference value by at most $CT^C\delta$.

For $E_c\ne0$, the amplitude equation contains the additional term
$-(\Id-2\frac{k}{|k|}\otimes\frac{k}{|k|})E_c\mathbf a$, of norm at
most $|E_c|\,|\mathbf a|$. By \eqref{cc:eq:scaled-pair},
\eqref{cc:eq:zeta-bounds} and \eqref{cc:eq:S-bounds},
$|\mathbf a|\le C|(\alpha_1,\alpha_2)|$, while
$\alpha_1=S\mathbf a\cdot E_1$ and $\alpha_2=\mathbf a\cdot E_2$.
In the time $\tau$ this adds a coefficient matrix of norm at most
$C(S/r)|E_c|\le C\delta$ to the reduced amplitude system.
\end{proof}

\subsubsection{The reference solution}\label{cc:sec:scalar}
At $\tau=0$, the reference system \eqref{cc:eq:model-system} is
$\alpha_1'=-\alpha_2$, $\alpha_2'=-\alpha_1$. Its growing solutions
have $\alpha_1=-\alpha_2$, the direction $-E_1/S+E_2$ of the amplitude
in \eqref{cc:eq:scaled-pair}. As $\varphi$ increases, the pressure
reduces the coefficient driving $\alpha_1$ from one to
$\cos^2\varphi(1-\varepsilon\tan\varphi)$ and introduces the term
$-\varepsilon\tau\sin2\varphi\,\alpha_1$. The driving coefficient
becomes negative when $\varepsilon\tan\varphi>1$.

The \emph{reference solution} is the solution of the reference system
with $(\alpha_1,\alpha_2)(0)=(-1,1)$, shown in
Figure~\ref{cc:fig:wave-turning}. A later wave is started with these
data in Subsection~\ref{cc:uf:sec:later-wave}. The second equation of
\eqref{cc:eq:model-system}, $\alpha_2'=-\alpha_1$, holds at every
$\tau$. Hence, wherever $\alpha_2>0$, the logarithmic derivative of
$\alpha_2$ is
\[
 Y=\frac{\alpha_2'}{\alpha_2}=-\frac{\alpha_1}{\alpha_2},
\]
and the first equation of \eqref{cc:eq:model-system} gives the Riccati
equation
\begin{equation}\label{cc:eq:model}
 Y'=\cos^2\varphi(1-\varepsilon\tan\varphi)
           -\varepsilon\tau\sin2\varphi\,Y-Y^2,\qquad Y(0)=1.
\end{equation}
Conversely, if $Y$ solves \eqref{cc:eq:model}, then
$(-Ya_{\rm ref},a_{\rm ref})$ is the reference solution, where
\[
 a_{\rm ref}(\tau)=\exp\left(\int_0^\tau Y(\tau')\,d\tau'\right).
\]
We call $a_{\rm ref}$ the \emph{reference amplitude}.

\begin{lemma}\label{cc:lem:model}
Let $0<\varepsilon\le1$. The solution $Y$ of \eqref{cc:eq:model}
exists for all $\tau\ge0$, satisfies $|Y|\le1$, and has the form
\begin{equation}\label{cc:eq:Y-ansatz}
 Y=(1+e)\cos\varphi-\tfrac14\varepsilon\tau\sin2\varphi,
 \qquad 0\le e\le\tfrac18\varepsilon.
\end{equation}
The reference amplitude satisfies
\begin{equation}\label{cc:eq:reference-amplitude}
 a_{\rm ref}(\tau)=\sqrt{\cos\varphi(\tau)}
       \exp\left(\int_0^\tau(1+e)\cos\varphi\,d\tau'\right).
\end{equation}
The function $Y$ has exactly one zero $\tau_m$, with
$2\le\varepsilon\tau_m\le3$. It is positive before this zero and
negative afterwards, so $a_{\rm ref}$ is largest at $\tau_m$, and
\[
 \log a_{\rm ref}(\tau_m)\ge\varepsilon^{-1/2}-1.
\]
\end{lemma}
\begin{proof}
Since $\tan\varphi=\varepsilon\tau^2/2$,
\[
 \varphi'=\varepsilon\tau\cos^2\varphi,\qquad
 (\log\cos\varphi)'=-\tfrac12\varepsilon\tau\sin2\varphi.
\]
We look for $Y$ in the form \eqref{cc:eq:Y-ansatz}, an ansatz with an
unknown function $e$. To motivate it, write $Y=\cos\varphi+w$, where
$\cos\varphi$ is the positive equilibrium of the equation
$Y'=\cos^2\varphi-Y^2$ with $\varphi$ frozen. Neglecting $w'$, $w^2$,
$\varepsilon\tau w\sin2\varphi$ and the term
$-\varepsilon\sin\varphi\cos\varphi$ of order $\varepsilon$,
\eqref{cc:eq:model} reduces to
$2w\cos\varphi=-\varepsilon\tau\sin\varphi\cos^2\varphi$, that is,
$w=-\tfrac14\varepsilon\tau\sin2\varphi=\tfrac12(\log\cos\varphi)'$.
The factor $1+e$ absorbs the neglected terms.

Integrating \eqref{cc:eq:Y-ansatz} and using the formula for
$(\log\cos\varphi)'$ gives \eqref{cc:eq:reference-amplitude}.
Substituting \eqref{cc:eq:Y-ansatz} into \eqref{cc:eq:model} gives
\[
 e'=\cos\varphi\bigl(-2e-e^2
                 +\tfrac12\varepsilon\sin\varphi\cos\varphi\bigr),
 \qquad e(0)=0.
\]
The right-hand side is nonnegative at $e=0$ and nonpositive at
$e=\varepsilon/8$, because $\sin\varphi\cos\varphi\le1/2$. Hence
$e$ remains in $[0,\varepsilon/8]$ and exists for all $\tau\ge0$.
This constructs the global solution $Y$.

Wherever $Y\ge0$, \eqref{cc:eq:model} gives $Y'\le1-Y^2$, because
$\cos^2\varphi(1-\varepsilon\tan\varphi)\le1$. Since $Y(0)=1$,
$Y(\tau)$ never exceeds one. Also
$\sin\varphi\cos\varphi\le\min\{1/2,\cot\varphi\}$, so for $\tau>0$,
\[
 \tfrac14\varepsilon\tau\sin2\varphi
 \le\min\{\varepsilon\tau/4,1/\tau\}\le\tfrac12.
\]
Equation \eqref{cc:eq:Y-ansatz} therefore gives $Y\ge-1/2$.

Write $Y=\cos\varphi(1+e-\tfrac12\varepsilon\tau\sin\varphi)$.
It is positive for $\varepsilon\tau\le2$, because $e\ge0$ and
$\sin\varphi<1$. For $\varepsilon\tau=3$, we have
$\tan\varphi=9/(2\varepsilon)\ge9/2$, hence $\sin\varphi>3/4$ and
$\tfrac12\varepsilon\tau\sin\varphi>9/8\ge1+e$. Thus $Y$ is negative
there. At any zero with $\varepsilon\tau\ge2$,
\[
 Y'=\cos^2\varphi(1-\varepsilon\tan\varphi)
    =\cos^2\varphi\bigl(1-(\varepsilon\tau)^2/2\bigr)<0.
\]
Consequently $Y$ has exactly one zero $\tau_m$, with
$2\le\varepsilon\tau_m\le3$, and has the stated signs on its two sides.

For $0\le\tau\le(2/\varepsilon)^{1/2}$, we have $\tan\varphi\le1$,
hence $\cos\varphi\ge2^{-1/2}$. This interval ends before
$\tau_m\ge2/\varepsilon$, where $a_{\rm ref}$ is largest. By
\eqref{cc:eq:reference-amplitude} and $e\ge0$,
\[
 \log a_{\rm ref}(\tau_m)
 \ge\log a_{\rm ref}\bigl((2/\varepsilon)^{1/2}\bigr)
 \ge2^{-1/2}(2/\varepsilon)^{1/2}-\tfrac14\log2
 \ge\varepsilon^{-1/2}-1.
 \qedhere
\]
\end{proof}

Relative to the reference solution, the reference system is triangular.
To compare a solution $(\alpha_1,\alpha_2)$ of
\eqref{cc:eq:model-system} with the reference solution
$(-Ya_{\rm ref},a_{\rm ref})$, we divide by the reference amplitude,
which removes the exponential growth common to the growing solutions,
and we separate the size of $\alpha_2$ from the deviation of the ratio
$\alpha_1/\alpha_2$ from its reference value $-Y$. We call
\begin{equation}\label{cc:eq:normalized}
 h=\frac{\alpha_2}{a_{\rm ref}},\qquad
 l=-\frac{\alpha_1+Y\alpha_2}{a_{\rm ref}\cos\varphi}
\end{equation}
the \emph{normalized coordinates} of $(\alpha_1,\alpha_2)$. The factor
$1/\cos\varphi$ in $l$ gives the equations for $h$ and $l$ below the
same factor $\cos\varphi$, which yields a propagator bound uniform in
$\varepsilon$ and in the length of the interval. By the second
equation of \eqref{cc:eq:model-system} and $a_{\rm ref}'=Ya_{\rm ref}$,
we have $h'=l\cos\varphi$. The trace of the reference coefficient matrix
is $-\varepsilon\tau\sin2\varphi=2(\log\cos\varphi)'$. By Liouville's
formula, the determinant
\[
 \det\begin{pmatrix}\alpha_1&-Ya_{\rm ref}\\
                      \alpha_2&a_{\rm ref}\end{pmatrix}
 =a_{\rm ref}(\alpha_1+Y\alpha_2)=-a_{\rm ref}^2l\cos\varphi
\]
is a constant multiple of $\cos^2\varphi$. Formula
\eqref{cc:eq:reference-amplitude} thus gives
\[
 l(\tau)=l(0)\exp\left(-2\int_0^\tau(1+e)\cos\varphi\,d\tau'\right).
\]
The reference system becomes
\begin{equation}\label{cc:eq:triangular-comparison}
 h'=l\cos\varphi,\qquad l'=-2(1+e)l\cos\varphi,
\end{equation}
and the reference solution is $(h,l)=(1,0)$. The second equation shows
that $|l|$ does not increase. Since $e\ge0$, for $\tau_1\le\tau_2$,
\[
 |h(\tau_2)-h(\tau_1)|\le\int_{\tau_1}^{\tau_2}|l|\cos\varphi
 \le\int_{\tau_1}^{\tau_2}(1+e)|l|\cos\varphi
 =\tfrac12\bigl(|l(\tau_1)|-|l(\tau_2)|\bigr).
\]
The triangle inequality and $|l(\tau_2)|\le|l(\tau_1)|$ give
\[
 |h(\tau_2)|+|l(\tau_2)|\le|h(\tau_1)|+|l(\tau_1)|.
\]
Since $|(h,l)|\le|h|+|l|\le\sqrt2\,|(h,l)|$, the propagator
$\Psi_0(\tau_2,\tau_1)$ of \eqref{cc:eq:triangular-comparison} has norm
at most $2$ for $\tau_1\le\tau_2$, uniformly in both times and in
$\varepsilon$.

\subsubsection{Comparison with the reference solution}
The next lemma shows that the wave follows the reference solution on
$0\le\tau\le T$. Dividing by $a_{\rm ref}$ removes the exponential
growth, and the bound for $\Psi_0$ gives errors polynomial in $T$.
A direct Gronwall estimate for $(\alpha_1,\alpha_2)$ would introduce
$\exp(CT)$, while \eqref{cc:eq:finite-smallness} controls polynomial
losses in $T$. The relative bounds give $|\mathbf a|$ up to a factor
$\exp(\pm\varepsilon)$, so the wave inherits the reference gain.
Recovering the normalized coordinates from an amplitude vector costs
a factor $S$, because its $E_1$ component is scaled by $S^{-1}$.

\begin{lemma}\label{cc:lem:finite-comparison}
Assume \eqref{cc:eq:finite-smallness}, and let $\mathbf a$ be the
amplitude vector with $\mathbf a(t_0)=-E_1(t_0)/S+E_2(t_0)$, that is,
with $(\alpha_1,\alpha_2)(0)=(-1,1)$ in \eqref{cc:eq:scaled-pair}, the
initial data of the reference solution.
Then $\alpha_2>0$ on $0\le\tau\le T$, and
\begin{equation}\label{cc:eq:amplitude-comparison}
 \left|\frac{\alpha_1}{\alpha_2}+Y\right|\le CT^C\delta,\qquad
 \alpha_2\le|\mathbf a|\le(1+CT^C\delta)\alpha_2.
\end{equation}
Moreover,
\begin{equation}\label{cc:eq:norm-rate-error}
 \left|\log\frac{|\mathbf a|}{a_{\rm ref}}\right|\le\varepsilon.
\end{equation}
For $0\le\tau_1\le\tau_2\le T$, let $\Phi(\tau_2,\tau_1)$ map
$\mathbf b(\tau_1)$ to $\mathbf b(\tau_2)$ for every amplitude vector
$\mathbf b$. It maps the plane perpendicular to $k(\tau_1)$ onto the
plane perpendicular to $k(\tau_2)$, and
\[
 \|\Phi(\tau_2,\tau_1)\|
 \le CS\,\frac{|\mathbf a(\tau_2)|}{|\mathbf a(\tau_1)|}.
\]
\end{lemma}
\begin{proof}
For an amplitude vector $\mathbf b$, let $\beta=(\beta_1,\beta_2,\beta_3)$
be its scaled components, defined as those of $\mathbf a$ in
\eqref{cc:eq:scaled-pair}:
\[
 \mathbf b=\frac{\beta_1}SE_1+\beta_2E_2+\frac{\beta_3}SE_3 .
\]
Its normalized coordinates are given by \eqref{cc:eq:normalized} with
$(\beta_1,\beta_2)$ in place of $(\alpha_1,\alpha_2)$, that is,
\[
 \binom hl=\frac1{a_{\rm ref}}R\binom{\beta_1}{\beta_2},\qquad
 R=\begin{pmatrix}0&1\\-(\cos\varphi)^{-1}&-Y(\cos\varphi)^{-1}\end{pmatrix},
 \qquad R^{-1}=\begin{pmatrix}-Y&-\cos\varphi\\1&0\end{pmatrix}.
\]
By Lemma~\ref{cc:uf:exact:lem:additive}, $(h,l)$ solves
\eqref{cc:eq:triangular-comparison} with the additional coefficient
matrix $R\mathcal ER^{-1}$. Since $|Y|\le1$ and
$(\cos\varphi)^{-1}\le T^2$, we have $|R|\le CT^2$, $|R^{-1}|\le C$
and $|R\mathcal ER^{-1}|\le CT^C\delta$.

Let $\Psi(\tau_2,\tau_1)$ be the propagator of this perturbed system.
Variation of constants gives
\[
 \Psi(\tau_2,\tau_1)=\Psi_0(\tau_2,\tau_1)
 +\int_{\tau_1}^{\tau_2}\Psi_0(\tau_2,\tau')
       (R\mathcal ER^{-1})(\tau')\Psi(\tau',\tau_1)\,d\tau'.
\]
Consequently $\sup|\Psi|\le2+2T\cdot CT^C\delta\sup|\Psi|$, and
\eqref{cc:eq:finite-smallness} gives $|\Psi(\tau_2,\tau_1)|\le4$.

For the prescribed initial amplitude $\mathbf a(t_0)$ the normalized
coordinates at $\tau=0$ are $(h,l)(0)=(1,0)$, a constant solution of the reference
system. Applying the preceding formula to this vector gives
\[
 \sup|(h-1,l)|\le2T\cdot CT^C\delta\cdot4\le CT^C\delta\le\tfrac12.
\]
In particular $h\ge1/2$, and the coordinate identities give
\[
 \alpha_2=a_{\rm ref}h>0,\qquad
 \frac{\alpha_1}{\alpha_2}+Y=-\frac{l\cos\varphi}{h},\qquad
 \left|\frac{\alpha_1}{\alpha_2}+Y\right|\le2|l|\le CT^C\delta.
\]
This is the first bound of \eqref{cc:eq:amplitude-comparison}, and
$|Y|\le1$ gives $|\alpha_1|\le2\alpha_2$. Since
$\alpha_3=-(\zeta_1\alpha_1+\zeta_2\alpha_2)/\zeta_3$, the bounds
\eqref{cc:eq:zeta-bounds} give $|\alpha_3|\le18T^2\alpha_2$. By
\eqref{cc:eq:scaled-pair},
\[
 |\mathbf a|^2=\alpha_2^2+S^{-2}(\alpha_1^2+\alpha_3^2).
\]
Since $S^{-2}\le\delta$ by \eqref{cc:eq:S-bounds}, this proves the
second bound of
\eqref{cc:eq:amplitude-comparison}. Thus
$\log(|\mathbf a|/a_{\rm ref})=\log h+\log(|\mathbf a|/\alpha_2)$
is at most $CT^C\delta$ in absolute value, and
\eqref{cc:eq:finite-smallness} gives \eqref{cc:eq:norm-rate-error}.

To bound $\Phi(\tau_2,\tau_1)$, note that for an arbitrary amplitude
vector $\mathbf b$,
$\beta_1=-a_{\rm ref}(Yh+l\cos\varphi)$, $\beta_2=a_{\rm ref}h$ and
$\beta_3=-(\zeta_1\beta_1+\zeta_2\beta_2)/\zeta_3$. Hence
$\mathbf b/a_{\rm ref}=hL_h+lL_l$, where
\[
 L_h=-\frac YSE_1+E_2+\frac{\zeta_1Y-\zeta_2}{S\zeta_3}E_3,\qquad
 L_l=\frac{\cos\varphi}{S}\left(-E_1+\frac{\zeta_1}{\zeta_3}E_3\right).
\]
By \eqref{cc:eq:zeta-bounds} and $S\ge T^2$, $|L_h|+|L_l|\le C$.
Conversely,
$h=(\mathbf b/a_{\rm ref})\cdot E_2$, and
$lL_l=\mathbf b/a_{\rm ref}-hL_h$. By
\eqref{cc:eq:wave-vector-comparison},
\[
 |L_l|=\frac{\cos\varphi}{S}
          \sqrt{1+\zeta_1^2/\zeta_3^2}\ge\frac cS,
 \qquad |h|+|l|\le CS|\mathbf b|/a_{\rm ref}.
\]
Together with $|\Psi|\le4$, these two coordinate bounds give
\[
 \|\Phi(\tau_2,\tau_1)\|
 \le CS\,\frac{a_{\rm ref}(\tau_2)}{a_{\rm ref}(\tau_1)}
 \le CS\,\frac{|\mathbf a(\tau_2)|}{|\mathbf a(\tau_1)|},
\]
where the last step uses \eqref{cc:eq:norm-rate-error}.
\end{proof}

\subsubsection{The couplings of the new wave}
Recall that the new wave evolves over the gradient
$M=B+\sigma_pE_2\otimes E_1+E_c$ of Subsection~\ref{cc:sec:comparison}.
The couplings $\nu_+,\eta_+$ of $M$ in the frame of the new wave drive
its successor. The first controls the sign of the pressure coefficient,
since
\[
 p_{\rm w}=\frac{2k\cdot M\mathbf a}{|k|^2}
          =2\nu_+\frac{|\mathbf a|}{|k|}.
\]
Subsection~\ref{cc:uf:sec:parent-step} uses this positivity throughout
the amplification interval.

\begin{lemma}\label{cc:lem:new-wave-couplings}
Assume \eqref{cc:eq:finite-smallness}, and let $\mathbf a$ be the
amplitude vector of Lemma~\ref{cc:lem:finite-comparison}. Put
$E_{1,+}=k/|k|$, $E_{2,+}=\mathbf a/|\mathbf a|$ and
$E_{3,+}=E_{1,+}\times E_{2,+}$, and define
\[
 \nu_+=E_{1,+}\cdot ME_{2,+},\qquad
 \eta_+=E_{3,+}\cdot ME_{2,+}.
\]
Then $\nu_+\ge\varepsilon\nu/2$ on $0\le\tau\le T$. Where
$\varepsilon\tan\varphi\ge1$, that is, for
$\sqrt2/\varepsilon\le\tau\le T$,
\begin{equation}\label{cc:eq:actual-ratio}
 \left|\frac{\nu_+}{\nu}-1\right|\le\varepsilon,\qquad
 \left|\frac{\eta_+/\nu_+}{\cos\varphi}-1\right|\le\frac12.
\end{equation}
\end{lemma}
\begin{proof}
Put
\[
 q=\frac{\alpha_1}{\alpha_2},\qquad
 z=\frac{\alpha_3}{\alpha_2}
   =-\frac{\zeta_1q+\zeta_2}{\zeta_3},\qquad n=\frac{|\mathbf a|}{\alpha_2}.
\]
Then
\[
 \frac{\mathbf a}{\alpha_2}=\frac qSE_1+E_2+\frac zSE_3,\qquad
 n^2=1+\frac{q^2+z^2}{S^2}.
\]
For $E_c=0$, the components of $M\mathbf a/(r\alpha_2)$ in the frame
$(E_1,E_2,E_3)$ of the preceding wave are $(f_1/S,f_2,f_3/S)$, where
\[
 f_1=b_{11}q+Sb_{12}+b_{13}z,\qquad
 f_2=(\sigma_p/\sigma+b_{21}/S)q+b_{22}+(b_{23}/S)z,\qquad
 f_3=b_{31}q+Sb_{32}+b_{33}z.
\]
In the same frame, the components of
$(k/\lambda)\times(\mathbf a/\alpha_2)$ are
\[
 \left(\frac{\zeta_2z}{S^2}-\zeta_3,
       \frac{\zeta_3q-\zeta_1z}{S},
       \zeta_1-\frac{\zeta_2q}{S^2}\right).
\]
Since $E_{1,+}=k/|k|$, $E_{2,+}=\mathbf a/|\mathbf a|$ and
$E_{3,+}=(k\times\mathbf a)/(|k||\mathbf a|)$,
\[
 \nu_+=\frac{k\cdot M\mathbf a}{|k||\mathbf a|},\qquad
 \eta_+=\frac{(k\times\mathbf a)\cdot M\mathbf a}{|k||\mathbf a|^2}.
\]
We compute these scalar products in the frame of the preceding wave,
pairing $k/\lambda$ and $(k/\lambda)\times(\mathbf a/\alpha_2)$ with
$M\mathbf a/(r\alpha_2)$. With $|k|=\lambda\sqrt m$,
$|\mathbf a|=n\alpha_2$ and $r/S=2\nu$ this gives
\begin{align}\label{cc:eq:coupling-rational}
 \nu_+/\nu&=\frac{2(\zeta_1f_1+\zeta_2f_2+\zeta_3f_3)}{\sqrt m\,n},\\
 \eta_+/\nu&=\frac{2}{\sqrt m\,n^2}
  \left\{\left(-\zeta_3+\frac{\zeta_2z}{S^2}\right)f_1
       +(\zeta_3q-\zeta_1z)f_2
       +\left(\zeta_1-\frac{\zeta_2q}{S^2}\right)f_3\right\}.
 \label{cc:eq:coupling-rational-end}
\end{align}
The matrix $E_c$ adds $E_{1,+}\cdot E_cE_{2,+}/\nu$ to $\nu_+/\nu$
and $E_{3,+}\cdot E_cE_{2,+}/\nu$ to $\eta_+/\nu$. Each is at most
$\delta$ in absolute value.

At the reference values,
\[
 \zeta=(\tan\varphi,-\varepsilon\tau,1),\qquad
 q=-Y,\qquad z=Y\tan\varphi+\varepsilon\tau,
\]
and $f_1=1/2$, $f_2=-Y$, $f_3=\varepsilon/2$, $n=1$ and
$m=\cos^{-2}\varphi$. The two expressions in
\eqref{cc:eq:coupling-rational}--\eqref{cc:eq:coupling-rational-end}
are therefore, respectively,
\[
 \sin\varphi+\varepsilon\cos\varphi+2\varepsilon\tau Y\cos\varphi,
 \qquad
 (\varepsilon\tan\varphi-1)\cos\varphi
    +\frac{2Y^2}{\cos\varphi}+2\varepsilon\tau Y\sin\varphi.
\]
In the second expression the terms $\varepsilon\tan\varphi\cos\varphi$,
$2Y^2/\cos\varphi$ and $2\varepsilon\tau Y\sin\varphi$ are large compared
with $\cos\varphi$ when $\varepsilon\tan\varphi$ is large, and
$\varepsilon\tan\varphi=\tfrac12(\varepsilon\tau)^2$ grows with $\tau$.
Their large parts cancel. To see this, write $Y=(1+e-u)\cos\varphi$
by \eqref{cc:eq:Y-ansatz}, with $u=\tfrac12\varepsilon\tau\sin\varphi$,
so that $2\varepsilon\tau\sin\varphi=4u$ and
$2u^2=\varepsilon\tan\varphi\sin^2\varphi$. Dividing the second
expression by $\cos\varphi$ gives
\begin{align*}
 &(\varepsilon\tan\varphi-1)+2(1+e-u)^2+4u(1+e-u)\\
 &\qquad=\varepsilon\tan\varphi-1
   +\bigl[2(1+e)^2-4u(1+e)+2u^2\bigr]+\bigl[4u(1+e)-4u^2\bigr]\\
 &\qquad=1+4e+2e^2+\bigl(\varepsilon\tan\varphi-2u^2\bigr).
\end{align*}
The terms $\pm4u(1+e)$ cancel, and the large parts of
$\varepsilon\tan\varphi$ and $2u^2$ cancel:
\[
 \varepsilon\tan\varphi-2u^2=\varepsilon\tan\varphi\,(1-\sin^2\varphi)
 =\varepsilon\sin\varphi\cos\varphi\le\tfrac12\varepsilon .
\]
Similarly,
$2\varepsilon\tau Y\cos\varphi
=2(1+e)\varepsilon\tau\cos^2\varphi
 -2\varepsilon\sin^2\varphi\cos\varphi$.
Thus the reference values of $\nu_+/\nu$ and $\eta_+/\nu$ are
\begin{equation}\label{cc:eq:reference-couplings}
 \begin{gathered}
 \sin\varphi+\varepsilon\cos\varphi\cos2\varphi
                   +2(1+e)\varepsilon\tau\cos^2\varphi,\\
 (1+4e+2e^2+\varepsilon\sin\varphi\cos\varphi)\cos\varphi.
 \end{gathered}
\end{equation}

Where $\varepsilon\tan\varphi\ge1$, we have $\cos\varphi\le\cot\varphi\le\varepsilon$, so
$1-\varepsilon^2\le\sin\varphi\le1$. Also $0\le e\le\varepsilon/8$ and
\[
 \varepsilon\tau\cos^2\varphi
 \le\varepsilon\tau\cot^2\varphi
 =\frac4{\varepsilon\tau^3}\le\sqrt2\,\varepsilon^2.
\]
The two reference values in \eqref{cc:eq:reference-couplings} are
$1+O(\varepsilon^2)$ and $(1+O(\varepsilon))\cos\varphi$.

Where $\varepsilon\tan\varphi\le1$, we have $\varepsilon\tau\le\sqrt2$, hence $Y>0$ by
Lemma~\ref{cc:lem:model}. The unsimplified reference value of
$\nu_+/\nu$ is at least $\sin\varphi+\varepsilon\cos\varphi
\ge\varepsilon$.

For the wave, in
\eqref{cc:eq:coupling-rational}--\eqref{cc:eq:coupling-rational-end},
the coefficients differ from their reference values by at most
$CT\delta$; $\zeta$ differs from $\zeta_0$ by at most $CT^C\delta$;
the terms with $S^{-2}$ in those formulas, in $m$ and in $n$, are at
most $CT^C\delta$; and $q$ differs from $-Y$ by at most $CT^C\delta$,
by \eqref{cc:eq:amplitude-comparison}. The formula for $z$ and the
bound $\zeta_3\ge1/2$ give the same bound for its change.
The numerators are polynomials of bounded degree in these quantities
and $1/\zeta_3$, which are at most $CT^2$, while the denominators
$\sqrt m\,n$ and $\sqrt m\,n^2$ are at least $1/2$.
Together with the contribution of $E_c$, the two ratios therefore
differ from their reference values by at most $CT^C\delta$, which is
at most $\varepsilon T^{-2}/8$ by \eqref{cc:eq:finite-smallness}.

In both ranges the reference value of $\nu_+/\nu$ is at least
$\varepsilon$; in the late range this uses $C_0\varepsilon\le1$.
Hence $\nu_+/\nu\ge\varepsilon/2$ on the whole interval. In the late
range, $\cos\varphi\ge T^{-2}$. Combining the perturbation bounds with
\eqref{cc:eq:reference-couplings} gives
\[
 \left|\frac{\nu_+}{\nu}-1\right|\le\varepsilon,\qquad
 \frac{\eta_+}{\nu_+}=(1+O(\varepsilon))\cos\varphi.
\]
The latter gives \eqref{cc:eq:actual-ratio} for $C_0$ sufficiently large.
\end{proof}

\subsubsection{The crossing}
After the peak of the reference solution at $\tau_m$, the new wave
becomes the shear for the next wave, whose coupling ratio
$\varepsilon_+$ is prescribed in advance. By
\eqref{cc:eq:actual-ratio}, the coupling ratio $\eta_+/\nu_+$ of the
new wave is comparable to $\cos\varphi$. It therefore reaches
$\varepsilon_+$ when $\cos\varphi$ is comparable to $\varepsilon_+$,
that is, when $\tan\varphi$ is comparable to $1/\varepsilon_+$, or
$\tau$ to $(\varepsilon\varepsilon_+)^{-1/2}$.
We call a time $\tau$ at which $\eta_+/\nu_+=\varepsilon_+$ a
\emph{crossing}: there the ratio reaches the threshold $\varepsilon_+$
prescribed for the next wave. The next wave is activated at a crossing. The reference
amplitude decreases after $\tau_m$, and the wave-vector length is
comparable to $\lambda/\cos\varphi$. The next proposition bounds their
product for the actual wave at a crossing, and shows that its amplitude
is small compared with its maximum during the first unit of normalized
time. The strength $\sigma_+=|\mathbf a||k|$ at the crossing and the
coupling $\nu_+$ supply the frozen parameters $\sigma,\nu$ for the next
wave.

\begin{proposition}\label{cc:prop:finite-return}
Assume \eqref{cc:eq:finite-smallness}, let $\mathbf a$ be the amplitude
vector of Lemma~\ref{cc:lem:finite-comparison}, with the couplings of
Lemma~\ref{cc:lem:new-wave-couplings}, and suppose that
$0<\varepsilon_+\le\varepsilon/32$ and
$T\ge2(\varepsilon\varepsilon_+)^{-1/2}$. Then the interval
$\tau_m\le\tau\le2(\varepsilon\varepsilon_+)^{-1/2}$, where
$\tau_m\ge2/\varepsilon$ is the zero of $Y$
(Lemma~\ref{cc:lem:model}), contains a crossing. At every crossing in
this interval the strength $\sigma_+=|\mathbf a(\tau)||k(\tau)|$ satisfies
\begin{equation}\label{cc:eq:general-return}
 C^{-1}(\varepsilon\varepsilon_+)^{-1/2}
 \le\frac{\sigma_+}{\lambda\max_{[0,T]}|\mathbf a|}
 \le C(\varepsilon\varepsilon_+)^{-1/2}.
\end{equation}
Moreover
\begin{equation}\label{cc:eq:gain}
 |\mathbf a(\tau)|\le\exp\bigl(3-\varepsilon^{-1/2}\bigr)
 \max_{[0,T]}|\mathbf a|\qquad(0\le\tau\le1).
\end{equation}
\end{proposition}
\begin{proof}
By Lemma~\ref{cc:lem:model}, $2/\varepsilon\le\tau_m\le3/\varepsilon$.
Since $\varepsilon_+\le\varepsilon/32$, we have
$(\varepsilon\varepsilon_+)^{-1/2}\ge4/\varepsilon$. Thus
$\tau_m<2(\varepsilon\varepsilon_+)^{-1/2}\le T$, and $[\tau_m,T]$ lies
in the range $\tau\ge\sqrt2/\varepsilon$ of \eqref{cc:eq:actual-ratio}.

For $\tau_m\le\tau\le T$, we have
$\tan\varphi=\varepsilon\tau^2/2\ge2/\varepsilon\ge1$, hence
\[
 \frac1{2\tan\varphi}\le\cos\varphi\le\frac1{\tan\varphi},\qquad
 \frac1{4\tan\varphi}\le\frac{\eta_+}{\nu_+}\le\frac3{2\tan\varphi},
\]
by \eqref{cc:eq:actual-ratio}. At $\tau_m$,
$\tan\varphi\le9/(2\varepsilon)$, so
$\eta_+/\nu_+\ge\varepsilon/18>\varepsilon_+$. At
$\tau=2(\varepsilon\varepsilon_+)^{-1/2}$, we have
$\tan\varphi=2/\varepsilon_+$, so
$\eta_+/\nu_+\le3\varepsilon_+/4<\varepsilon_+$.
Continuity therefore gives a crossing in this interval. At every
crossing in $[\tau_m,T]$ the displayed bounds give
$1/(4\varepsilon_+)\le\tan\varphi\le3/(2\varepsilon_+)$, so
$\cos\varphi$ is comparable to $\varepsilon_+$.

To estimate the strength at a crossing, we compare the amplitude there
with its maximum on $[0,T]$. By \eqref{cc:eq:norm-rate-error},
$|\mathbf a|$ lies between $\exp(-\varepsilon)a_{\rm ref}$ and
$\exp(\varepsilon)a_{\rm ref}$. The maximum of $a_{\rm ref}$ on
$[0,T]$ is attained at $\tau_m$, so for $\tau_m\le\tau\le T$,
$|\mathbf a(\tau)|/\max_{[0,T]}|\mathbf a|$ is bounded above and below
by absolute multiples of $a_{\rm ref}(\tau)/a_{\rm ref}(\tau_m)$.
Formula \eqref{cc:eq:reference-amplitude} gives
\[
 \frac{a_{\rm ref}(\tau)}{a_{\rm ref}(\tau_m)}
 =\left(\frac{\cos\varphi(\tau)}{\cos\varphi(\tau_m)}\right)^{1/2}
   \exp\left(\int_{\tau_m}^{\tau}(1+e)\cos\varphi\,d\tau'\right).
\]
The integral lies between zero and two, since
$\int_{\tau_m}^\infty\cos\varphi\,d\tau\le2/(\varepsilon\tau_m)\le1$.
Also $|k|$ is comparable to $\lambda/\cos\varphi$, by
\eqref{cc:eq:wave-vector-comparison}. Thus
$\sigma_+/(\lambda\max_{[0,T]}|\mathbf a|)$ is comparable to
$(\cos\varphi(\tau)\cos\varphi(\tau_m))^{-1/2}$.
Here $\cos\varphi(\tau_m)$ is comparable to $\varepsilon$, because
$2/\varepsilon\le\tan\varphi(\tau_m)\le9/(2\varepsilon)$, and at a
crossing $\cos\varphi$ is comparable to $\varepsilon_+$. This proves
\eqref{cc:eq:general-return}.

For $0\le\tau\le1$, the bound $Y\le1$ gives
$a_{\rm ref}(\tau)\le\exp(\tau)\le\exp(1)$, while
$\log a_{\rm ref}(\tau_m)\ge\varepsilon^{-1/2}-1$ by Lemma~\ref{cc:lem:model}. By \eqref{cc:eq:norm-rate-error} at such $\tau$ and at $\tau_m$,
$|\mathbf a(\tau)|\le\exp(1+\varepsilon)$ and
$|\mathbf a(\tau_m)|\ge\exp(\varepsilon^{-1/2}-1-\varepsilon)$.
Since $2\varepsilon\le1$, this gives \eqref{cc:eq:gain}.
\end{proof}

\section{Localized waves over a smooth parent}\label{cc:sec:localized}
We transport the phase and envelope of a Craik--Criminale wave by the
parent flow. The amplitude evolves along each particle path as in the
affine model. Its variation from particle to particle, together with the
spatial cutoff, requires a correction to keep the velocity divergence free.

\subsection{Lagrangian coordinates and the leading wave}\label{cc:sec:material}
Let $I=[0,T]$, $T>0$, and let $U$ be a smooth odd divergence-free
velocity on $I\times\T$. We need only its flow $X$ and its velocity
gradient. Write $a$ for a particle's position at time zero and
$x=X(t,a)$ for its position at time $t$. Set
\[
 X(0,a)=a,\qquad M=(D_xU)\circ X,\qquad
 \mathbf d=(D_aX)^{-T}\nabla_a .
\]
Time derivatives in these coordinates hold $a$ fixed.
On the lift to $\R^3$, $X(t,a+n)=X(t,a)+n$ for $n\in\mathbb Z^3$.

Fix the time-zero wave vector $k_*\in\R^3\setminus\{0\}$ and give the
particle with label $a$ the phase $\phi(a)=k_*\cdot a$.
The transported phase $\phi\circ X(t,\cdot)^{-1}$
has Eulerian gradient $\kappa=\mathbf d\phi=(D_aX)^{-T}k_*$ at $X(t,a)$.
For a flow linear in $a$ the transported phase is $\kappa\cdot x$, as in
Lemma~\ref{cc:lem:wave}.
A spatial derivative of a wave can fall on its amplitude or on its
oscillatory factor. To distinguish the two, write a profile as $f(t,a,s)$
with an independent phase $s\in\R/\mathbb Z$, measured in cycles, and set
\[
 f^\natural(t,X(t,a))=f(t,a,\phi(a)),\qquad
 \widetilde{\mathbf d}=\mathbf d+\kappa\partial_s .
\]
Here $\mathbf d$ differentiates at fixed $s$, and $\kappa\partial_s$
differentiates through the phase. Along the parent flow the phase is
constant, so the two chain rules are
\begin{equation}\label{cc:uf:hier:eq:graph}
 \nabla_x f^\natural=(\widetilde{\mathbf d}f)^\natural,\qquad
 (\partial_t+U\cdot\nabla_x)f^\natural=(f_t)^\natural.
\end{equation}
The phase modulo one is periodic in $a$ exactly when $k_*\in\mathbb Z^3$.
To allow arbitrary $k_*$, we evaluate the oscillatory profiles on the
lifted label cell $[-1/2,1/2]^3$ and extend periodically. They vanish near
the cell boundary, so the extension is smooth. Phase-independent profiles
are already periodic in $a$.

Along each particle path the local velocity gradient $M$ plays the role of
the matrix $D$ in \eqref{cc:eq:wave-ode}. Indeed, differentiating
$\kappa=(D_aX)^{-T}k_*$ and using $\partial_tD_aX=MD_aX$ gives
$\kappa_t=-M^T\kappa$. Choose a smooth envelope $g(a)$ supported in the
interior of the label cell and an \emph{activation time} $t_0\in[0,T)$.
On a neighborhood of $\supp g$ the amplitude, called the
\emph{seed}, solves the other Craik--Criminale equation,
\begin{equation}\label{cc:uf:step:eq:seed}
 \mathbf a_t=-M\mathbf a+2\kappa\frac{\kappa\cdot M\mathbf a}{|\kappa|^2},
 \qquad \mathbf a(t_0,a)=\mathbf a_*(a),\qquad
 \kappa(t_0,a)\cdot\mathbf a_*(a)=0 .
\end{equation}
The initial amplitude is smooth and may depend on $a$. For the first wave
it is constant and $t_0=0$. By \eqref{cc:eq:wave-ode}, the seed remains
perpendicular to $\kappa$ for $t_0\le t\le T$.

Let $F_0$ be a smooth odd $1$-periodic profile and $P_0$ its zero-mean
periodic primitive. The \emph{leading oscillation} and its pressure are
\begin{equation}\label{cc:uf:hier:eq:primary}
 A_1=g\mathbf a F_0,\qquad \mathsf p[A_1]=-p_{\rm w}gP_0,\qquad
 p_{\rm w}=2\frac{\kappa\cdot M\mathbf a}{|\kappa|^2}.
\end{equation}
Along the particle with label $a$ the leading oscillation is
$A_1^\natural(t,X(t,a))=g(a)\mathbf a(t,a)F_0(\phi(a))$.
The pressure balances the acceleration along $\kappa$. Multiplying
\eqref{cc:uf:step:eq:seed} by $gF_0$ gives
$A_{1,t}+MA_1=p_{\rm w}\kappa gF_0$, and $P_0'=F_0$ gives
\begin{equation}\label{cc:uf:hier:eq:leading}
 A_{1,t}+MA_1+\kappa\partial_s\mathsf p[A_1]=0 .
\end{equation}

The seed is transverse, so the phase derivative contributes zero to
$\nabla_x\cdot A_1^\natural$. Differentiating $g\mathbf a$ still leaves
$(\mathbf d\cdot A_1)^\natural$. To cancel this divergence while keeping
the wave localized, we look for a potential whose curl has fast part $A_1$.
By \eqref{cc:uf:hier:eq:graph}, the Eulerian curl of $Q^\natural$ is
$(\widetilde{\mathbf d}\times Q)^\natural$. We need the same correction for each later
oscillation. For a profile $V$ with $\kappa\cdot V=0$ and phase mean
$\av V=0$, the fast part of the curl is $\kappa\times\partial_s Q$.
To make it equal to $V$, take
\begin{equation}\label{cc:uf:hier:eq:completion}
 Q[V]=-\partial_s^{-1}\frac{\kappa\times V}{|\kappa|^2},
 \qquad C[V]=\mathbf d\times Q[V],
\end{equation}
where $\partial_s^{-1}$ is the zero-mean periodic primitive in $s$.
Since $\kappa\cdot V=0$, the vector triple product gives
$\kappa\times\partial_s Q[V]=V$. The other part of the curl is $C[V]$,
the \emph{divergence correction}. Its component along $\kappa$ determines
how it differentiates the phase of another oscillation. The identities
for the completed velocity and for this component are
\begin{equation}\label{cc:uf:hier:eq:slaving}
 V+C[V]=\widetilde{\mathbf d}\times Q[V],\qquad
 \kappa\cdot C[V]=-\partial_s^{-1}(\mathbf d\cdot V).
\end{equation}
For the second identity, apply $\widetilde{\mathbf d}\cdot$ to the first.
The $\mathbf d_i$ commute as pulled-back coordinate derivatives. Since
$\kappa_j=\mathbf d_j\phi$ and $\kappa$ is independent of $s$, the
$\widetilde{\mathbf d}_i$ commute as well. Thus both
$\widetilde{\mathbf d}\cdot(\widetilde{\mathbf d}\times Q[V])$ and
$\mathbf d\cdot C[V]$ vanish. Transversality leaves
$\partial_s(\kappa\cdot C[V])=-\mathbf d\cdot V$. Taking the zero-mean
primitive gives \eqref{cc:uf:hier:eq:slaving}, since $Q[V]$ and $C[V]$
have zero phase mean.

For the leading oscillation the potential is
\[
 Q[A_1]=gP_0\,\frac{\mathbf a\times\kappa}{|\kappa|^2}.
\]
The potential supplies a factor $|\kappa|^{-1}$, and its curl at fixed
phase costs one label derivative, as quantified in \eqref{cc:uf:hier:eq:completion-bound}.
Both $Q[V]$ and $C[V]$ vanish outside the label support of $V$ at each
time, since the primitive acts at fixed label and the curl is local in $a$.
Thus the corrected velocity $(A_1+C[A_1])^\natural$ is a curl supported
in $X(t,\supp g)$. The leading pressure $\mathsf p[A_1]^\natural$ also
vanishes outside this set.
In the construction $g$ is even and equals one near zero,
and $\mathbf a_*$ is even in $a$. The odd flow gives even $M$ and
$\kappa$, so the seed is even in $a$. Since $F_0$ is odd and $P_0$ is
even, the completed Eulerian velocity is odd.

A later wave must switch on at its activation time. Take a smooth time
cutoff $0\le\widetilde\chi\le1$ on $[0,T]$ vanishing for $t<t_0$,
with $\widetilde\chi\equiv1$ for the first wave. Set
\begin{equation}\label{cc:uf:hier:eq:switched}
 V_1=\widetilde\chi A_1,\qquad
 \mathsf p_1=\widetilde\chi\mathsf p[A_1],\qquad
 C_2=\widetilde\chi C[A_1]=C[V_1].
\end{equation}
For $t_0>0$ these products extend smoothly by zero across $t_0$
because the cutoff is flat there. Multiplication by a time cutoff preserves
the spatial divergence cancellation. Its derivative along the parent
flow contributes the term
\begin{equation}\label{cc:uf:hier:eq:switch}
 R_{\widetilde\chi}=\widetilde\chi'\,(A_1+C[A_1])
\end{equation}
to the residual, so the cutoff error depends on the size of the wave
while it is being switched on.

\subsection{The profile}\label{cc:sec:profile}
We choose $F_0=F_h$ to prescribe the central velocity gradient and
control the wave pressure. Since $g(0)=1$ and $F_h(0)=0$, taking
$F_h'(0)=1$ gives the leading oscillation $A_1^\natural$ the gradient
$\mathbf a\otimes\kappa$ at the origin. In the affine model the wave
pressure has Hessian $-p_{\rm w}F_h'(k\cdot x)k\otimes k$. When $p_{\rm w}>0$,
the bound $F_h'\ge-Ch$ controls its largest eigenvalue, as in
\eqref{cc:intro:profile-curvature}.

\begin{lemma}\label{cc:uf:exact:lem:profile}
There are odd smooth \(1\)-periodic functions
\(F_h\), \(0<h\le1\), with zero-mean primitives \(P_h\), such that
\[
 F_h'(0)=1,\qquad -Ch\le F_h'\le1,\qquad
 \norm{F_h}_\infty\le Ch,\qquad
 \norm{\partial_s^jF_h}_\infty\le C_jh^{1-j}\quad(j\ge2).
\]
\end{lemma}
\begin{proof}
Fix an even $\psi\in C_c^\infty((-1/4,1/4))$ with $0\le\psi\le1$ and
$\psi(0)=1$. The bump $\psi(s/h)$ has integral
$c_h=h\int_{\R}\psi<h/2$ over one period. Subtract its mean to obtain a
periodic primitive, and normalize the slope at zero by setting
\[
 F_h'(s)=\frac{\psi(s/h)-c_h}{1-c_h}\quad(|s|\le1/2),\qquad F_h(0)=0.
\]
The derivative is even with mean zero, so $F_h$ is odd and periodic
and has a zero-mean periodic primitive $P_h$.
Since $c_h<1/2$, the slope equals one at zero and lies between $-2c_h$
and one. Integrating the lower bound from $0$ to $s$ and from $s$ to $1$,
using $F_h(0)=F_h(1)=0$, gives $\norm{F_h}_\infty\le 2c_h$.
Differentiating $\psi(s/h)$ another $j-1$ times gives the bound $C_jh^{1-j}$.
\end{proof}

\section{Selecting initial data by the endpoint problem}\label{cc:sec:endpoint}
For a prescribed force $f$ in the linearized Euler equation
\eqref{cc:intro:linearized-momentum}, we choose the initial value $b(0)$
of the phase-independent correction to preserve each particle's final
position to first order. The first-order displacement $\xi$ induced by
$b$ vanishes initially because the starting positions are fixed.
We impose $\xi(T)=0$ to select $b(0)$.

Fix $0<T\le1$ and $I=[0,T]$. Let the parent $u$ be a smooth
divergence-free field on $I\times\mathbb T^3$, not assumed to solve any
equation, and let $X(t,a)$, $X(0,a)=a$, be its volume-preserving flow. Set
\begin{equation}\label{cc:uf:mean:data}
 M=(\nabla u)\circ X,\quad
 H^\#=-(M_t+M^2),\quad
 K=(D_aX)^TD_aX,\quad \Omega=(D_aX)^T(M^T-M)D_aX.
\end{equation}
Time derivatives keep $a$ fixed. Choose
$D_-\ge\max\{1,\sup_{I,a}|(D_aX)^{-1}|\}$.
The flow equation gives
\begin{equation}\label{cc:uf:mean:flow-identities}
 \partial_t D_aX=MD_aX,\qquad
 K_t-\Omega=2(D_aX)^T\partial_t D_aX.
\end{equation}
Let $p$ be a smooth function and write $R=u_t+u\cdot\nabla u+\nabla p$.
The gradient of this identity is
$(\partial_t+u\cdot\nabla)\nabla u+(\nabla u)^2+D^2p=\nabla R$, with
$(\nabla u)_{ij}=\partial_ju_i$. Along trajectories the material
derivative becomes $\partial_t$, so
\begin{equation}\label{cc:uf:mean:eq:Hsharp}
 H^\#=(D^2p-\nabla R)\circ X.
\end{equation}
For a solution of the Euler equations, $H^\#$ is the pressure Hessian
along the flow: it is what prevents the velocity gradient from evolving
by the Riccati equation $M_t=-M^2$, and it is symmetric; for a linear
flow $Dx$ this is the condition that $D'+D^2$ be symmetric in
Subsection~\ref{cc:sec:affine}. In general $H^\#$ need not be symmetric.
The matrix $M^T-M$ is the vorticity of $u$ along the flow:
$(M^T-M)y=-(\curl u\circ X)\times y$. Thus $\Omega$ is the vorticity
pulled back to labels. When $H^\#\le\Id$ as quadratic forms, the
displacement energy is coercive and the inverse in Appendix~\ref{cc:app:mean}
gives a unique correction with a H\"older bound.

\subsection{Displacement and the two end conditions}

Here $H^1_0(I;\R^3)$ is the space of $z\in H^1(I;\R^3)$ with
$z(0)=z(T)=0$, the closure of $C^\infty_c((0,T);\R^3)$ in $H^1$. Fix a
label $a$, let $z\in H^1_0(I;\R^3)$, and put $\xi=D_aX(\cdot,a)z$. Since
$\partial_t D_aX=MD_aX$,
\begin{equation}\label{cc:uf:mean:xi-t}
 \xi_t=(\partial_t D_aX)z+D_aX z_t=M\xi+D_aX z_t,\qquad
 z_t=(D_aX)^{-1}(\xi_t-M\xi).
\end{equation}

By the first variation formula \eqref{cc:intro:first-variation}, the
displacement induced by $b$ satisfies $\xi_t=M\xi+b\circ X$.
Thus \eqref{cc:uf:mean:xi-t} identifies $D_aX z_t=b\circ X$.
Since $D_aX$ is invertible, the prescribed final displacement
$\xi(T)=0$ is equivalent to $z(T)=0$.

Cauchy--Schwarz in $z(t)=-\int_t^T z_t(s)\dd s$ gives
\begin{equation}\label{cc:uf:mean:poincare}
 |z(t)|^2\le(T-t)\int_I|z_t|^2\dd s,\qquad
 \int_I|z|^2\dd t\le\frac{T^2}{2}\int_I|z_t|^2\dd t.
\end{equation}
Since $D_aX$ is smooth, $\xi\in H^1_0(I;\R^3)$ as well, so
\eqref{cc:uf:mean:poincare} also holds for $\xi$; in particular
$\|\xi\|_{L^2(I)}\le\|\xi_t\|_{L^2(I)}$, since $T\le1$. By
\eqref{cc:uf:mean:xi-t} and $|(D_aX)^{-1}|\le D_-$,
\[
 \|z_t\|_{L^2(I)}
 \le D_-\bigl(\|\xi_t\|_{L^2(I)}+\|M\|_\infty\|\xi\|_{L^2(I)}\bigr)
 \le D_-(1+\|M\|_\infty)\|\xi_t\|_{L^2(I)}.
\]
Squaring gives
\begin{equation}\label{cc:uf:mean:equivalence}
 \int_I|z_t|^2\dd t\le D_-^2(1+\|M\|_\infty)^2\int_I|\xi_t|^2\dd t.
\end{equation}

\subsection{Coercivity}

The inverse in Appendix~\ref{cc:app:mean} requires the label form
$\int_I(Kz_t\cdot z_t+\Omega z_t\cdot z)\dd t$ to control
$\int_I|z_t|^2\dd t$. In the displacement variable $\xi$, the integrand
contains $H^\#$ and a total derivative that integrates to zero because
$\xi$ vanishes at both ends. We can then use the upper bound on $H^\#$ and
\eqref{cc:uf:mean:equivalence} to obtain the required estimate.

Substitute $D_aX z_t=\xi_t-M\xi$ and $D_aX z=\xi$ in the two terms. First,
$Kz_t\cdot z_t=|D_aX z_t|^2=|\xi_t|^2-2\xi_t\cdot M\xi+|M\xi|^2$. Second,
$\Omega z_t\cdot z=(M^T-M)(\xi_t-M\xi)\cdot\xi
=\xi_t\cdot M\xi-\xi_t\cdot M^T\xi-|M\xi|^2+M^2\xi\cdot\xi$. Adding,
and using
$\partial_t(\xi\cdot M\xi)=\xi_t\cdot M\xi+\xi_t\cdot M^T\xi+M_t\xi\cdot\xi$
and $M_t+M^2=-H^\#$, we obtain
\[
 Kz_t\cdot z_t+\Omega z_t\cdot z
 =|\xi_t|^2-H^\#\xi\cdot\xi-\partial_t(\xi\cdot M\xi).
\]
Since $\xi$ vanishes at both ends, integration and density give
\begin{equation}\label{cc:uf:mean:jacobi-id}
 \int_I(Kz_t\cdot z_t+\Omega z_t\cdot z)\dd t
 =\int_I(|\xi_t|^2-H^\#\xi\cdot\xi)\dd t.
\end{equation}

\begin{lemma}\label{cc:uf:mean:hardy}
Assume $H^\#\le\Id$ on $I\times\mathbb T^3$ as quadratic forms. Then
for every $a\in\mathbb T^3$, $z\in H^1_0(I;\R^3)$ and
$\xi=D_aX(\cdot,a)z$, one has
$\int_I H^\#\xi\cdot\xi\dd t\le\frac12\int_I|\xi_t|^2\dd t$. Hence
\begin{equation}\label{cc:uf:mean:coercive}
 \int_I(Kz_t\cdot z_t+\Omega z_t\cdot z)\dd t
 \ge\frac12\int_I|\xi_t|^2\dd t.
\end{equation}
\end{lemma}
\begin{proof}
The assumption $H^\#\le\Id$, the second inequality of
\eqref{cc:uf:mean:poincare} for $\xi$ and $T\le1$ give
\[
 \int_I H^\#\xi\cdot\xi\dd t
 \le\int_I|\xi|^2\dd t
 \le\frac{T^2}2\int_I|\xi_t|^2\dd t
 \le\frac12\int_I|\xi_t|^2\dd t.
\]
Equation~\eqref{cc:uf:mean:jacobi-id} gives
\eqref{cc:uf:mean:coercive}.
\end{proof}

\subsection{The inverse and its pressure}\label{cc:uf:sec:mean}

Given a source $L(t,a)$, we solve the system
\begin{equation}\label{cc:uf:mean:integrated}
 K(t)v(t)-\int_0^t(\Omega v)(s)\dd s+L(t)=c+\nabla_a\phi,
 \qquad v=z_t,\quad z(0)=z(T)=0,\quad\diver_a z=0,
\end{equation}
for the displacement $z$, a time-independent divergence-free field $c$
and a label pressure gradient $\nabla_a\phi$. For the force sources of
Theorem~\ref{cc:uf:mean:inverse}, \eqref{cc:uf:mean:integrated} is the
linearized momentum equation \eqref{cc:uf:mean:strong} integrated once
in time. For these sources $L(0)=0$. Since $K(0)=\Id$, the equation
at time zero reads $v(0)=c+\nabla_a\phi(0)$. Projecting
onto divergence-free fields gives $c=v(0)$. Since $b\circ X=D_aXv$
and $X(0,a)=a$, the constant $c$ is the initial correction $b(0)$.

\begin{theorem}\label{cc:uf:mean:inverse}
Assume $H^\#\le\Id$ on $I\times\mathbb T^3$ as quadratic forms. For
every $L\in C(I;\Holder{0})$ the system
\eqref{cc:uf:mean:integrated} has a unique solution with
$z\in C^1(I;\Holder{0})$, $c\in\Holder{0}$ divergence free and
$\nabla_a\phi\in C(I;\Holder{0})$. There are absolute $C,p$ such that, with
\begin{equation}\label{cc:uf:mean:base-constant}
 A=C\bigl(D_-+\|M\|_\infty+\norm K_{C_t\Holder{1}}
   +\norm\Omega_{C_t\Holder{1}}\bigr)^p,
\end{equation}
one has $\norm v_{C_t\Holder{0}}\le A\norm L_{C_t\Holder{0}}$.

Let $f$ be a smooth field on $I\times\mathbb T^3$, and let
$L=-\int_0^tg(s)\dd s$ with $g=(D_aX)^T(f\circ X)$. Then $z$ is smooth.
Put $B=D_aX v$ and $\mathcal D=2(D_aX)^T\partial_t D_aX$.
There is a smooth zero-mean label pressure $\mathsf q$ such that
\begin{equation}\label{cc:uf:mean:strong}
 Kv_t+\nabla_a\mathsf q=g-\mathcal Dv.
\end{equation}
Thus $b=B\circ X^{-1}$ and $q=\mathsf q\circ X^{-1}$ are smooth and solve
\begin{equation}\label{cc:uf:mean:physical}
 \diver b=0,\qquad
 b_t+u\cdot\nabla b+b\cdot\nabla u+\nabla q=f.
\end{equation}
The endpoint conditions select the initial velocity $b(0)$.

The number $A$ also bounds, uniformly in $t$, the operator norm on
$\Holder{0}$ of the fixed-time inverse $(\PP K(t))^{-1}\PP$, where
$(\PP K(t))^{-1}$ inverts the map $v\mapsto\PP(K(t)v)$ on the
divergence-free fields of $\Holder{0}$.
\end{theorem}
\begin{proof}
Rescale time by $\tau=t/T$ and put $Z(\tau)=z(t)/T$, so
$Z_\tau=v$. Theorem~\ref{cc:holder:linear} applies with coefficients
$K,T\Omega$, ellipticity $D_-^{-2}$ and source $L$.
Its coercivity constant is $1/\bigl(2D_-^2(1+\|M\|_\infty)^2\bigr)$,
by Lemma~\ref{cc:uf:mean:hardy} and \eqref{cc:uf:mean:equivalence}.
It gives the unique solution. Since $D_-\ge1$ and $T\le1$, the quantity
$Q$ of Appendix~\ref{cc:app:mean} is at most an absolute constant times
a fixed power of the sum in parentheses in
\eqref{cc:uf:mean:base-constant}, so the estimate
\eqref{cc:holder:eq:base-linear} gives the asserted bound.
For a smooth $f$ the source $L$ is smooth, and so is the solution, by
the last assertion of Theorem~\ref{cc:holder:linear}.
Differentiating \eqref{cc:uf:mean:integrated}, using
$K_t-\Omega=\mathcal D$ from \eqref{cc:uf:mean:flow-identities}, and
setting $\mathsf q=-\phi_t$ gives
\eqref{cc:uf:mean:strong}. Then
$B_t+MB=D_aX v_t+2(\partial_t D_aX)v$ and the Piola identity
$\diver_x(B\circ X^{-1})=(\diver_a v)\circ X^{-1}$ give
\eqref{cc:uf:mean:physical}.
Finally, by Lemma~\ref{cc:holder:freezing} with $H=\R^3$, $\mathcal R=\PP$ and
$B=K(t)$, whose coercivity is $D_-^{-2}$, the fixed-time inverse is
bounded by a fixed power of $2+D_-^2+\|K(t)\|_{C^1}$. Since
$D_-\ge1$, this number is also at most an absolute constant times a
fixed power of the sum in parentheses in
\eqref{cc:uf:mean:base-constant}, so $A$ bounds the fixed-time inverse
once $C$ and $p$ are enlarged.
\end{proof}

\section{From a localized wave to an approximate solution}\label{cc:sec:exact}
We correct the localized wave of Section~\ref{cc:sec:localized} to obtain
the parent velocity for the next wave. The corrections must preserve the
leading wave's rank-one contribution to the velocity gradient at the
origin, where the next wave is amplified. The new parent must also have
controlled flow deformation and satisfy the pressure and residual bound
\eqref{cc:uf:tower:eq:parent-coercive} for the next endpoint problem.

At each level of the expansion we first remove the phase average of the
error using Theorem~\ref{cc:uf:mean:inverse}. The first mean correction
also cancels the residual of the parent. Transport of the leading wave
by the new mean then enters the equation for the oscillatory correction.
The residual of the new parent has three sources. The depth controls the
remainder of the expansion. The wave is small compared with its peak
while it is switched on, so its gain controls the error from the time
cutoff. The final Fourier cutoff contributes the error estimated in
\eqref{cc:uf:step:eq:truncated}. Section~\ref{cc:uf:sec:tower} chooses
the parameters to make all three errors small on the scale of the next
wave.

\subsection{A step over a truncated parent}\label{cc:uf:sec:parent-step}
Throughout this section the parent is the truncated velocity to which
we add the new wave. We prescribe the wave vector and seed at the origin
and extend the seed to the labels in the envelope.
Section~\ref{cc:sec:iteration} verifies the two assumptions
\eqref{cc:uf:tower:eq:wave-origin} on the central pair for the first
wave and for each later wave.

We use the Fourier modes $e^{2\pi i\xi\cdot x}$ for $\xi\in\mathbb Z^3$.
A field is \emph{supported on frequencies at most $\Lambda$} if its
Fourier coefficients vanish for $|\xi|>\Lambda$. Let
the parent $U$ be a smooth odd divergence-free field on
$[0,T]\times\mathbb T^3$, $T\le1$, and use its flow $X$, the matrix $M$
and the operators $\mathbf d$, $\widetilde{\mathbf d}$ of
Subsection~\ref{cc:sec:material}. Let $\mathsf D\ge2$ be a number such
that $U$ is supported on frequencies at most $2\mathsf D$ and
\begin{equation}\label{cc:uf:tower:eq:parent-bound}
 \norm{\nabla U}_\infty\le\mathsf D,\qquad
 \sup|(D_aX)^{\pm1}|\le\mathsf D .
\end{equation}
Since $U$ is odd, the flow fixes the origin, and $D_aX$ and $M$ are even
in $a$. Let $p_U$ be the function of zero mean with
$\Delta p_U=-\tr((\nabla U)^2)$, and let
$R_U=U_t+U\cdot\nabla U+\nabla p_U$ be the residual of the parent; no
other equation is assumed. Thus $R_U=U_t+\PP\diver(U\otimes U)$ is odd,
divergence free and supported on frequencies at most $4\mathsf D$.
To apply the endpoint inverse of Theorem~\ref{cc:uf:mean:inverse}
to the mean correction over $U$, assume that
\begin{equation}\label{cc:uf:tower:eq:parent-coercive}
 D^2p_U-\nabla R_U\le\Id\qquad\hbox{on }[0,T]\times\mathbb T^3
\end{equation}
as quadratic forms. By \eqref{cc:uf:mean:eq:Hsharp} this is the
hypothesis $H^\#\le\Id$ of that theorem.

The wave is that of Subsection~\ref{cc:sec:material} over $U$, with its
constant vector $k_*\ne0$ and transported field
$\kappa=(D_aX)^{-T}k_*$, its activation time $t_0\in[0,T)$ and its cutoff
in time $\widetilde\chi$. Let $\lambda\ge1$ be the frequency of the wave:
its wave vector at the origin has length $|\kappa(t_0,0)|=\lambda$ at the
activation time. Since $k_*=(D_aX)^T\kappa$ at every time and label, the
bound \eqref{cc:uf:tower:eq:parent-bound} for the deformation gives
\begin{equation}\label{cc:uf:tower:eq:wave-vector}
 \mathsf D^{-1}\lambda\le|k_*|\le\mathsf D\lambda,\qquad
 \mathsf D^{-2}\lambda\le|\kappa|\le\mathsf D^2\lambda .
\end{equation}
The profile is $F_0=F_h$ of Lemma~\ref{cc:uf:exact:lem:profile}, with
$h=\lambda^{-1/2}$. At the origin the seed $\mathbf a(t,0)$,
$t_0\le t\le T$, is a solution of the seed equation
\eqref{cc:uf:step:eq:seed} at the label $a=0$, perpendicular to
$\kappa(t,0)$ and normalized by
$\max_{t_0\le t\le T}|\mathbf a(t,0)|=\lambda^{-3/4}$. It does not
vanish, because it is a nonzero solution of a linear ordinary
differential equation. Put
\[
 \gamma(t)=\lambda^{3/4}|\mathbf a(t,0)|,\qquad
 \nu(t,0)=\Bigl(\frac{\kappa}{|\kappa|}\cdot
          M\frac{\mathbf a}{|\mathbf a|}\Bigr)(t,0)
 \qquad(t_0\le t\le T).
\]
As in Lemma~\ref{cc:lem:frame}, $\nu(t,0)$ is the first coupling of the
pair $(\kappa,\mathbf a)$ at the origin over the parent, and
$|\mathbf a(t,0)||\kappa(t,0)|$ is its strength. At the other labels
the seed is the solution of \eqref{cc:uf:step:eq:seed} on $[t_0,T]$ with
\begin{equation}\label{cc:uf:tower:eq:start}
 \mathbf a(t_0,a)=\bigl(\Id-P_{\kappa(t_0,a)}\bigr)\mathbf a(t_0,0),
 \qquad P_\kappa=\frac{\kappa\otimes\kappa}{|\kappa|^2},
\end{equation}
the projection of the vector at the origin onto the plane perpendicular
to $\kappa(t_0,a)$. It is even in $a$.

Fix a smooth function $\chi$ on $[0,\infty)$ with $0\le\chi\le1$,
$\chi=1$ on $[0,1]$ and $\chi=0$ on $[2,\infty)$, and take the even
envelope $g(a)=\chi(|a|/\ell_*)$, equal to one on the ball of radius
$\ell_*$ specified below. The leading oscillation is the field
$A_1=g\mathbf aF_h$ of \eqref{cc:uf:hier:eq:primary}, and it enters
through $V_1=\widetilde\chi A_1$ of \eqref{cc:uf:hier:eq:switched}.

At the origin, assume that for $t_0\le t'\le t\le T$,
\begin{equation}\label{cc:uf:tower:eq:wave-origin}
 \frac{|\mathbf b(t)|}{\gamma(t)}
 \le\mathsf D\,\frac{|\mathbf b(t')|}{\gamma(t')},\qquad
 \nu(t,0)\ge\mathsf D^{-1},
\end{equation}
the first for every solution $\mathbf b$ of the seed equation at the
origin that is perpendicular to $\kappa(\cdot,0)$. The first inequality
bounds the growth of every transverse solution by that of the seed,
up to the factor $\mathsf D$. On a sufficiently small envelope this
gives the weighted estimate for the transverse corrections, and the
second inequality keeps $\kappa\cdot M\mathbf a$ positive there
(Lemma~\ref{cc:uf:tower:lem:inverse}). With this sign, the lower bound
$F_h'\ge-Ch$ controls the leading pressure tensor in
\eqref{cc:uf:exact:eq:primary-hess} from above.
Both inequalities can be imposed before choosing the frequency and the
depth of the expansion.

On the envelope, \eqref{cc:uf:tower:eq:near-origin} bounds the seed
by twice its size at the origin. Since $\|F_h\|_\infty\le Ch$, the leading
velocity has size at most $C\lambda^{-3/4}h=C\lambda^{-5/4}$.
The first mean equation removes both the quadratic wave forcing and
$R_U$, so we require the parent residual to be at most the square of
this scale,
\begin{equation}\label{cc:uf:tower:eq:parent-small}
 \norm{R_U}_{C_tL^\infty}\le\lambda^{-5/2}.
\end{equation}
The radius of the envelope is a negative power of $\mathsf D$ with an
absolute exponent, $\ell_*\le\mathsf D^{-3}$. The exponent is
fixed sufficiently large for the finitely many estimates of
Subsections~\ref{cc:uf:trunc:sec:parent}--\ref{cc:uf:sec:transverse}.
In particular the support of $g$ lies in the interior of
$[-1/2,1/2]^3$.

Constants and exponents written $C$ are absolute: they do not depend on
the parent, the wave, an order of differentiation or the depth of the
expansion. A constant written $C_n$ depends only on the order $n$, on
the fixed profiles, and on the constants of the same kind in the
hypotheses of the statement in which it appears; $C_{n,p}$ and $C_{m,P}$
depend also on a level $p$ or on the depth $P$. These constants may
change from one estimate to the next, and no bound on their growth in
$n$ is used. In the iteration these constants and the powers of
$\mathsf D$ are at most a small power of $\lambda$; see
\eqref{cc:par:step}.

\subsection{Coefficient bounds for the parent}\label{cc:uf:trunc:sec:parent}
The coefficients of the two linear problems are built from the
deformation $D_aX$ of the parent, its velocity gradient $M$ and the wave
vector $\kappa$. Bernstein's inequality \eqref{cc:uf:trunc:eq:lp} for
the truncated parent and the differentiated flow equation bound each
label derivative of these coefficients at the cost of a fixed power
of $\mathsf D$. The norm $\|\cdot\|_{[n]}$ of
\eqref{cc:uf:trunc:eq:field-norm} weights each label derivative by
the radius $\ell_*$ and each phase derivative by the width $h$.
For a coefficient $c(t,a)$ it is $\sum_{m\le n}\ell_*^m|c|_m$, where
$|c|_m$ sums the suprema of its label derivatives of order $m$.
Taking $\ell_*$ sufficiently small absorbs the cost of each derivative,
so the exponent of $\mathsf D$ in \eqref{cc:uf:hier:eq:coefficient-bound}
is independent of $n$.

\begin{proposition}\label{cc:uf:trunc:prop:parent}
There is an absolute exponent $C$ such that each coefficient $c$ among
\[
 (D_aX)^{\pm1},\qquad M,\qquad
 \kappa/\lambda,\qquad \lambda^2/|\kappa|^2,\qquad P_\kappa
\]
satisfies \eqref{cc:uf:trunc:eq:power}. Let $f(t,x)$ be a field with
values in a normed space such that $|f|_m\le C_m'A_f\mathsf D^m$ for
every $m\ge0$ and a number $A_f$. Then
\begin{equation}\label{cc:uf:trunc:eq:parent-composition}
 |f\circ X|_m\le C_mA_f\mathsf D^{Cm}\qquad(m\ge0).
\end{equation}
Consequently, for every $n\ge0$,
\begin{equation}\label{cc:uf:hier:eq:coefficient-bound}
 \|c\|_{[n]}\le C_n\mathsf D^C,\qquad
 \|f\circ X\|_{[n]}\le C_nA_f .
\end{equation}
\end{proposition}
\begin{proof}
For $n\ge1$ let $\upsilon_n=\max_i|X_i|_n$, so
$\upsilon_1\le3\mathsf D$. By \eqref{cc:uf:trunc:eq:lp} for $\nabla U$,
with values in the continuous functions of time, and
\eqref{cc:uf:tower:eq:parent-bound},
$|U|_k\le C_k\mathsf D^{k-1}\norm{\nabla U}_\infty\le C_k\mathsf D^k$
for $k\ge1$. For $|\beta|=n\ge2$, the chain rule gives
\[
 \partial_t\partial^\beta X=(\nabla U\circ X)\partial^\beta X+N_\beta,
\]
where $N_\beta$ is a sum, with coefficients that depend only on $\beta$,
of the terms
\[
 \bigl((\partial^{\beta'}U)\circ X\bigr)
 \prod_{j=1}^k\partial^{\beta_j}X_{l_j},\qquad
 2\le k=|\beta'|\le n,\quad \beta_1+\dots+\beta_k=\beta,\quad|\beta_j|\ge1.
\]
Since $\partial_t D_aX=(\nabla U\circ X)D_aX$ and
$\partial^\beta X(0)=0$,
\[
 \partial^\beta X(t)
 =\int_0^tD_aX(t)(D_aX(\sigma))^{-1}N_\beta(\sigma)\dd\sigma.
\]
By \eqref{cc:uf:tower:eq:parent-bound} the propagator has norm at most
$\mathsf D^2$, and $T\le1$. Hence, for $n\ge2$,
\begin{equation}\label{cc:uf:trunc:eq:recursion}
 \upsilon_n\le C_n\mathsf D^2\sum_{k=2}^n\mathsf D^k
 \sum_{n_1+\dots+n_k=n}\upsilon_{n_1}\cdots\upsilon_{n_k},
\end{equation}
with every $n_i\ge1$. Induction on $n$ gives
$\upsilon_n\le C_n\mathsf D^{5n-4}$: by the bounds at the orders
$n_i<n$, a term of \eqref{cc:uf:trunc:eq:recursion} with $k$ factors is
at most $C_n$ times
\[
 \mathsf D^{2+k}\,\mathsf D^{5n-4k}
 =\mathsf D^{5n-4-3(k-2)}\le\mathsf D^{5n-4} .
\]
Thus
$|D_aX|_m\le C_m\mathsf D^{5m+1}$ for every $m\ge0$. For $f$ as in
the statement, the chain rule bounds $|f\circ X|_n$, $n\ge1$, by $C_n$
times the terms $A_f\mathsf D^k\mathsf D^{5n-4k}$ with
$1\le k\le n$, which proves \eqref{cc:uf:trunc:eq:parent-composition}.
By \eqref{cc:uf:trunc:eq:lp}, $f=\nabla U$ satisfies the hypothesis with
$A_f=\norm{\nabla U}_\infty\le\mathsf D$, which bounds $M$.

The product and inverse rules of Lemma~\ref{cc:uf:trunc:lem:coefficients}(a),(b)
now apply successively to $(D_aX)^{-1}$, $\kappa/\lambda=(D_aX)^{-T}k_*/\lambda$,
$|\kappa/\lambda|^{-2}$ and $P_\kappa$.
Indeed, \eqref{cc:uf:tower:eq:parent-bound} bounds the inverse deformation,
and \eqref{cc:uf:tower:eq:wave-vector} gives $|k_*|/\lambda\le\mathsf D$
and $|\kappa|/\lambda\ge\mathsf D^{-2}$.
Finally, $\ell_*$ is at most the power $\mathsf D^{-C}$ that occurs
here. So part (c) gives the first bound of
\eqref{cc:uf:hier:eq:coefficient-bound}, and
$\ell_*^m|f\circ X|_m\le C_mA_f(\mathsf D^C\ell_*)^m\le C_mA_f$ gives
the second.
\end{proof}

In particular the residual of the parent, which is supported on
frequencies at most $4\mathsf D$, satisfies the hypothesis of
\eqref{cc:uf:trunc:eq:parent-composition} with $A_f=\norm{R_U}_{C_tL^\infty}$,
by \eqref{cc:uf:trunc:eq:lp}. With \eqref{cc:uf:tower:eq:parent-small}
this gives, for every $n\ge0$,
\begin{equation}\label{cc:uf:hier:eq:parent-residual}
 \|R_U\circ X\|_{[n]}\le C_n\lambda^{-5/2}.
\end{equation}

\subsection{Finite-order estimates for the mean}\label{cc:uf:mean:sec:finite}
The phase-independent corrections of the expansion solve the linearized
Euler equations over the parent, with zero displacement at both ends of
$[0,T]$. For a smooth phase-independent source \(\bar f(t,a)\),
Theorem~\ref{cc:uf:mean:inverse} over $U$
gives \((B,\bar q)\) with
\begin{equation}\label{cc:uf:hier:eq:mean}
 B_t+MB+\mathbf d\bar q=\bar f,\qquad \mathbf d\cdot B=0,\qquad
 z(0)=z(T)=0,\quad B=D_aX z_t.
\end{equation}
In the notation of that theorem, $\bar f=f\circ X$,
$g=(D_aX)^T\bar f$, $v=z_t$, and $\bar q$ is the pressure $\mathsf q$ of
\eqref{cc:uf:mean:strong}. We call \eqref{cc:uf:hier:eq:mean} the
\emph{mean equation}. The terminal condition selects the initial
velocity $B(0)$, so the mean can already be present before activation.
We estimate it on all of $[0,T]$ in the unweighted label norms.
The inverse preserves the derivative order in the scaled H\"older
norms. Passing to the supremum norms by \eqref{cc:uf:trunc:eq:holder-sup}
costs one derivative of the source. At each level of the expansion
the mean source is formed from terms already estimated at all orders.

\begin{proposition}\label{cc:uf:mean:columns}
There is an absolute exponent $C$ such that, for every $n\ge0$, every
solution of the mean equation \eqref{cc:uf:hier:eq:mean} satisfies
\begin{equation}\label{cc:uf:mean:bound}
 \norm B_{[n]}\le C_n\mathsf D^C\norm{\bar f}_{[n+1]}.
\end{equation}
\end{proposition}
\begin{proof}
By Proposition~\ref{cc:uf:trunc:prop:parent} and
Lemma~\ref{cc:uf:trunc:lem:coefficients}(a), the coefficients $K$ and
$\Omega$ of \eqref{cc:uf:mean:data}, which enter
\eqref{cc:uf:mean:integrated}, satisfy \eqref{cc:uf:trunc:eq:power}.
Theorem~\ref{cc:uf:mean:inverse} applies over $U$ with $D_-=\mathsf D$,
by \eqref{cc:uf:tower:eq:parent-bound}. Its constant $A$ of
\eqref{cc:uf:mean:base-constant} contains $D_-$,
$\|M\|_\infty\le\mathsf D$ and the norms of $K$ and $\Omega$ in
$C_t\Holder{1}$, so it is at most $\mathsf D^C$. Since $\ell_*$ is a
sufficiently large negative power of $\mathsf D$, the first bound of
\eqref{cc:uf:trunc:eq:holder-sup} and
Lemma~\ref{cc:uf:trunc:lem:coefficients}(c) give
\begin{equation}\label{cc:uf:mean:thresholds}
 \ell_*^m\bigl(|K|_{m,1/2}+|\Omega|_{m,1/2}\bigr)
 \le C_m\mathsf D^C\ell_*^{1/2}\le C_m/A\qquad(m\ge1).
\end{equation}
Label derivatives preserve the endpoint and divergence constraints,
the time independence of the divergence-free field, and the class of
label gradients in \eqref{cc:uf:mean:integrated}. Theorem~\ref{cc:uf:mean:inverse}
bounds the map $L\mapsto v$ by $A$ in $C_t\Holder{0}$. The term $Kv$
enters directly, and $\Omega v$ enters through time integration,
an operator of norm at most $T\le1$.
Thus \eqref{cc:uf:mean:thresholds} gives the coefficient bounds for
Lemma~\ref{cc:uf:trunc:lem:lifting} with $C_{\rm low}=A$.
Applying its H\"older estimate to the source $L=-\int_0^tg(t')\dd t'$
of \eqref{cc:uf:mean:integrated} gives
\begin{equation}\label{cc:uf:mean:holder-bound}
 \norm v_{[n],1/2}\le C_n\mathsf D^C\norm g_{[n],1/2}.
\end{equation}
By \eqref{cc:uf:trunc:eq:holder-sup} and \eqref{cc:uf:hier:eq:coefficient-bound},
$\|D_aX\|_{[n],1/2}\le C_n\ell_*^{-1/2}\mathsf D^C\le C_n\mathsf D^{C'}$,
with $C'$ absolute. Using $B=D_aXv$ and $g=(D_aX)^T\bar f$ in
\eqref{cc:uf:mean:holder-bound}, the product estimate \eqref{cc:uf:trunc:eq:product}
gives $\|B\|_{[n],1/2}\le C_n\mathsf D^C\|\bar f\|_{[n],1/2}$.
The comparison \eqref{cc:uf:trunc:eq:holder-sup} now gives \eqref{cc:uf:mean:bound}.
\end{proof}

\subsection{The transverse equation}\label{cc:uf:sec:transverse}
Each oscillatory correction solves a forced seed equation forward from
zero initial data at each label and phase. Its forcing vanishes before
activation, so the correction does too. From activation onward we measure
its size relative to the amplitude of the seed at the origin, using the
weight $\gamma$.

We call
\begin{equation}\label{cc:uf:hier:eq:trans}
 V_t+MV+\kappa\partial_s\mathsf p=f,\qquad
 \kappa\cdot V=0,\qquad \av V=\av f=0,
\end{equation}
the \emph{transverse equation}, because its unknown $V$ is perpendicular
to the wave vector $\kappa$. The pressure term keeps $V$
perpendicular to $\kappa$: differentiating $\kappa\cdot V=0$
and using $\kappa_t=-M^T\kappa$ in \eqref{cc:uf:hier:eq:trans} gives
\begin{equation}\label{cc:uf:hier:eq:pressure-from-trans}
 |\kappa|^2\partial_s\mathsf p=\kappa\cdot f-2\kappa\cdot MV.
\end{equation}
The transverse equation is therefore the linear ordinary differential
equation
\begin{equation}\label{cc:uf:hier:eq:trans-ode}
 V_t=-(\Id-2P_\kappa)MV+(\Id-P_\kappa)f
\end{equation}
at each fixed $(a,s)$; for $f=0$ it is the seed equation
\eqref{cc:uf:step:eq:seed}. For
every source we take its solution with $V(0,a,s)=0$, and $\mathsf p$ is
the zero-mean primitive in $s$ of
\eqref{cc:uf:hier:eq:pressure-from-trans}. For smooth $f$ the solution
is unique and smooth, and $V$ satisfies the two constraints in
\eqref{cc:uf:hier:eq:trans}. Both $V$ and $\mathsf p$ vanish at every
label where $f$ vanishes for all $(t,s)$. If $f$ vanishes for $t<t_0$,
then $V=0$ for $t\le t_0$.

The plane $\kappa^\perp$ varies with the label. To compare the equation
on the envelope with the one at the origin and to differentiate it in
the labels, we express $V$ in coordinates in a fixed $\R^2$.
Choose a constant $3\times2$ matrix with orthonormal columns spanning
$k_*^\perp$, let $A_\perp$ be $D_aX$ times this matrix, and let
$A_\perp^\dagger=(A_\perp^TA_\perp)^{-1}A_\perp^T$.
The columns of $A_\perp$ span $\kappa^\perp$, since
$\kappa=(D_aX)^{-T}k_*$, and $A_\perp^\dagger$ takes a vector to its
coordinates in those columns after orthogonal projection onto their
span: $A_\perp^\dagger A_\perp=\Id_{\R^2}$ and
$A_\perp A_\perp^\dagger=\Id-P_\kappa$. Write $V=A_\perp v$ for a
transverse velocity. Since $\partial_tA_\perp=MA_\perp$, the equation
\eqref{cc:uf:hier:eq:trans-ode} becomes
\begin{equation}\label{cc:uf:tower:eq:coordinate}
 v_t+\mathcal M_\perp v=\psi,\qquad
 \mathcal M_\perp=2A_\perp^\dagger MA_\perp,\qquad
 \psi=A_\perp^\dagger f.
\end{equation}

\begin{lemma}\label{cc:uf:tower:lem:inverse}
For $t_0\le t\le T$ and $|a|\le2\ell_*$, the seed satisfies
\begin{equation}\label{cc:uf:tower:eq:near-origin}
 |\mathbf a(t,a)|\le2|\mathbf a(t,0)|,\qquad
 (\kappa\cdot M\mathbf a)(t,a)>0 .
\end{equation}
There is an absolute exponent $C$ such that the following holds at every
order $n\ge0$.
\begin{enumerate}[label=(\roman*),leftmargin=*]
\item Let $V$ and $f$ be smooth fields on $[t_0,T]$ that satisfy
\eqref{cc:uf:hier:eq:trans-ode}, with $V$ perpendicular to $\kappa$, and
suppose that $V(t_0)$ and $f$ vanish at the labels outside $\supp g$.
Then
\begin{equation}\label{cc:uf:tower:eq:seed-norm}
 \|V\|_{[n],\gamma}\le C_n\mathsf D^C
 \bigl(\gamma(t_0)^{-1}\|V(t_0)\|_{[n]}+\|f\|_{[n],\gamma}\bigr).
\end{equation}
\item Let $f$ be smooth on $[0,T]$ with zero phase mean, and suppose
that $f$ vanishes for $t<t_0$ and at the labels outside $\supp g$.
The solution of \eqref{cc:uf:hier:eq:trans} with $V(0)=0$ and its
pressure satisfy
\begin{equation}\label{cc:uf:hier:eq:trans-bound}
 \|V\|_{[n],\gamma}+\lambda\|\mathsf p\|_{[n],\gamma}
 \le C_n\mathsf D^C\|f\|_{[n],\gamma}.
\end{equation}
\item On $[t_0,T]$ the leading oscillation and its pressure satisfy
\begin{equation}\label{cc:uf:hier:eq:size}
 \|A_1\|_{[n],\gamma}+\lambda\|\mathsf p[A_1]\|_{[n],\gamma}
 \le C_n\mathsf D^C\lambda^{-5/4}.
\end{equation}
\end{enumerate}
\end{lemma}
\begin{proof}
By \eqref{cc:uf:tower:eq:parent-bound},
$\mathsf D^{-1}|w|\le|A_\perp w|\le\mathsf D|w|$ for $w\in\R^2$. Hence
$|A_\perp|\le\mathsf D$ and $|A_\perp^\dagger|\le\mathsf D$, and the
matrix $K_\perp=A_\perp^TA_\perp$ satisfies
$\sup|K_\perp^{-1}|\le\mathsf D^2$. By
Proposition~\ref{cc:uf:trunc:prop:parent} and
Lemma~\ref{cc:uf:trunc:lem:coefficients}(a),(b), the coefficients
$A_\perp$, $A_\perp^\dagger=K_\perp^{-1}A_\perp^T$ and $\mathcal M_\perp$
satisfy \eqref{cc:uf:trunc:eq:power}.

In the coordinates $V=A_\perp v$, the seed equation at the label $a$ is
$v_t+\mathcal M_\perp(t,a)v=0$. With the bounds for $A_\perp$ and
$A_\perp^\dagger$, the first condition
in \eqref{cc:uf:tower:eq:wave-origin} gives
$|v(t)|/\gamma(t)\le\mathsf D^3|v(t')|/\gamma(t')$ for its solutions at
$a=0$. By the bound for the first label
derivatives of $\mathcal M_\perp$ and the choice of $\ell_*$,
\[
 \sup_{t_0\le t\le T}
 |\mathcal M_\perp(t,a)-\mathcal M_\perp(t,0)|
 \le C\mathsf D^C|a|\le\mathsf D^{-3}\qquad(|a|\le2\ell_*).
\]
Let $v$ solve the equation at such a label $a$ on $[t',T]$. It solves
the equation at the origin with
the source $-(\mathcal M_\perp(\cdot,a)-\mathcal M_\perp(\cdot,0))v$, so
variation of constants gives
\[
 \frac{|v(t)|}{\gamma(t)}
 \le\mathsf D^3\frac{|v(t')|}{\gamma(t')}
 +\int_{t'}^t\frac{|v(t'')|}{\gamma(t'')}\dd t''.
\]
The weights cancel because the bound at the origin is a ratio of values
of $\gamma$. Since $t-t'\le T\le1$, Gronwall's inequality gives
$|v(t)|/\gamma(t)\le e\mathsf D^3|v(t')|/\gamma(t')$. For a solution of
\eqref{cc:uf:tower:eq:coordinate} on $[t_0,T]$, variation of constants
therefore gives, at each phase and each label with $|a|\le2\ell_*$,
\begin{equation}\label{cc:uf:tower:eq:weighted-forward}
 \frac{|v(t,a)|}{\gamma(t)}
 \le e\mathsf D^3\left(
 \frac{|v(t_0,a)|}{\gamma(t_0)}
 +\int_{t_0}^t\frac{|\psi(t',a)|}{\gamma(t')}\,dt'\right).
\end{equation}

To compare the seed at $a$ with its value at the origin, write
$\mathbf a=A_\perp v_a$, with $v_a=A_\perp^\dagger\mathbf a$ at $a$. The
difference $v_a-v_0$ solves the equation at the origin with the source
$-(\mathcal M_\perp(\cdot,a)-\mathcal M_\perp(\cdot,0))v_a$. Its value at
$t_0$ is $(A_\perp^\dagger(t_0,a)-A_\perp^\dagger(t_0,0))\mathbf a(t_0,0)$,
by \eqref{cc:uf:tower:eq:start} and $A_\perp^\dagger P_\kappa=0$, and
$|\mathbf a(t_0,0)|=\gamma(t_0)\lambda^{-3/4}$. With the bound just
proved for $v_a$, variation of constants gives
$|v_a(t)-v_0(t)|\le C\mathsf D^C|a|\gamma(t)\lambda^{-3/4}$. Since
$\gamma(t)\lambda^{-3/4}=|\mathbf a(t,0)|$, it follows that
\[
 |\mathbf a(t,a)-\mathbf a(t,0)|
 \le C\mathsf D^C|a|\,|\mathbf a(t,0)|,
 \qquad t_0\le t\le T.
\]
The label derivatives also give
$|\kappa(t,a)-\kappa(t,0)|\le C\mathsf D^C|a||\kappa(t,0)|$
and $|M(t,a)-M(t,0)|\le C\mathsf D^C|a|$. Since $|M|\le\mathsf D$, the
three bounds give
\[
 \bigl|(\kappa\cdot M\mathbf a)(t,a)-(\kappa\cdot M\mathbf a)(t,0)\bigr|
 \le C\mathsf D^C|a|\,|\kappa(t,0)|\,|\mathbf a(t,0)| ,
\]
while $(\kappa\cdot M\mathbf a)(t,0)
=\nu(t,0)|\kappa(t,0)||\mathbf a(t,0)|
\ge\mathsf D^{-1}|\kappa(t,0)||\mathbf a(t,0)|$, by the second condition
in \eqref{cc:uf:tower:eq:wave-origin}. None of these exponents depends
on the exponent of $\ell_*$, so the choice of $\ell_*$ proves
\eqref{cc:uf:tower:eq:near-origin}.

Let $V$ and $f$ be as in (i). The solution vanishes at the labels
outside $\supp g$, where its initial value and source vanish, and
$\supp g$ lies in $\{|a|\le2\ell_*\}$. Taking the supremum in
\eqref{cc:uf:tower:eq:weighted-forward} gives, because
$t-t_0\le T\le1$,
\[
 \|v\|_{[0],\gamma}\le e\mathsf D^3\bigl(\gamma(t_0)^{-1}\|v(t_0)\|_{[0]}
       +\|\psi\|_{[0],\gamma}\bigr).
\]
Since $\supp g$ is closed, the label and phase derivatives of $v$ and
$\psi$ also vanish outside it. So this is the zero-order
estimate of Lemma~\ref{cc:uf:trunc:lem:lifting} on $[t_0,T]$, for the
initial value and the source together, with the one coefficient
$\mathcal M_\perp$, which enters through the identity. By the choice of
$\ell_*$, condition \eqref{cc:uf:trunc:eq:slow} holds, and
the lemma gives
$\|v\|_{[n],\gamma}\le C_n\mathsf D^3(\gamma(t_0)^{-1}\|v(t_0)\|_{[n]}
+\|\psi\|_{[n],\gamma})$.
The estimate divides by $\gamma$ after differentiation in label and
phase. Since $V=A_\perp v$, $v=A_\perp^\dagger V$ and
$\psi=A_\perp^\dagger f$, the product estimate
\eqref{cc:uf:trunc:eq:product} and
Lemma~\ref{cc:uf:trunc:lem:coefficients}(c) for $A_\perp$ and
$A_\perp^\dagger$ give \eqref{cc:uf:tower:eq:seed-norm}.

In (ii) the forcing vanishes for $t<t_0$, so $V(t_0)=0$.
Applying \eqref{cc:uf:tower:eq:seed-norm} to
\eqref{cc:uf:hier:eq:trans-ode} bounds $V$. By
\eqref{cc:uf:hier:eq:pressure-from-trans}, $\lambda\mathsf p$ is the
phase primitive of a product of the coefficients $\kappa/\lambda$ and
$\lambda^2/|\kappa|^2$ with $f-2MV$, so
\eqref{cc:uf:trunc:eq:phase-ops}, \eqref{cc:uf:trunc:eq:product} and
\eqref{cc:uf:hier:eq:coefficient-bound} give
\begin{equation}\label{cc:uf:hier:eq:pressure-bound}
 \lambda\|\mathsf p\|_{[n],\gamma}
 \le C_n\mathsf D^C\bigl(\|f\|_{[n],\gamma}+\|V\|_{[n],\gamma}\bigr),
\end{equation}
which proves \eqref{cc:uf:hier:eq:trans-bound}.

For (iii), the time independence of the envelope implies that $g\mathbf a$
solves the seed equation on $[t_0,T]$ at every label. It is perpendicular
to $\kappa$ and vanishes outside $\supp g$. By
\eqref{cc:uf:tower:eq:start}, its value at
the time $t_0$ is $g(\Id-P_\kappa)\mathbf a(t_0,0)$. The chain rule
gives $|g|_m\le C_m\ell_*^{-m}$, so $\|g\|_{[n]}\le C_n$. With
\eqref{cc:uf:hier:eq:coefficient-bound} for
$P_\kappa$ and \eqref{cc:uf:trunc:eq:product}, this gives
$\|g\mathbf a(t_0,\cdot)\|_{[n]}
\le C_n\mathsf D^C|\mathbf a(t_0,0)|=C_n\mathsf D^C\gamma(t_0)\lambda^{-3/4}$.
Hence \eqref{cc:uf:tower:eq:seed-norm} with $f=0$ gives
$\|g\mathbf a\|_{[n],\gamma}\le C_n\mathsf D^C\lambda^{-3/4}$.
For a function of the phase alone,
$\|F\|_{[n]}=\sum_{j\le n}h^j\|\partial_s^jF\|_\infty$, so
Lemma~\ref{cc:uf:exact:lem:profile} gives $\|F_h\|_{[n]}\le C_nh$.
By \eqref{cc:uf:trunc:eq:product},
$\|A_1\|_{[n],\gamma}\le C_n\mathsf D^C h\lambda^{-3/4}
=C_n\mathsf D^C\lambda^{-5/4}$. By \eqref{cc:uf:hier:eq:leading}, on
$[t_0,T]$ the leading oscillation solves \eqref{cc:uf:hier:eq:trans}
with $f=0$ and $\mathsf p=\mathsf p[A_1]$, so
\eqref{cc:uf:hier:eq:pressure-bound} gives \eqref{cc:uf:hier:eq:size}.
\end{proof}

\subsection{The finite expansion}\label{cc:uf:sec:hierarchy}
At depth $P$ we seek a perturbation
\[
 W=\sum_{p=1}^P Z_p+C_{P+1},\qquad
 Z_p=Y_p+B_p,\qquad Y_p=V_p+C_p,\qquad C_1=0,
\]
where $V_p$ is perpendicular to $\kappa$ and $B_p$ is independent of phase.
The divergence correction $C_{p+1}=C[V_p]$ makes
$V_p+C_{p+1}=\widetilde{\mathbf d}\times Q[V_p]$ divergence free by
\eqref{cc:uf:hier:eq:slaving}, which accounts for the last term $C_{P+1}$.
We call $p$ the \emph{level} of $Z_p$ and of $Y_p$, $V_p$, $C_p$ and
$B_p$. The factor $|\kappa|^{-1}$ in the vector potential gives the
divergence correction $C_p$ of $V_{p-1}$ the bound of a term of
level $p$, so we group it with $V_p$.

\subsubsection{The equations at each level}

Start from $V_1,\mathsf p_1,C_2$ in
\eqref{cc:uf:hier:eq:switched} and put $B_1=0$, so that $Y_1=Z_1=V_1$.
Set aside the activation error $R_{\widetilde\chi}$ of
\eqref{cc:uf:hier:eq:switch}. The terms left by the leading wave at level
two are
\begin{equation}\label{cc:uf:hier:eq:H2}
 H_2=\widetilde\chi\bigl[(\partial_t+M)C[A_1]+\mathbf d\mathsf p[A_1]\bigr]
   +(V_1\cdot\mathbf d)V_1+(C_2\cdot\kappa)\partial_sV_1 .
\end{equation}
The linear term has zero phase average because $C[A_1]$ and
$\mathsf p[A_1]$ do. In the phase average of the last term,
\eqref{cc:uf:hier:eq:slaving} and integration by parts over a period give
\[
 \av{(C_2\cdot\kappa)\partial_sV_1}
 =-\av{\partial_s(C_2\cdot\kappa)V_1}
 =\av{(\mathbf d\cdot V_1)V_1}.
\]
Thus both quadratic terms in $\av{H_2}$ have a single label
derivative and two factors of $V_1$. Each factor has size
$C_n\mathsf D^C\lambda^{-5/4}$ by \eqref{cc:uf:hier:eq:size}, and the
label derivative costs a power of $\mathsf D$ by \eqref{cc:uf:trunc:eq:deriv}.
Hence $\|\av{H_2}\|_{[n]}\le C_n\mathsf D^C\lambda^{-5/2}$.
The parent residual $R_U\circ X$ satisfies a bound of this same order by
\eqref{cc:uf:hier:eq:parent-residual}, so it can be included in the
first mean equation.

In general, a product $(Z_i\cdot\mathbf d)Z_j$ has level $i+j$.
The fast derivative costs one level, so $(Z_i\cdot\kappa)\partial_sY_j$
has level $i+j-1$. Since $V_i\cdot\kappa=0$, only $C_i$ and $B_i$ can
differentiate the phase. We collect the terms involving
$(B_i\cdot\kappa)\partial_sY_j$ after the mean solve. The other
terms at level $p\ge3$ give
\begin{align}
 H_p={}&(\partial_t+M)C_p+\mathbf d\mathsf p_{p-1}\nonumber\\
 &+\sum_{i+j=p}(Z_i\cdot\mathbf d)Z_j
   +\sum_{\substack{i+j=p+1\\i\ge2}}
                  (C_i\cdot\kappa)\partial_sY_j,\label{cc:uf:hier:eq:Hp}
\end{align}
where all summation indices are positive integers.
For every $p\ge2$ put
\begin{equation}\label{cc:uf:hier:eq:source}
 f_p=-(H_p-\av{H_p}),\qquad \bar f_p=-\av{H_p}\quad(p\ge3),\qquad
 \bar f_2=-\av{H_2}-R_U\circ X.
\end{equation}
Every profile in $H_p$ is determined by the fields $V_i$,
$\mathsf p_i$ and $B_i$ with $i<p$.
Apply \eqref{cc:uf:hier:eq:mean} to $\bar f_p$, obtaining $B_p,\bar q_p$.
Then put
\begin{equation}\label{cc:uf:hier:eq:supp}
 f_p^S=-\sum_{i=2}^p(B_i\cdot\kappa)\partial_sY_{p+1-i},
\end{equation}
and apply \eqref{cc:uf:hier:eq:trans} to $f_p+f_p^S$, obtaining $V_p,\mathsf p_p$.
Finally set $Q_p=Q[V_p]$ and $C_{p+1}=C[V_p]$.
The term $-(B_p\cdot\kappa)\partial_sV_1$ in
\eqref{cc:uf:hier:eq:supp} accounts for the transport of the leading wave
by the new mean, so $B_p$ must be found before $V_p$. Every term in
\eqref{cc:uf:hier:eq:supp} has zero phase average since $B_i$ and
$\kappa$ are phase independent.

Since $D_aX$, $M$ and $\kappa$ are even in $a$, $R_U\circ X$ is odd in
$a$ and $A_1=g\mathbf aF_h$ changes sign under the reflection
$(a,s)\mapsto(-a,-s)$, the construction is symmetric under this
reflection: every $V_p$, $C_p$ and $B_p$ changes sign under it, and
$W^\natural$ is odd because $X$ is odd.

\begin{proposition}\label{cc:uf:hier:prop:identity}
For every integer \(P\ge2\), the profiles obtained from
\eqref{cc:uf:hier:eq:mean}, \eqref{cc:uf:hier:eq:source} and
\eqref{cc:uf:hier:eq:supp} are well defined, and
\[
 W=\sum_{p=1}^PZ_p+C_{P+1}=\sum_{p=1}^P(V_p+C_{p+1}+B_p),\qquad
 \mathsf P=\sum_{p=1}^P(\mathsf p_p+\bar q_p),\quad \bar q_1=0,
\]
satisfy \(\widetilde {\mathbf d}\cdot W=0\) and
\begin{equation}\label{cc:uf:hier:eq:residual}
 W_t+MW+(W\cdot\widetilde {\mathbf d})W+\widetilde {\mathbf d}\mathsf P
       +R_U\circ X=R_P+R_{\widetilde\chi}.
\end{equation}
Since $\av{V_p}=\av{C_p}=0$, the oscillatory part
\(W_{\rm osc}=W-\av W=\sum_{p=1}^P(V_p+C_{p+1})\) is
supported in \(\supp g\), and so are the sources $f_p+f_p^S$.
The fields $V_p$, $C_p$, $\mathsf p_p$ and $f_p+f_p^S$ vanish for
$t<t_0$.
Set \(V_{P+1}=B_{P+1}=0\), so that \(Y_{P+1}=Z_{P+1}=C_{P+1}\). The
remainder of the expansion is
\begin{align}\label{cc:uf:hier:eq:tail}
 R_P={}&(\partial_t+M)C_{P+1}+\mathbf d\mathsf p_P\\
 &+\sum_{\substack{i,j\le P+1\\i+j\ge P+1}}
                  (Z_i\cdot \mathbf d)Z_j
 +\sum_{\substack{i,j\le P+1\\i+j\ge P+2}}
                  (Z_i\cdot\kappa)\partial_sY_j.\nonumber
\end{align}
At the labels outside $\supp g$ it is independent of the phase.
The construction is jointly smooth.
\end{proposition}
\begin{proof}
By \eqref{cc:uf:hier:eq:slaving} and \eqref{cc:uf:hier:eq:mean},
\[
 \widetilde{\mathbf d}\cdot W
 =\sum_{p=1}^P\bigl[\widetilde{\mathbf d}\cdot
   (\widetilde{\mathbf d}\times Q_p)+\mathbf d\cdot B_p\bigr]=0.
\]
With $W=\sum_{i=1}^{P+1}Z_i$, the quadratic term is
\[
 (W\cdot\widetilde{\mathbf d})W
 =\sum_{i,j}(Z_i\cdot\mathbf d)Z_j
  +\sum_{i,j}((C_i+B_i)\cdot\kappa)\partial_sY_j
\]
with $1\le i,j\le P+1$. For $2\le p\le P$, the terms to be
canceled at level $p$ sum to $H_p-f_p^S$, together with $R_U\circ X$
at $p=2$. The new linear terms
\[
 (\partial_t+M)(V_p+B_p)+\kappa\partial_s\mathsf p_p
    +\mathbf d\bar q_p=f_p+f_p^S+\bar f_p
\]
cancel them, since the right side is $-H_p+f_p^S$ for $p\ge3$ and
$-H_2+f_2^S-R_U\circ X$ for $p=2$. For $V_1$, \eqref{cc:uf:hier:eq:leading}
and \eqref{cc:uf:hier:eq:switched} give
\[
 V_{1,t}+MV_1+\kappa\partial_s\mathsf p_1=\widetilde\chi'A_1,\qquad
 (\partial_t+M)C_2+\mathbf d\mathsf p_1
 =\widetilde\chi\bigl[(\partial_t+M)C[A_1]+\mathbf d\mathsf p[A_1]\bigr]
  +\widetilde\chi'C[A_1].
\]
The two terms with $\widetilde\chi'$ sum to $R_{\widetilde\chi}$, and the
remaining term is canceled in $H_2$. The surviving
products have $i+j>P$ for $\mathbf d$ and $i+j-1>P$ for
$\kappa\partial_s$. Together with the two linear terms left at
level $P+1$, these give \eqref{cc:uf:hier:eq:tail}.

The sources at each level are smooth when the preceding profiles are.
The transverse equation is a linear ordinary differential equation at
fixed $(a,s)$ with smooth coefficients, the mean equation has smooth
solutions by Theorem~\ref{cc:uf:mean:inverse}, and taking the curl
preserves smoothness. By induction on the level, the profiles are well
defined and jointly smooth.

For the supports, let $\mathcal O$ be the open set of pairs $(t,a)$ with
$t<t_0$ or $a\notin\supp g$. The fields
$V_1=\widetilde\chi g\mathbf a F_h$, $C_2$ and $\mathsf p_1$ vanish on
$\mathcal O$, and so does the first term of $H_2$. If $V_i$,
$C_{i+1}$ and $\mathsf p_i$ vanish on $\mathcal O$ for every $i<p$, then
on $\mathcal O$ the field $H_p$ contains only the phase-independent
products $(B_i\cdot\mathbf d)B_j$, and $f_p=f_p^S=0$ there. The solution
of \eqref{cc:uf:hier:eq:trans} then vanishes on $\mathcal O$, by the two
properties recorded after \eqref{cc:uf:hier:eq:trans-ode}, and so do
$C_{p+1}=\mathbf d\times Q[V_p]$ and, by
\eqref{cc:uf:hier:eq:pressure-from-trans}, $\mathsf p_p$, the zero-mean
primitive of a function that vanishes on $\mathcal O$. This proves the
assertions on the supports. On $\mathcal O$ the remainder \eqref{cc:uf:hier:eq:tail}
therefore consists of the products $(B_i\cdot\mathbf d)B_j$, which are
independent of the phase. The means $B_p$ do not vanish before $t_0$
in general.
\end{proof}

\subsubsection{Bounds at every order}

To bound the linear terms in \eqref{cc:uf:hier:eq:H2} and
\eqref{cc:uf:hier:eq:Hp}, we use the transverse equation to evaluate
the time derivative of the curl correction.
\begin{lemma}\label{cc:uf:hier:lem:relay}
There is an absolute exponent $C$ such that, for every \(n\ge0\) and
every field $V$ with $\av V=0$ and $\kappa\cdot V=0$,
\begin{equation}\label{cc:uf:hier:eq:completion-bound}
 \|C[V]\|_{[n],\gamma}\le C_n\mathsf D^C\lambda^{-1}\|V\|_{[n+1],\gamma}.
\end{equation}
If moreover \(V_t+MV+\kappa\partial_s\mathsf p=f\) with $\av f=0$ and
$\av{\mathsf p}=0$, then
\begin{equation}\label{cc:uf:hier:eq:relay-bound}
 \|(\partial_t+M)C[V]+\mathbf d\mathsf p\|_{[n],\gamma}
 \le C_n\mathsf D^C\lambda^{-1}
             \bigl(\|V\|_{[n+1],\gamma}+\|f\|_{[n+1],\gamma}\bigr).
\end{equation}
For the leading oscillation $V=A_1$ on $[t_0,T]$, the source in this estimate is
$f=0$.
\end{lemma}
\begin{proof}
By \eqref{cc:uf:hier:eq:completion}, $Q[V]$ is $\lambda^{-1}$ times the
phase primitive of $V$, multiplied by the coefficients $\kappa/\lambda$
and $\lambda^2/|\kappa|^2$, and $C[V]=\mathbf d\times Q[V]$.
Since $\ell_*^{-1}$ is a power of $\mathsf D$,
\eqref{cc:uf:trunc:eq:product}, \eqref{cc:uf:trunc:eq:deriv},
\eqref{cc:uf:trunc:eq:phase-ops} and
\eqref{cc:uf:hier:eq:coefficient-bound} give
\eqref{cc:uf:hier:eq:completion-bound}.

For \eqref{cc:uf:hier:eq:relay-bound} use \(\kappa_t=-M^T\kappa\),
\((|\kappa|^{-2})_t=2(\kappa\cdot M\kappa)|\kappa|^{-4}\),
and \([\partial_t,\mathbf d_j]=-M_{\ell j}\mathbf d_\ell\).
Differentiating $Q[V]$ in time and inserting the equation for $V_t$, in
which the term \(\kappa\times\kappa\,\partial_s\mathsf p\) vanishes, gives
\begin{equation}\label{cc:uf:hier:eq:potential-time}
 \partial_tQ[V]=\partial_s^{-1}\left[
  \frac{(M^T\kappa)\times V+\kappa\times MV-\kappa\times f}
       {|\kappa|^2}
  -2\frac{(\kappa\times V)(\kappa\cdot M\kappa)}
             {|\kappa|^4}\right].
\end{equation}
Commuting $\partial_t$ past the curl gives
$(\partial_t+M)C[V]=\mathbf d\times\partial_tQ[V]
+[\partial_t,\mathbf d\times]Q[V]+MC[V]$.
By \eqref{cc:uf:hier:eq:potential-time} and the formula for
$\partial_s\mathsf p$ in \eqref{cc:uf:hier:eq:pressure-from-trans},
each term, including $\mathbf d\mathsf p$, carries a factor $\lambda^{-1}$
and at most one label derivative. The product and primitive estimates
therefore give \eqref{cc:uf:hier:eq:relay-bound} without a time derivative
of $M$, $f$ or a mean. The weight $\gamma$ is applied after the time
differentiation.
\end{proof}

Up to fixed powers of $\mathsf D$, the leading oscillation has size
$\lambda^{-5/4}$ by \eqref{cc:uf:hier:eq:size}, its divergence correction
has size $\lambda^{-9/4}$ by \eqref{cc:uf:hier:eq:completion-bound},
and the first mean has size $\lambda^{-5/2}$ by the calculation after
\eqref{cc:uf:hier:eq:H2} and \eqref{cc:uf:mean:bound}.
A label derivative costs a power of $\mathsf D$, while a phase derivative
in Eulerian advection costs $|\kappa|h^{-1}\le\mathsf D^2\lambda^{3/2}$.
The smallest gain in frequency comes from advection of the phase by the divergence
correction, which contributes $\lambda^{-9/4}\lambda^{3/2}=\lambda^{-3/4}$.
Advection of the phase by the mean gives $\lambda^{-1}$,
as does taking a divergence correction. Advection in the labels gives
$\lambda^{-5/4}$. Allowing $\lambda^{-2/3}$ per level leaves a factor
$\lambda^{-1/12}$ to absorb the parent bounds.

We induct on the level at all derivative orders at once. When a source
requires a higher derivative of a previously computed profile, that
order is already available. The new mean $B_p$ enters the transverse
source \eqref{cc:uf:hier:eq:supp} undifferentiated.

\Needspace{10\baselineskip}
\begin{theorem}\label{cc:uf:hier:thm:depth}
There is an absolute exponent $C$ such that the following holds if
$\mathsf D^C\lambda^{-1/12}\le1$. For all integers \(n\ge0\) and
\(p\ge2\), the construction satisfies
\begin{align}
 \|V_p\|_{[n],\gamma}+\|f_p+f_p^S\|_{[n],\gamma}
        +\lambda\|\mathsf p_p\|_{[n],\gamma}
 &\le C_{n,p}\mathsf D^C\lambda^{-5/4}\lambda^{-2(p-1)/3},
 \label{cc:uf:hier:eq:Vbound}\\
 \|B_p\|_{[n]}+\|\bar f_p\|_{[n]}
 &\le C_{n,p}\mathsf D^C\lambda^{-5/2}\lambda^{-2(p-2)/3},
 \label{cc:uf:hier:eq:Bbound}\\
 \|C_{p+1}\|_{[n],\gamma}
 &\le C_{n,p}\mathsf D^C\lambda^{-9/4}\lambda^{-2(p-1)/3}.
 \label{cc:uf:hier:eq:Cbound}
\end{align}
The bounds \eqref{cc:uf:hier:eq:Vbound} and
\eqref{cc:uf:hier:eq:Cbound} also hold for $p=1$, with $f_1=f_1^S=0$.
For all integers \(P\ge2\) and \(m\ge0\), the expansion of
Proposition~\ref{cc:uf:hier:prop:identity} at the depth $P$ satisfies
\begin{equation}\label{cc:uf:hier:eq:sum}
 \|W_{\rm osc}\|_{[m],\gamma}\le C_{m,P}\mathsf D^C\lambda^{-5/4},\qquad
 \Big\|\sum_{p=2}^PB_p\Big\|_{[m]}+\Big\|\sum_{p=2}^P\bar f_p\Big\|_{[m]}
 \le C_{m,P}\mathsf D^C\lambda^{-5/2},
\end{equation}
and its remainder satisfies
\begin{equation}\label{cc:uf:hier:eq:tail-est}
 \|R_P\|_{[m]}\le C_{m,P}\lambda^{-5/4}\lambda^{-2P/3}.
\end{equation}
\end{theorem}
\begin{proof}
By \eqref{cc:uf:trunc:eq:product}, \eqref{cc:uf:trunc:eq:deriv} and
\eqref{cc:uf:hier:eq:coefficient-bound}, since $\ell_*^{-1}$ is a power
of $\mathsf D$ and $\kappa\partial_s$ costs
$|\kappa|h^{-1}\le\mathsf D^2\lambda^{3/2}$ by
\eqref{cc:uf:tower:eq:wave-vector}, fields $Z$, $Z'$, $Y$
satisfy, with the weight on either factor of a product or on neither,
\begin{equation}\label{cc:uf:hier:eq:advection}
 \begin{aligned}
 \|(Z\cdot\mathbf d)Z'\|_{[n],\gamma}+\|(\mathbf d\cdot Z)Z'\|_{[n],\gamma}
 &\le C_n\mathsf D^C\|Z\|_{[n+1]}\|Z'\|_{[n+1],\gamma},\\
 \|(Z\cdot\kappa)\partial_sY\|_{[n],\gamma}
 &\le C_n\mathsf D^C\lambda^{3/2}\|Z\|_{[n]}\|Y\|_{[n+1],\gamma}.
 \end{aligned}
\end{equation}
For $p=1$, \eqref{cc:uf:hier:eq:size} and
\eqref{cc:uf:hier:eq:completion-bound} give
\eqref{cc:uf:hier:eq:Vbound} and \eqref{cc:uf:hier:eq:Cbound} for
$A_1$, $\mathsf p[A_1]$, $C[A_1]$ on $[t_0,T]$, and
$0\le\widetilde\chi\le1$ gives them for $V_1$, $\mathsf p_1$, $C_2$.

To keep powers of $\mathsf D$ from accumulating with the depth, we
propagate a bound for the transverse source and the mean that contains
only powers of $\lambda$. For $p\ge2$ the induction claim at every order is
\begin{equation}\label{cc:uf:hier:eq:induction}
 \|f_p+f_p^S\|_{[n],\gamma}+\lambda^{1/2}\|B_p\|_{[n]}
 \le C_{n,p}\lambda^{-5/4}\lambda^{-2(p-1)/3}.
\end{equation}
For $B_p$ this is \eqref{cc:uf:hier:eq:Bbound} with $\lambda^{1/12}$ in
place of $\mathsf D^C$. Since \eqref{cc:uf:hier:eq:induction} contains
no power of $\mathsf D$, every exponent of $\mathsf D$ below is a sum
of a bounded number of the exponents of \eqref{cc:uf:mean:bound},
\eqref{cc:uf:hier:eq:trans-bound}, \eqref{cc:uf:hier:eq:size},
\eqref{cc:uf:hier:eq:completion-bound},
\eqref{cc:uf:hier:eq:relay-bound} and \eqref{cc:uf:hier:eq:advection},
and it depends neither on the level nor on the depth. The theorem holds
with the largest of these exponents, since $\mathsf D\ge2$. In the rest
of the proof every bound holds up to a factor $C_{n,p}$ or $C_{m,P}$,
which we do not write.

The bound \eqref{cc:uf:hier:eq:induction} at a level $i\ge2$ gives
\eqref{cc:uf:hier:eq:Vbound} and \eqref{cc:uf:hier:eq:Cbound} for $V_i$,
$\mathsf p_i$ and $C_{i+1}$. Indeed, by Proposition~\ref{cc:uf:hier:prop:identity} the source
$f_i+f_i^S$ vanishes for $t<t_0$ and outside $\supp g$, so
\eqref{cc:uf:hier:eq:trans-bound} bounds $V_i$ and $\lambda\mathsf p_i$
by $\mathsf D^C$ times the right side of
\eqref{cc:uf:hier:eq:induction}, and
\eqref{cc:uf:hier:eq:completion-bound} then bounds $C_{i+1}$. For the bounds
at each level it remains to prove \eqref{cc:uf:hier:eq:induction} and
\eqref{cc:uf:hier:eq:Bbound}.

Let $p\ge2$, and suppose that \eqref{cc:uf:hier:eq:induction} holds at
the levels $2\le i<p$ and at every order. Then
\eqref{cc:uf:hier:eq:Vbound} and \eqref{cc:uf:hier:eq:Cbound} hold for
$V_i$ and $C_{i+1}$, $1\le i<p$. Recall $Y_i=V_i+C_i$ and $Z_i=Y_i+B_i$.
The right side of \eqref{cc:uf:hier:eq:Cbound} for $C_i$ and the
bound \eqref{cc:uf:hier:eq:induction} for $B_i$ are at most the right
side of \eqref{cc:uf:hier:eq:Vbound} for $V_i$. Since
$\gamma\le1$, it follows that, for $i<p$,
\begin{equation}\label{cc:uf:hier:eq:YZ}
 \|Y_i\|_{[n],\gamma}+\|Z_i\|_{[n]}
 \le\mathsf D^C\lambda^{-5/4}\lambda^{-2(i-1)/3}.
\end{equation}

The source $f_p$ in \eqref{cc:uf:hier:eq:source} uses only the
fields $V_i$, $\mathsf p_i$ and $B_i$ with $i<p$, and the corrections
$C_i$ with $i\le p$. Subtracting the phase mean of $(Z_i\cdot\mathbf d)Z_j$
removes $(B_i\cdot\mathbf d)B_j$, and every other product has a
factor $Y_i$ or $Y_j$, to which we apply the $\gamma$ norm.
By \eqref{cc:uf:hier:eq:relay-bound} and \eqref{cc:uf:hier:eq:Vbound}
for $V_{p-1}$, the linear term of $H_p$ is at most
$\mathsf D^C\lambda^{-9/4}\lambda^{-2(p-2)/3}$; for $p=2$
this is the bound for $A_1$, with $f=0$ and $\widetilde\chi\le1$. By
\eqref{cc:uf:hier:eq:advection} and \eqref{cc:uf:hier:eq:YZ}, a product
$(Z_i\cdot\mathbf d)Z_j$ with $i+j=p$ is at most
$\mathsf D^C\lambda^{-5/2}\lambda^{-2(p-2)/3}$. By
\eqref{cc:uf:hier:eq:advection}, \eqref{cc:uf:hier:eq:YZ} and
\eqref{cc:uf:hier:eq:Cbound} for $C_i$, a product
$(C_i\cdot\kappa)\partial_sY_j$ with $i+j=p+1$ is at most
$\mathsf D^C\lambda^{3/2}\lambda^{-9/4}\lambda^{-5/4}\lambda^{-2(p-2)/3}
=\mathsf D^C\lambda^{-2}\lambda^{-2(p-2)/3}$. The right side of
\eqref{cc:uf:hier:eq:induction} is
$\lambda^{-23/12}\lambda^{-2(p-2)/3}$, so the three bounds are this
number times
\[
 \mathsf D^C\lambda^{-1/3},\qquad
 \mathsf D^C\lambda^{-7/12},\qquad
 \mathsf D^C\lambda^{-1/12},
\]
respectively. Thus $\|f_p\|_{[n],\gamma}$ is at most
$\mathsf D^C\lambda^{-1/12}$ times the right side of
\eqref{cc:uf:hier:eq:induction}.

In the mean source $\bar f_p$, integration by parts in $s$ removes the
fast derivative just as at level two. By \eqref{cc:uf:hier:eq:slaving},
for $i\ge2$ we have
\begin{equation}\label{cc:uf:hier:eq:fast-mean}
 \av{(C_i\cdot\kappa)\partial_sY_j}
                 =\av{(\mathbf d\cdot V_{i-1})Y_j}.
\end{equation}
The linear terms have zero phase average because $C_p$ and
$\mathsf p_{p-1}$ do. Products of one mean and one oscillation also
average to zero, since $B_i$ and $\kappa$ are phase independent and
$\av{Y_j}=0$. Thus the mean source is
\[
 -\bar f_p=
 \sum_{i+j=p}\av{(Y_i\cdot\mathbf d)Y_j}
 +\sum_{i+j=p}\av{(\mathbf d\cdot V_i)Y_j}
 +\sum_{i+j=p}(B_i\cdot\mathbf d)B_j,
\]
with the additional term $R_U\circ X$ at $p=2$. By
\eqref{cc:uf:hier:eq:advection} and \eqref{cc:uf:hier:eq:YZ}, the first
two sums are at most $\mathsf D^C\lambda^{-5/2}\lambda^{-2(p-2)/3}$.
By \eqref{cc:uf:hier:eq:advection} and
\eqref{cc:uf:hier:eq:induction} at the lower levels, the third is at
most $\mathsf D^C\lambda^{-7/2}\lambda^{-2(p-2)/3}$, which is smaller.
The term $R_U\circ X$ is at most $C_n\lambda^{-5/2}$, by
\eqref{cc:uf:hier:eq:parent-residual}. Hence
$\|\bar f_p\|_{[n]}\le\mathsf D^C\lambda^{-5/2}\lambda^{-2(p-2)/3}$
at every order. With this estimate at the order $n+1$,
\eqref{cc:uf:mean:bound} bounds $\|B_p\|_{[n]}$ by $\mathsf D^C$ times
as much. This proves \eqref{cc:uf:hier:eq:Bbound} at the level $p$, and
with the hypothesis $\mathsf D^C\le\lambda^{1/12}$ it gives
\eqref{cc:uf:hier:eq:induction} for $B_p$.

For $f_p^S$ in \eqref{cc:uf:hier:eq:supp}, the means
$B_i$, $2\le i\le p$, are now bounded by
\eqref{cc:uf:hier:eq:induction}, and $Y_{p+1-i}$ is bounded by
\eqref{cc:uf:hier:eq:YZ}. By \eqref{cc:uf:hier:eq:advection}, each
product $(B_i\cdot\kappa)\partial_sY_{p+1-i}$ is therefore at most
$\mathsf D^C\lambda^{-3/2}\lambda^{-2(p-1)/3}$, which is
$\mathsf D^C\lambda^{-1/4}$ times the right side of
\eqref{cc:uf:hier:eq:induction}. With the bound for $f_p$,
$\|f_p+f_p^S\|_{[n],\gamma}$ is at most $\mathsf D^C\lambda^{-1/12}$
times that right side, and the hypothesis gives
\eqref{cc:uf:hier:eq:induction} at the level $p$. This completes the
induction. The number of products in a level and the constants of the
lower levels at higher orders enter $C_{n,p}$.

Since $\lambda\ge1$, summing
\eqref{cc:uf:hier:eq:Vbound}--\eqref{cc:uf:hier:eq:Cbound} over
$p\le P$ gives \eqref{cc:uf:hier:eq:sum}. For the remainder,
recall $Y_{P+1}=Z_{P+1}=C_{P+1}$. By \eqref{cc:uf:hier:eq:Cbound} for
$C_{P+1}$, the bound \eqref{cc:uf:hier:eq:YZ} holds for $i\le P+1$. The
products in \eqref{cc:uf:hier:eq:tail} are therefore bounded as those
of $H_p$ and $f_p^S$ were: the products $(Z_i\cdot\mathbf d)Z_j$ with
$i+j=q$, and $(C_i\cdot\kappa)\partial_sY_j$ and
$(B_i\cdot\kappa)\partial_sY_j$ with $i+j=q+1$, are at most
$\mathsf D^C\lambda^{-1/12}$ times $\lambda^{-5/4}\lambda^{-2(q-1)/3}$,
and in \eqref{cc:uf:hier:eq:tail} they have $q\ge P+1$. By
\eqref{cc:uf:hier:eq:relay-bound} and \eqref{cc:uf:hier:eq:Vbound} for
$V_P$, the linear term
\((\partial_t+M)C_{P+1}+\mathbf d\mathsf p_P\) is at most
$\mathsf D^C\lambda^{-1/3}$ times $\lambda^{-5/4}\lambda^{-2P/3}$. The
products of two means are bounded in the
unweighted norm, and every other term in the weighted norm, which
dominates the unweighted one because $\gamma\le1$. With the hypothesis,
this is \eqref{cc:uf:hier:eq:tail-est}.
\end{proof}

\subsection{The approximate solution}\label{cc:uf:sec:approximate}
Fix an integer $P\ge2$. The profiles $W,\mathsf P$ of
Proposition~\ref{cc:uf:hier:prop:identity} at depth $P$ give the
Eulerian velocity and pressure
\[
 V=U+W^\natural,\qquad \Pi=p_U+\mathsf P^\natural,\qquad
 W=W_{\rm osc}+\bar W,\quad \bar W=\sum_{p=2}^PB_p .
\]
Fix an absolute integer $b\ge1$, at least as large
as the exponents $C$ of Theorem~\ref{cc:uf:hier:thm:depth} and of
Lemma~\ref{cc:uf:hier:lem:physical}, and sufficiently large for the
finitely many estimates of this subsection and the next. Put
\[
 C_{\rm step}=\Gamma_P\mathsf D^b,
\]
where $\Gamma_P\ge1$ depends only on $P$ and the fixed profiles.
We take it sufficiently large compared with every $C_{2,P}$ in this
subsection and the next. The estimates of the residual, the pressure
Hessian and the flow use the profile bounds through order two. The
truncation uses the velocity bound in $\Holder{P}$ and the tail estimate
\eqref{cc:uf:trunc:eq:tail}. From here on we assume
\begin{equation}\label{cc:uf:step:eq:small}
 C_{\rm step}\lambda^{-1/12}\le1 .
\end{equation}
Since $C_{\rm step}\ge\mathsf D^C$ for the exponent $C$ of
Theorem~\ref{cc:uf:hier:thm:depth}, the hypothesis of that theorem
holds.

To estimate the velocity and residual in Eulerian coordinates, we evaluate
each profile at $a=X(t,\cdot)^{-1}(x)$ and $s=k_*\cdot a$.

\begin{lemma}\label{cc:uf:hier:lem:physical}
Let \(g_0(t,a,s)\) be smooth and independent of the phase at the labels
outside \(\supp g\). Assume its phase average extends smoothly and
periodically in $a$ across the cell faces. There is an absolute exponent
$C$ such that, for every integer \(m\ge0\) and every $t\in[t_0,T]$,
\begin{align}
 \sum_{|\beta|=m}\|\partial_x^\beta g_0^\natural(t)\|_{L^\infty}
 &\le C_m(\mathsf D^C\lambda^{3/2})^m\gamma(t)
                  \|g_0\|_{[m],\gamma},
 \label{cc:uf:hier:eq:physical-infty}\\
 \|g_0^\natural(t)\|_{\Holder{m}}
 &\le C_m(\mathsf D^C\lambda^{3/2})^{m+1}\gamma(t)
                  \|g_0\|_{[m+1],\gamma}.
 \label{cc:uf:hier:eq:physical}
\end{align}
Both bounds also hold on $[0,T]$ with \(\gamma\) replaced by one, and for
a phase-independent field with \(\mathsf D^C\lambda^{3/2}\) replaced by
\(\mathsf D^C\).
\end{lemma}
\begin{proof}
The field \(\widetilde{\mathbf d}_ig_0\) is again independent of the
phase outside \(\supp g\), which lies in the interior of the cell, so
\(g_0^\natural\) and its derivatives extend smoothly and periodically.
The spatial chain rule in \eqref{cc:uf:hier:eq:graph} and
\(\widetilde{\mathbf d}=(D_aX)^{-T}(\nabla_a+k_*\partial_s)\) give
\[
 \sum_i\|\widetilde{\mathbf d}_ig_0\|_{[n],\gamma}
 \le C_n\mathsf D^C\lambda^{3/2}\|g_0\|_{[n+1],\gamma},\qquad n\ge0,
\]
by Lemma~\ref{cc:uf:trunc:lem:fields} and
\eqref{cc:uf:hier:eq:coefficient-bound}: a label derivative costs
$\ell_*^{-1}$, a fixed power of $\mathsf D$, and $k_*\partial_s$ costs
$|k_*|h^{-1}\le\mathsf D\lambda^{3/2}$, by
\eqref{cc:uf:tower:eq:wave-vector}.
Iterate the displayed bound; at each time the supremum of a field is at
most $\gamma(t)$ times its weighted norm of order zero. This proves
\eqref{cc:uf:hier:eq:physical-infty}. The bound
\(\|v\|_{\Holder{0}}\le C(\|v\|_\infty+\|\nabla v\|_\infty)\),
applied to derivatives through order \(m\), proves
\eqref{cc:uf:hier:eq:physical} with one additional derivative. A
phase-independent field has
\(\widetilde{\mathbf d}g_0=\mathbf dg_0\), proving the last assertion.
\end{proof}

In Section~\ref{cc:uf:sec:tower} the oscillations are summed in the
supremum norm, while their higher derivatives control the Fourier tail
in \eqref{cc:uf:tower:eq:tail-small}. The means can extend across the
torus and back to time zero, so their sum must converge at every fixed
derivative order.

\begin{lemma}\label{cc:uf:step:lem:approximate}
The field \(V\) is smooth, odd and
divergence free, \(\Pi\) is smooth, and
\begin{equation}\label{cc:uf:step:eq:approximate}
 V_t+V\cdot\nabla V+\nabla\Pi=R_V,\qquad
 R_V=R_P^\natural+R_{\widetilde\chi}^\natural .
\end{equation}
The oscillation satisfies
\begin{equation}\label{cc:uf:step:eq:physical-sup}
 \|W_{\rm osc}^\natural\|_{C_tL^\infty}
 \le C_{0,P}\mathsf D^b\lambda^{-5/4},
\end{equation}
and for every \(m\ge0\),
\begin{equation}\label{cc:uf:step:eq:physical}
\begin{gathered}
 \|W_{\rm osc}^\natural\|_{C_t\Holder{m}}
 \le C_{m,P}\mathsf D^{b(m+2)}\lambda^{3(m+1)/2-5/4},\qquad
 \|\bar W^\natural\|_{C_t\Holder{m}}
 \le C_{m,P}\mathsf D^{b(m+2)}\lambda^{-5/2},\\
 \|R_V\|_{C_t\Holder{1}}
 \le C_{\rm step}\lambda^2
   \bigl(\lambda^{-2P/3}+\|\widetilde\chi'\gamma\|_\infty\bigr),
\end{gathered}
\end{equation}
where the last supremum is over $[t_0,T]$.
\end{lemma}
\begin{proof}
By \eqref{cc:uf:hier:eq:graph}, at \(x=X(t,a)\) one has
\((\partial_t+U\cdot\nabla)W^\natural=(W_t)^\natural\),
\(W^\natural\cdot\nabla U=(MW)^\natural\),
\(W^\natural\cdot\nabla W^\natural=((W\cdot\widetilde{\mathbf d})W)^\natural\) and
\(\nabla\mathsf P^\natural=(\widetilde{\mathbf d}\mathsf P)^\natural\). Adding
the definition of \(R_U\) and \eqref{cc:uf:hier:eq:residual} gives
\eqref{cc:uf:step:eq:approximate}, because \((R_U\circ X)^\natural=R_U\):
the residual of the parent is canceled exactly. Also
\(\widetilde{\mathbf d}\cdot W=0\) gives \(\diver V=0\).
Since \(\supp g\) lies in the interior of the cell, the oscillatory
terms in \(W^\natural\) and \(\mathsf P^\natural\) extend smoothly and
periodically for every \(k_*\). Proposition~\ref{cc:uf:hier:prop:identity}
gives the smoothness of \(V\) and \(\Pi\), and \(W^\natural\) is odd by
the symmetry of the construction of Subsection~\ref{cc:uf:sec:hierarchy}
under \((a,s)\mapsto(-a,-s)\).

Since $\gamma\le1$, the supremum of
$W_{\rm osc}^\natural$ is at most $\|W_{\rm osc}\|_{[0],\gamma}$, and
\eqref{cc:uf:hier:eq:sum} gives \eqref{cc:uf:step:eq:physical-sup}.
The bound \eqref{cc:uf:hier:eq:physical} of
Lemma~\ref{cc:uf:hier:lem:physical}, applied to $W_{\rm osc}$ and
$\bar W$ with \eqref{cc:uf:hier:eq:sum} at the order $m+1$, proves the
first two bounds of \eqref{cc:uf:step:eq:physical}.
By Proposition~\ref{cc:uf:hier:prop:identity}, $R_P$ is independent of
the phase outside $\supp g$, so \eqref{cc:uf:hier:eq:physical} at $m=1$,
with $\gamma$ replaced by one, and \eqref{cc:uf:hier:eq:tail-est} bound
$R_P^\natural$ in $C_t\Holder{1}$ by
$C_{2,P}\mathsf D^C\lambda^{7/4}\lambda^{-2P/3}$. By
\eqref{cc:uf:hier:eq:size} and \eqref{cc:uf:hier:eq:completion-bound},
$\|A_1+C[A_1]\|_{[2],\gamma}\le C\mathsf D^C\lambda^{-5/4}$. At each time
$t$, $R_{\widetilde\chi}$ is $\widetilde\chi'(t)$ times this field, so
\eqref{cc:uf:hier:eq:physical} at that time bounds
$R_{\widetilde\chi}^\natural(t)$ in $\Holder{1}$ by
$C\mathsf D^C\lambda^{7/4}|\widetilde\chi'(t)|\gamma(t)$. This proves
the last bound of \eqref{cc:uf:step:eq:physical}.
\end{proof}

In \eqref{cc:uf:step:eq:physical} the depth $P$ controls the expansion
remainder. For a later wave, \eqref{cc:uf:tower:eq:switch-small}
bounds $\gamma$ on the interval where the time cutoff varies.

The next mean equation requires an upper bound on the added pressure
Hessian, since the pressure Hessian enters the displacement form
\eqref{cc:uf:mean:jacobi-id} with a minus sign. The next wave also needs
the velocity gradient at the origin. Once $\widetilde\chi=1$, the
rank-one term in \eqref{cc:uf:tower:eq:central-return} supplies the
wave contribution to the gradient in \eqref{cc:uf:tower:eq:split}.

\Needspace{8\baselineskip}
\begin{lemma}\label{cc:uf:step:lem:pressure-center}
On \([0,T]\times\mathbb T^3\),
\begin{equation}\label{cc:uf:step:eq:hessian}
 D_x^2\mathsf P^\natural\le C_{\rm step}\lambda^{-1/4}\Id .
\end{equation}
At \(x=0\), for \(t_0\le t\le T\),
\begin{equation}\label{cc:uf:tower:eq:central-return}
 \nabla V(t,0)=\nabla U(t,0)
   +\widetilde\chi(t)\,\mathbf a(t,0)\otimes\kappa(t,0)+E_c(t),
 \qquad
 |E_c|\le\lambda^{-1/3} .
\end{equation}
\end{lemma}
\begin{proof}
Write the leading pressure as \(\mathsf p_1=-c_{\mathbf a}P_h\), with
\(c_{\mathbf a}=2\widetilde\chi g(\kappa\cdot M\mathbf a)/|\kappa|^2\). Since
$D_x^2\mathsf p_1^\natural
=(\widetilde{\mathbf d}\,\widetilde{\mathbf d}\,\mathsf p_1)^\natural$ and
\(\mathbf d_i\kappa_j=\mathbf d_j\kappa_i\), with \(a=X(t,\cdot)^{-1}(x)\)
and \(s=k_*\cdot a\) on the right,
\begin{equation}\label{cc:uf:exact:eq:primary-hess}
 D_x^2\mathsf p_1^\natural=
 -c_{\mathbf a}F_h'\kappa\otimes\kappa
 +(\mathbf d\partial_s\mathsf p_1)\otimes\kappa+
 \kappa\otimes(\mathbf d\partial_s\mathsf p_1)
 +(\partial_s\mathsf p_1)\mathbf d\kappa+\mathbf d^2\mathsf p_1 .
\end{equation}
The product $\widetilde\chi g$ takes values in $[0,1]$ and vanishes off
$[t_0,T]\times\supp g$, where $\kappa\cdot M\mathbf a>0$ by
\eqref{cc:uf:tower:eq:near-origin}. Hence
\(0\le c_{\mathbf a}|\kappa|^2\le2|\kappa||M||\mathbf a|\). With
\(F_h'\ge-Ch\), the first tensor on the right is therefore at most
\(Chc_{\mathbf a}|\kappa|^2\Id\). On $[t_0,T]\times\supp g$,
\eqref{cc:uf:tower:eq:near-origin} gives
$|\mathbf a|\le2\lambda^{-3/4}$, while $|\kappa|\le\mathsf D^2\lambda$
and $|M|\le\mathsf D$ by \eqref{cc:uf:tower:eq:wave-vector} and
\eqref{cc:uf:tower:eq:parent-bound}. So
\(c_{\mathbf a}|\kappa|^2\le4\mathsf D^3\lambda^{1/4}\) and, since
$h=\lambda^{-1/2}$, the first tensor is at most
\(C\mathsf D^3\lambda^{-1/4}\Id\). By
\eqref{cc:uf:trunc:eq:product}, \eqref{cc:uf:trunc:eq:deriv} and
\eqref{cc:uf:hier:eq:coefficient-bound}, the sum of the remaining
tensors is at most
$C\mathsf D^C\lambda^{3/2}\|\mathsf p_1\|_{[2]}$ in supremum norm, and
\(\lambda\|\mathsf p_1\|_{[2]}\le C\mathsf D^C\lambda^{-5/4}\)
by \eqref{cc:uf:hier:eq:size}. So their sum is at most
$C\mathsf D^C\lambda^{-3/4}$.

For the higher oscillatory pressures, \eqref{cc:uf:hier:eq:Vbound}
gives the least decay in $\lambda$ at $p=2$, where it bounds
$\mathsf p_2$ by $C_{2,P}\mathsf D^C\lambda^{-1-23/12}$.
Two Eulerian derivatives in \eqref{cc:uf:hier:eq:physical-infty}
would cost $\lambda^3$, leaving $\lambda^{1/12}$. To bound the first
derivative, use the transverse equation \eqref{cc:uf:hier:eq:pressure-from-trans}.
For $2\le p\le P$ it gives
\[
 \nabla_x\mathsf p_p^\natural
   =(\mathbf d\mathsf p_p+\kappa\partial_s\mathsf p_p)^\natural,
 \qquad
 \kappa\partial_s\mathsf p_p=P_\kappa(f_p+f_p^S-2MV_p).
\]
The label derivative costs a power of $\mathsf D$, and the phase term
is bounded by the source and velocity in \eqref{cc:uf:hier:eq:Vbound}.
The coefficient estimates therefore give
\[
 \Big\|\widetilde{\mathbf d}\sum_{p=2}^P\mathsf p_p\Big\|_{[1],\gamma}
 \le C_{2,P}\mathsf D^C\lambda^{-23/12}.
\]
The oscillatory pressures are supported in $\supp g$ in Lagrangian coordinates.
One further Eulerian derivative costs $\lambda^{3/2}$ up to a power
of $\mathsf D$, so \eqref{cc:uf:hier:eq:physical-infty} gives
\begin{equation}\label{cc:uf:exact:eq:higher-hess}
 \Big\|D_x^2\sum_{p=2}^P\mathsf p_p^\natural\Big\|_\infty
 \le C_{2,P}\mathsf D^C\lambda^{-5/12}.
\end{equation}
For the mean pressures, add the mean equations
\eqref{cc:uf:hier:eq:mean} of the levels $2\le p\le P$. In Eulerian coordinates the
sum is the linearized Euler equation \eqref{cc:uf:mean:physical} for
$\bar W^\natural$ over $U$, with the pressure
$\sum_p\bar q_p^\natural$ and the source $(\sum_p\bar f_p)^\natural$,
and its divergence is
\[
 \Delta\sum_{p=2}^P\bar q_p^\natural
 =\diver\Bigl(\sum_{p=2}^P\bar f_p\Bigr)^\natural
  -2\tr\bigl(\nabla U\,\nabla\bar W^\natural\bigr).
\]
The operator $D^2\Delta^{-1}$ is bounded on $\Holder{0}$. By
\eqref{cc:uf:hier:eq:physical} for phase-independent fields and
\eqref{cc:uf:hier:eq:sum}, the source and $\bar W^\natural$ are at most
$C_{2,P}\mathsf D^C\lambda^{-5/2}$ in $C_t\Holder{1}$, and
$\|\nabla U\|_{C_t\Holder{0}}\le C\mathsf D\norm{\nabla U}_\infty
\le C\mathsf D^2$ by \eqref{cc:uf:trunc:eq:lp} and
\eqref{cc:uf:tower:eq:parent-bound}. Hence the Hessian of the
sum of the mean pressures is at most
$C_{2,P}\mathsf D^C\lambda^{-5/2}$ in supremum norm. Together with the
bounds for the tensors of \eqref{cc:uf:exact:eq:primary-hess} and with
\eqref{cc:uf:exact:eq:higher-hess}, this proves
\eqref{cc:uf:step:eq:hessian}.

To compute the velocity gradient at $x=0$, evaluate at $a=0$ and $s=0$,
since the flow of $U$ fixes the origin. As
\(g(0)=F_h'(0)=1\) and \(F_h(0)=0\), the chain rule
\eqref{cc:uf:hier:eq:graph} gives
\(D_xV_1^\natural(t,0)=\widetilde\chi(t)\mathbf a(t,0)\otimes\kappa(t,0)\).
Hence
\[
 E_c(t)=D_x(W_{\rm osc}-V_1)^\natural(t,0)+D_x\bar W^\natural(t,0).
\]
The field \(W_{\rm osc}-V_1=\sum_{p=2}^PV_p+\sum_{p=2}^{P+1}C_p\) is supported in
\(\supp g\). By \eqref{cc:uf:hier:eq:Vbound} and
\eqref{cc:uf:hier:eq:Cbound}, its norm \(\|\cdot\|_{[1],\gamma}\) is at
most \(C_{2,P}\mathsf D^C\lambda^{-23/12}\), the bound for $V_2$,
and \(\|\bar W\|_{[1]}\le C_{2,P}\mathsf D^C\lambda^{-5/2}\) by
\eqref{cc:uf:hier:eq:sum}. The bound \eqref{cc:uf:hier:eq:physical-infty}
at \(m=1\) therefore gives
\(|E_c|\le C_{2,P}\mathsf D^C\lambda^{-5/12}
=C_{2,P}\mathsf D^C\lambda^{-1/12}\cdot\lambda^{-1/3}\). Since
$b$ is at least the exponent $C$ of this bound and
$\Gamma_P\ge C_{2,P}$, \eqref{cc:uf:step:eq:small} makes it at
most $\lambda^{-1/3}$.
This proves \eqref{cc:uf:tower:eq:central-return}.
\end{proof}

To allow for the difference introduced by frequency truncation, let $Y$
be the flow of $V+v$ for a small divergence-free field $v$. A direct Gronwall
estimate for $D_aY$ puts the large velocity gradient in an exponential.
Write $Y(t,a)=X(t,\mathcal A_t(a))$ and let the label and phase evolve
as independent variables. Incompressibility expresses the oscillatory
phase velocity through a label divergence, giving a small velocity
gradient in the scaled variables. We then differentiate the change of
label on the graph $s=k_*\cdot a$ to recover $D\mathcal A_t$.

\begin{lemma}\label{cc:uf:step:lem:flow-bounds}
Let $v$ be a smooth divergence-free field on $[0,T]\times\mathbb T^3$
with $\norm v_{C_tC^2}\le\lambda^{-5/2}$. The flow $Y$ of $V+v$ from
time zero satisfies
\begin{equation}\label{cc:uf:step:eq:deformation}
 \sup|(D_aY)^{\pm1}|\le C_{\rm step}\lambda^{1/4}.
\end{equation}
\end{lemma}
\begin{proof}
The field $v\circ X$ is independent of the phase, and
$\mathbf d\cdot(v\circ X)=(\diver v)\circ X=0$, so
$\widetilde{\mathbf d}\cdot(W+v\circ X)=0$. By the chain rule and
\eqref{cc:uf:hier:eq:coefficient-bound} for $D_aX$,
$\|v\circ X\|_{[2]}\le C\mathsf D^C\norm v_{C_tC^2}$.
Put $c=(D_aX)^{-1}(W+v\circ X)$ and write $c=c_{\rm osc}+\bar c$ for its
oscillatory and phase-independent parts, with $c_{\rm osc}$ repeated over
the lattice by $(a,s)\mapsto(a+n,s+k_*\cdot n)$,
$n\in\mathbb Z^3$. On the graph $s=k_*\cdot a$, the field $c$ is then
$(D_aX)^{-1}((W^\natural+v)\circ X)$. The Piola identity turns the
solenoidality of $W+v\circ X$ into $(\nabla_a+k_*\partial_s)\cdot c=0$,
and the zero-mean phase primitive of its oscillatory part gives
$k_*\cdot c_{\rm osc}=-\partial_s^{-1}\diver_a c_{\rm osc}$. With
\eqref{cc:uf:hier:eq:sum} and \eqref{cc:uf:hier:eq:coefficient-bound},
the norms without weight therefore satisfy
\[
 \|c_{\rm osc}\|_{[2]}\le C_{2,P}\mathsf D^C\lambda^{-5/4},\qquad
 \|\bar c\|_{[2]}\le C_{2,P}\mathsf D^C\lambda^{-5/2},\qquad
 \|k_*\cdot c\|_{[1]}\le C_{2,P}\mathsf D^C\lambda^{-5/4};
\]
in the last bound the oscillatory part costs one label derivative and
the mean part the factor $|k_*|\le\mathsf D\lambda$ of
\eqref{cc:uf:tower:eq:wave-vector}.

Let \(\mathscr L_t\), with first component \(a+\eta_a(t,a,s)\), be the
flow from time zero of
\begin{equation}\label{cc:uf:exact:eq:particle}
 \dot a=c(t,a,s),\qquad \dot s=k_*\cdot c(t,a,s).
\end{equation}
In the coordinates $(a/\ell_*,s/h)$ this vector field is
$(c/\ell_*,k_*\cdot c/h)$, and its gradient is at most
$C(\ell_*^{-1}\|c\|_{[1]}+h^{-1}\|k_*\cdot c\|_{[1]})
\le C_{2,P}\mathsf D^C\lambda^{-3/4}
\le C_{2,P}\mathsf D^C\lambda^{-1/12}$. Since $b$ is at least the
exponent $C$ of this bound, \eqref{cc:uf:step:eq:small} makes it at
most $C_{2,P}/\Gamma_P$, a small absolute constant since $\Gamma_P$ is
large.
Since $t\le T\le1$, Gronwall's inequality then bounds the
gradients of $\mathscr L_t$ and of $\mathscr L_t^{-1}$ in these
coordinates by $2$. The scaled gradient of $c$ is at most
$\|c\|_{[1]}$, so integrating $\dot\eta_a=c$ along the flow gives
\[
 \ell_*|\partial_a\eta_a|+h|\partial_s\eta_a|
 \le C\|c\|_{[1]}\le C_{2,P}\mathsf D^C\lambda^{-5/4},
\]
and the displacement of $\mathscr L_t^{-1}$ satisfies the same bound.

Differentiating $Y=X(t,\mathcal A)$ gives
$\partial_t\mathcal A=((D_aX)^{-1}((W^\natural+v)\circ X))(t,\mathcal A)$.
Since $s-k_*\cdot a$ is preserved by \eqref{cc:uf:exact:eq:particle},
the graph $s=k_*\cdot a$ is invariant, so
$\mathcal A_t(a)=a+\eta_a(t,a,k_*\cdot a)$, the first
component of \(\mathscr L_t\) on this graph. The same formula with the
displacement of \(\mathscr L_t^{-1}\) gives the inverse map.
Thus $D\mathcal A_t-\Id=\partial_a\eta_a+\partial_s\eta_a\otimes k_*$,
and with the corresponding identity for the inverse,
\[
 |D\mathcal A_t^{\pm1}-\Id|
 \le\bigl(\ell_*^{-1}+|k_*|h^{-1}\bigr)C_{2,P}\mathsf D^C\lambda^{-5/4}
 \le C_{2,P}\mathsf D^C\lambda^{1/4}.
\]
Then \(D_aY=(D_aX\circ\mathcal A)D\mathcal A\) and
\((D_aY)^{-1}=(D\mathcal A)^{-1}((D_aX)^{-1}\circ\mathcal A)\), with
the bound \eqref{cc:uf:tower:eq:parent-bound} for the deformation of
the parent, give \eqref{cc:uf:step:eq:deformation}.
\end{proof}

\Needspace{11\baselineskip}
\subsection{Frequency truncation}\label{cc:uf:trunc:sec:truncation}
We truncate $V$ with the smooth frequency cutoff $S_\Lambda$ of
Appendix~\ref{cc:app:cutoff}, as in the mollification step
of~\cite[Section~2.4]{BDLSVOnsager}. The cutoff bounds all higher
derivatives of the next parent in terms of $\Lambda$ and its supremum.
We choose its pressure $p$ so that its residual is divergence free.
The equation for the expansion pressure $\Pi$ contains $\diver R_V$
by \eqref{cc:uf:step:eq:approximate}, so the bound for $D^2(p-\Pi)$
includes both the residual of $V$ and the Fourier tail. The iteration
uses the gradient bound on the new residual and the upper bound on
$D^2(p-p_U)$ to propagate \eqref{cc:uf:tower:eq:parent-coercive}.

\begin{proposition}\label{cc:uf:step:thm:step}
Let $\Lambda\ge1$, and with the depth $P$ of the expansion put
\begin{equation}\label{cc:uf:step:eq:truncation}
 u=S_\Lambda V,\qquad v_\Lambda=u-V,\qquad
 \varepsilon_\Lambda=\Lambda^{2-P}
 \bigl(1+\norm V_{C_t\Holder{P}}\bigr)^2 .
\end{equation}
Let $p$ be the function of zero mean with
$\Delta p=-\tr((\nabla u)^2)$, and let
$R_u=u_t+u\cdot\nabla u+\nabla p$. Then $u$ is smooth, odd, divergence
free and supported on frequencies at most $2\Lambda$, its residual $R_u$ is odd,
divergence free and supported on frequencies at most $4\Lambda$, and
\begin{equation}\label{cc:uf:step:eq:truncated}
 \norm{v_\Lambda}_{C_t\Holder{2}}\le C_{2,P}\varepsilon_\Lambda,\qquad
 \norm{R_u}_{C_t\Holder{1}}+\norm{D^2(p-\Pi)}_{C_t\Holder{0}}
 \le C\norm{R_V}_{C_t\Holder{1}}+C_{2,P}\varepsilon_\Lambda .
\end{equation}
Moreover, on $[0,T]\times\mathbb T^3$,
\begin{equation}\label{cc:uf:step:eq:pressure-error}
 D^2(p-p_U)\le\bigl(C_{\rm step}\lambda^{-1/4}
   +C\norm{R_V}_{C_t\Holder{1}}+C_{2,P}\varepsilon_\Lambda\bigr)\Id .
\end{equation}
\end{proposition}
\begin{proof}
By Lemma~\ref{cc:uf:trunc:lem:cutoff}, $u$ and $v_\Lambda$ are odd and
divergence free and $u$ is supported on frequencies at most $2\Lambda$.
Since $\diver(u\cdot\nabla u)=\tr((\nabla u)^2)$, the residual is
$R_u=\partial_tu+\PP\diver(u\otimes u)$; it is odd, divergence free and
supported on frequencies at most $4\Lambda$.
Applying $S_\Lambda$ and then $\PP$ to \eqref{cc:uf:step:eq:approximate}
gives
$\partial_tu+\PP\diver S_\Lambda(V\otimes V)=\PP S_\Lambda R_V$.
Subtracting this from the formula for $R_u$ gives
\begin{equation}\label{cc:uf:trunc:eq:residual}
 R_u=\PP\bigl[S_\Lambda R_V
   +\diver\bigl(v_\Lambda\otimes u+V\otimes v_\Lambda
                +(\Id-S_\Lambda)(V\otimes V)\bigr)\bigr].
\end{equation}
To bound $R_u$ in $\Holder{1}$ we need the three quadratic tensors in
\eqref{cc:uf:trunc:eq:residual} in $\Holder{2}$. This accounts for
the power $\Lambda^{2-P}$ and the square in $\varepsilon_\Lambda$.
The tail bound \eqref{cc:uf:trunc:eq:tail} at the order $P$ gives
$\norm{v_\Lambda}_{\Holder{s}}
\le C_{2,P}\Lambda^{s-P}\norm V_{\Holder{P}}$ for
$s\le2$, and in particular the first bound of
\eqref{cc:uf:step:eq:truncated}. By Leibniz' rule and the algebra
inequality for the norm of $\Holder{0}$, for every integer $m\ge0$,
\begin{equation}\label{cc:holder:product-growth}
 \norm{fg}_{\Holder{m}}\le C_m\norm f_{\Holder{m}}\norm g_{\Holder{m}}.
\end{equation}
Hence the three tensors on the right of
\eqref{cc:uf:trunc:eq:residual} are at most
$C_{2,P}\Lambda^{2-P}\norm V_{\Holder{P}}^2\le C_{2,P}\varepsilon_\Lambda$ in
$\Holder{2}$. The periodic
multipliers $\PP$ and $D^2\Delta^{-1}$ of order zero are bounded
on $\Holder{0}$ and commute with derivatives, as does $S_\Lambda$; this
bounds $R_u$ in $\Holder{1}$.
For the pressure, the divergence of \eqref{cc:uf:step:eq:approximate} is
$\Delta\Pi=\diver R_V-\tr((\nabla V)^2)$. Hence
\[
 \Delta(p-\Pi)
 =-\tr(\nabla v_\Lambda\nabla V)-\tr(\nabla u\nabla v_\Lambda)
  -\diver R_V,
\]
and $\norm{D^2(p-\Pi)}_{\Holder{0}}
\le C\norm{v_\Lambda}_{\Holder{1}}(\norm V_{\Holder{1}}+\norm u_{\Holder{1}})
+C\norm{R_V}_{\Holder{1}}$. This proves \eqref{cc:uf:step:eq:truncated}.
Since $p-p_U=(p-\Pi)+\mathsf P^\natural$,
Lemma~\ref{cc:uf:step:lem:pressure-center} then gives
\eqref{cc:uf:step:eq:pressure-error}.
\end{proof}

\section{The initial wave and the inductive step}\label{cc:sec:iteration}
The bounds \eqref{cc:uf:tower:eq:wave-origin} in
Section~\ref{cc:sec:exact} require a seed whose growth controls every
transverse amplitude at the origin and whose first coupling stays positive.
The first wave supplies the shear
for the second at time zero. Each later wave grows in the shear of its
predecessor. Its gain allows us to switch it on while its amplitude is
still small compared with its maximum.

\subsection{The initial wave over a stationary shear}\label{cc:uf:sec:first-wave}
For $0<\rho\le1$ let
\begin{equation}\label{cc:uf:base:eq:stationary}
 U_0(x)=\sin(2\pi x_2)(e_1+\rho e_3),\qquad B_\rho=2\pi(e_1+\rho e_3)\otimes e_2
\end{equation}
be the stationary shear of Section~\ref{cc:sec:problem} and its gradient
at the origin. The field $U_0$ is smooth, odd, divergence free and
supported on the Fourier modes $\pm e_2$. As it depends only on $x_2$
and has no second component, $U_0\cdot\nabla U_0=0$. Its pressure and
residual therefore vanish. Each particle keeps its second coordinate,
so the flow and its deformation are
\begin{equation}\label{cc:uf:base:eq:deformation}
 X_0(t,a)=a+t\sin(2\pi a_2)(e_1+\rho e_3),\qquad
 (D_aX_0)^{\pm1}(t,a)=\Id\pm t\cos(2\pi a_2) B_\rho.
\end{equation}
The inverse formula uses $B_\rho^2=0$. The origin is a fixed particle and
$|B_\rho|\le4\pi$.

Take the first wave of frequency $\lambda$ over $U_0$ on $[0,T]$ with
\begin{equation}\label{cc:uf:base:eq:first-wave}
 k_*=\lambda e_1,\qquad \widetilde\chi\equiv1,
\end{equation}
so it is fully present at time zero. At the origin solve the seed
equation \eqref{cc:uf:step:eq:seed} with initial direction $e_2$ and
scale the solution so that $\max_{[0,T]}|\mathbf a(t,0)|=\lambda^{-3/4}$.
At time zero $\kappa(0,a)=\lambda e_1$ at every label, so the projection
in \eqref{cc:uf:tower:eq:start} leaves this initial amplitude unchanged.
The initial frame at the origin is $(e_1,e_2,e_3)$, in which
\begin{equation}\label{cc:uf:base:eq:shear-couplings}
 e_1\cdot B_\rho e_2=2\pi,\qquad e_3\cdot B_\rho e_2=2\pi\rho .
\end{equation}
The second wave starts at time zero, so its amplification parameter is
the ratio $\eta(0)/\nu(0)=\rho$ of these couplings in the first wave's
frame. The shear that amplifies it is the leading
gradient $\mathbf a(0,0)\otimes\kappa(0,0)$ contributed by the first
wave in \eqref{cc:uf:tower:eq:central-return}.

\begin{lemma}\label{cc:uf:base:lem:first}
If $T|B_\rho|\le1/8$, the first wave satisfies the two conditions
\eqref{cc:uf:tower:eq:wave-origin} for every $\mathsf D\ge2$, and its
strength at the origin at time zero lies between $\lambda^{1/4}/2$ and
$\lambda^{1/4}$.
\end{lemma}
\begin{proof}
At the origin the seed equation reads
$\mathbf b_t=-(\Id-2P_\kappa)B_\rho\mathbf b$. Since $\Id-2P_\kappa$
is a reflection, every nonzero solution satisfies
$\bigl|\partial_t\log|\mathbf b|\bigr|\le|B_\rho|$. Applying this estimate
to $\mathbf b$ and to the seed $\mathbf a$, whose length is proportional
to $\gamma$, gives
\[
 \frac{|\mathbf b(t)|}{\gamma(t)}
 \le e^{2T|B_\rho|}\frac{|\mathbf b(t')|}{\gamma(t')}
 \le 2\frac{|\mathbf b(t')|}{\gamma(t')}\qquad(0\le t'\le t\le T).
\]
The same estimate on the seed gives
$e^{-T|B_\rho|}\lambda^{-3/4}\le|\mathbf a(0,0)|\le\lambda^{-3/4}$.
Multiplying by $|\kappa(0,0)|=\lambda$ gives the initial strength.

For the sign of the coupling, the frame equations \eqref{cc:eq:frame}
give $|\partial_tE_1|,|\partial_tE_2|\le|B_\rho|$. Thus the first coupling
can change by at most $2t|B_\rho|^2$ from its value in
\eqref{cc:uf:base:eq:shear-couplings}, and
\[
 \nu(t,0)=E_1(t)\cdot B_\rho E_2(t)\ge 2\pi-2t|B_\rho|^2\ge\pi .
\]
\end{proof}

\subsection{Choosing the next wave}\label{cc:uf:sec:later-wave}
Let $U$ and its flow $X$ be as in Subsection~\ref{cc:uf:sec:parent-step}.
The new wave starts at $t_0$ with frequency $\lambda$. The origin is a
fixed particle. Once the preceding wave's cutoff
equals one, \eqref{cc:uf:tower:eq:central-return} and the truncation
estimate \eqref{cc:uf:step:eq:truncated} give its contribution to the
gradient there. Suppose, accordingly, that on $[t_0,T]$ the gradient
$M=\nabla U(t,0)$ has the form
\begin{equation}\label{cc:uf:tower:eq:split}
 M=B+\mathbf a_p\otimes\kappa_p+E_c ,
\end{equation}
where $B$ is the smooth trace-free gradient of the older velocity, $(\mathbf a_p,\kappa_p)$ is
the preceding wave over $B$, and $E_c$ contains the errors from its
corrections and truncation. The preceding pair solves
\eqref{cc:eq:wave-ode} with $D=B$. In its moving frame $(E_i)$ write
$B_{ij}=E_i\cdot BE_j$. At $t_0$ freeze its strength
$\sigma=|\mathbf a_p||\kappa_p|$ and its couplings $\nu=B_{12}>0$ and
$\eta=B_{32}>0$. As in Subsection~\ref{cc:sec:comparison}, put
$r=\sqrt{2\sigma\nu}$, $S=\sigma/r$, $\varepsilon=\eta/\nu$ and
$\tau=r(t-t_0)$.

Choose the wave vector at activation perpendicular to the preceding
wave vector and amplitude. Since $\kappa=(D_aX)^{-T}k_*$, the
time-zero vector that gives $\kappa(t_0,0)=\lambda E_3(t_0)$ is
\[
 k_*=\lambda (D_aX(t_0,0))^TE_3(t_0).
\]
For the amplitude take the growing direction
$-E_1(t_0)/S+E_2(t_0)$. In the coordinates of \eqref{cc:eq:scaled-pair}
this is $(\alpha_1,\alpha_2)=(-1,1)$ at $\tau=0$. Scale the seed by a positive
constant so that $\max_{[t_0,T]}|\mathbf a(t,0)|=\lambda^{-3/4}$, and
write $\gamma(t)=\lambda^{3/4}|\mathbf a(t,0)|$ as in
Subsection~\ref{cc:uf:sec:parent-step}.

To apply the comparison of Subsection~\ref{cc:sec:amplification}, we need
to bound the variation of the two couplings of $B$ in the moving frame.
When $\varepsilon\le1$, the comparison error
\eqref{cc:uf:exact:eq:additive-error} satisfies
\begin{equation}\label{cc:uf:tower:eq:active-error}
 \delta\le C\left[
 \frac{(T-t_0)\bigl(\norm B_\infty^2+\norm{B_t}_\infty\bigr)}{\eta}
 +\frac{\norm B_\infty}{r}+\frac{\norm{E_c}_\infty}{\nu}\right].
\end{equation}
Here the suprema are over $[t_0,T]$. Indeed, the frame equations
\eqref{cc:eq:frame} give $|\partial_t E_i|\le2\norm B_\infty$, so
$|\partial_t B_{ij}|\le4\norm B_\infty^2+\norm{B_t}_\infty$.
Integration from $t_0$ bounds the changes in $B_{12}$ and $B_{32}$.
Dividing by their frozen values $\nu$ and $\eta$, with $\eta\le\nu$,
gives the first term of \eqref{cc:uf:tower:eq:active-error}. The other
two terms occur directly in \eqref{cc:uf:exact:eq:additive-error}.

Suppose the normalized interval has length $r(T-t_0)\ge2$ and that
\eqref{cc:eq:finite-smallness} holds with this length in place of $T$.
Lemmas~\ref{cc:lem:finite-comparison} and~\ref{cc:lem:new-wave-couplings}
then apply to the chosen direction. The positive rescaling of the seed
preserves its frame and all ratios of its lengths. Write
$\nu_+,\eta_+$ for the couplings of $M$ in the frame of the new wave.
For $t_0\le t'\le t\le T$ and every transverse solution $\mathbf b$
of the seed equation at the origin, the two lemmas give
\begin{equation}\label{cc:uf:tower:eq:later-origin}
 \frac{|\mathbf b(t)|}{\gamma(t)}
 \le CS\,\frac{|\mathbf b(t')|}{\gamma(t')},\qquad
 \nu_+(t)\ge\tfrac12\nu\varepsilon.
\end{equation}
These are the two estimates \eqref{cc:uf:tower:eq:wave-origin} when
$\mathsf D$ bounds $CS$ and $2/(\nu\varepsilon)$.

We switch on the wave during the first unit of normalized time. With
the function $\chi$ of Subsection~\ref{cc:uf:sec:parent-step}, take
\[
 \widetilde\chi(t)=1-\chi\bigl(1+r(t-t_0)\bigr)\quad(t\ge t_0),\qquad
 \widetilde\chi(t)=0\quad(t<t_0).
\]
This cutoff is smooth, equals one for $\tau\ge1$ and satisfies
$|\widetilde\chi'|\le Cr$. Prescribe the coupling ratio
$0<\varepsilon_+\le\varepsilon/32$ for the following wave, and suppose
$r(T-t_0)\ge2(\varepsilon\varepsilon_+)^{-1/2}$. By \eqref{cc:eq:gain},
\begin{equation}\label{cc:uf:tower:eq:switch-small}
 \gamma\le\exp\bigl(3-\varepsilon^{-1/2}\bigr)
 \qquad\hbox{on }[t_0,t_0+1/r].
\end{equation}
Thus the factor $\|\widetilde\chi'\gamma\|_\infty$ in the residual
estimate \eqref{cc:uf:step:eq:physical} is at most
$Cr\exp(3-\varepsilon^{-1/2})$.

Proposition~\ref{cc:prop:finite-return} also gives a time in
$\tau_m\le\tau\le2(\varepsilon\varepsilon_+)^{-1/2}$ at which
$\eta_+/\nu_+=\varepsilon_+$. We activate the following wave at the
first such crossing in this interval. Since the reference amplitude
peaks at $\tau_m\ge2/\varepsilon>1$, the cutoff already equals one
there, and \eqref{cc:eq:general-return} gives the strength of the
shear available to amplify that wave.

\section{A common sequence of approximate solutions}\label{cc:uf:sec:tower}

We construct a single sequence of approximate solutions for all
derivative orders. The parent for wave $1$ is the stationary shear
\eqref{cc:uf:base:eq:stationary} with $\rho=\varepsilon_1$. At step $j$
we choose the wave over $U_{j-1}$ as in Section~\ref{cc:sec:iteration}
and apply Proposition~\ref{cc:uf:step:thm:step} to obtain $U_j$.
The frequencies and orders are fixed by \eqref{cc:uf:tower:eq:schedule}.
The residual of $U_j$ must already be small on the scale of the
frequency $\lambda_{j+1}$ of the next wave.

Write $X_j$ for the flow of $U_j$ from time zero, and $\pi_j$ and $R_j$
for its pressure and its residual, defined as $p_U$ and $R_U$ are for
the parent in Subsection~\ref{cc:uf:sec:parent-step}.
At the center, let $\mathbf a_j,\kappa_j,\sigma_j$ be the seed,
wave vector and strength of wave $j$, and let
\begin{equation}\label{cc:uf:tower:eq:couplings}
 \nu_j=\frac{\kappa_j}{|\kappa_j|}\cdot\nabla U_{j-1}(t,0)
        \frac{\mathbf a_j}{|\mathbf a_j|},\qquad
 \eta_j=\frac{\kappa_j\times\mathbf a_j}{|\kappa_j||\mathbf a_j|}
        \cdot\nabla U_{j-1}(t,0)\frac{\mathbf a_j}{|\mathbf a_j|}
\end{equation}
be the couplings of $\nabla U_{j-1}(t,0)$ in the frame of this pair.
For $j\ge2$ these are $\nu_+,\eta_+$ of Subsection~\ref{cc:uf:sec:later-wave}.
Let $\widetilde\chi_j$ be the time cutoff of wave $j$ and
$O_j=W_{\rm osc}^\natural$ its oscillatory velocity increment before
truncation, as in Lemma~\ref{cc:uf:step:lem:approximate}.

The first two waves start at time zero. Given the activation time $t_j$,
we fix the final time $S_j$ of the endpoint equations by
\[
 t_1=t_2=0,\qquad
 \tau_j=\lambda_{j-1}^{-1/8}L_{j+1}^3,\qquad
 S_0=1,\qquad S_j=t_j+2\tau_j\quad(j\ge1).
\]
For $j\ge2$ we normalize the seed on $[t_j,S_j]$. The central wave
equations over $U_{j-1}$ then determine $t_{j+1}$ on that interval.
We take the first crossing
of $\eta_j/\nu_j=\varepsilon_j$ in the interval of
Proposition~\ref{cc:prop:finite-return}. It will occur by $t_j+\tau_j$,
leaving room before $S_j$ for the subsequent activation times.
To make the remainder of the expansion small at frequency $\lambda_{j+1}$,
we take the depth
\begin{equation}\label{cc:uf:tower:eq:step-orders}
 P=6K_{j+1}.
\end{equation}

\Needspace{8\baselineskip}
\begin{theorem}\label{cc:uf:tower:thm:finite}
For every sufficiently large $N$ there are fields $U_j$, $j\ge1$, and
times $0=t_1=t_2<t_3<\cdots$ with the following properties, in which
$U_0$ is the shear \eqref{cc:uf:base:eq:stationary} at $\rho=\varepsilon_1$:
\textup{(T1)}--\textup{(T3)} for $j\ge0$ and \textup{(T4)} for $j\ge1$.
\begin{enumerate}[label=(T\arabic*),leftmargin=*]
\item $U_j$ is a smooth odd divergence-free field on
$[0,S_j]\times\mathbb T^3$, supported on frequencies at most $2\Lambda_j$.
\item $\norm{U_j}_{C_tL^\infty}\le\lambda_j$,
$\sup|(D_aX_j)^{\pm1}|\le\lambda_j$ and
$\norm{R_j}_{C_tL^\infty}\le\lambda_{j+1}^{-3}$.
\item $D^2\pi_j\le\bigl(\sum_{i=1}^j\lambda_{i-1}^{-1}\bigr)\Id$ on
$[0,S_j]\times\mathbb T^3$.
\item At the center on $[t_j,S_j]$,
$\nabla U_j=\nabla U_{j-1}
+\widetilde\chi_j\,\mathbf a_j\otimes\kappa_j+E_{c,j}$ with
$|E_{c,j}|\le\lambda_j^{-1/4}$. Moreover
$0\le t_{j+1}-t_j\le\tau_j$ and $\widetilde\chi_j=1$ on
$[t_{j+1},S_j]$, so that the gradient of $U_j$ has the splitting
\eqref{cc:uf:tower:eq:split} there. At $t_{j+1}$, with an absolute
constant $C\ge2$ for which \eqref{cc:eq:general-return} holds,
\begin{equation}\label{cc:uf:tower:eq:late-values}
 \frac{\eta_j}{\nu_j}=\varepsilon_j,\qquad
 \pi\le\nu_j\le4\pi,\qquad
 C^{-1}\lambda_j^{1/4}\le\sigma_j\le C\lambda_j^{1/4}L_{j+1}^2 .
\end{equation}
\end{enumerate}
The activation times converge to a limit
$T_\infty\in(0,\infty)$, and $T_\infty<S_j$ for every $j$. For $j\ge1$
the oscillation $O_j$ is supported, at each time of $[0,S_j]$, in the
ball of radius $\lambda_{j-1}^{-1}$ about the origin, and
\begin{equation}\label{cc:uf:tower:eq:increment}
 \|O_j\|_{C([0,S_j];L^\infty)}\le\lambda_j^{-1},\qquad
 \|U_j-U_{j-1}-O_j\|_{C([0,S_j];\Holder{m})}\le\lambda_j^{-1}
 \quad(0\le m\le q_j).
\end{equation}
For $j\ge2$ the oscillation $O_j$ vanishes on $[0,t_j]$.
\end{theorem}

\begin{figure}[!htbp]
\centering
\includegraphics[width=\textwidth]{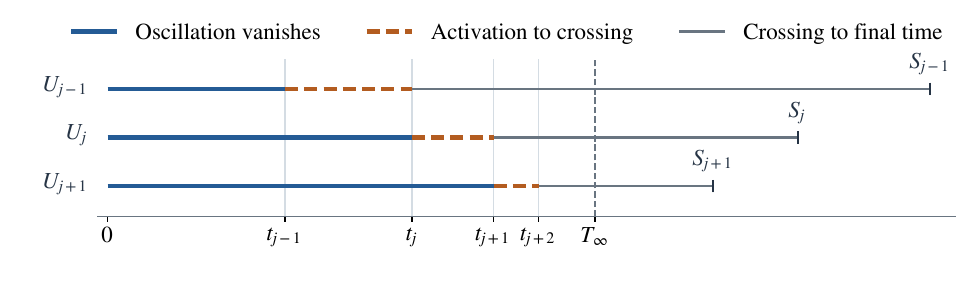}
\caption{Three successive approximate solutions for $j\ge4$. In each row
the oscillation vanishes on the thick solid segment $[0,t_i]$. The
dashed segment ends at the crossing $t_{i+1}$, when the next wave is
activated. Each solution continues to its final time $S_i>T_\infty$.}
\label{cc:fig:construction-timeline}
\end{figure}

For each fixed $m$, the orders $q_j$ eventually exceed $m$. By
\eqref{cc:uf:tower:eq:increment}, the increments $U_j-U_{j-1}-O_j$
are then summable in $\Holder{m}$ on the common interval $[0,T_\infty]$.
On $[0,t_j]$ they give the entire velocity increment for $j\ge2$.
This is the estimate used for convergence in
Section~\ref{cc:sec:limit}. Figure~\ref{cc:fig:construction-timeline}
shows the intervals on which the oscillations are present.

\subsection{The data of a step}

At any index $j$ for which \textup{(T1)} and \textup{(T2)} hold, the
frequency bound gives estimates for the gradient and its time derivative.
As in Subsection~\ref{cc:uf:sec:parent-step}, the residual is supported
on frequencies at most $4\Lambda_j$ and
\begin{equation}\label{cc:uf:tower:eq:leray}
 \partial_tU_j=-\PP\diver(U_j\otimes U_j)+R_j .
\end{equation}
By \eqref{cc:uf:trunc:eq:lp},
$\|\nabla U_j\|_\infty\le C\Lambda_j\lambda_j$ and
$\|\nabla R_j\|_\infty\le C\Lambda_j\lambda_{j+1}^{-3}$.
For the time derivative of the gradient, apply \eqref{cc:uf:trunc:eq:lp}
to $U_j\otimes U_j$ and use the boundedness of $\PP$ on $\Holder{0}$
in \eqref{cc:uf:tower:eq:leray}. Two spatial derivatives of the product
in $\Holder{0}$ cost at most $C\Lambda_j^3\lambda_j^2=C\lambda_j^{17}$.
Thus for large $N$,
\begin{equation}\label{cc:uf:tower:eq:derived-low}
 \|\nabla U_j\|_\infty\le\lambda_j^9,\qquad
 \|\partial_t\nabla U_j\|_\infty\le\lambda_j^{18}
 \quad\hbox{on }[0,S_j].
\end{equation}
At step $j$ we apply Section~\ref{cc:sec:exact} with $U=U_{j-1}$,
$\lambda=\lambda_j$, activation time $t_j$ and final time $S_j$. Write
$\mu=\lambda_{j-1}$ and $\widetilde\chi=\widetilde\chi_j$. By
\eqref{cc:uf:tower:eq:derived-low} at $j-1$ we can take
$\mathsf D=\mu^9$ for the parent gradient. When $j\ge2$ we also use
\eqref{cc:uf:tower:eq:derived-low} at $j-2$ for the older gradient in
the wave comparison. The final Fourier cutoff is $\Lambda_j=\lambda^5$.

To compare successive final times, observe that by \eqref{cc:uf:tower:eq:uniform},
$\tau_j=\mu^{-1/8}L_{j+1}^3\le\mu^{-1/16}$ for large $N$.
Applying this bound at $j+1$ and using $\lambda_j\ge\mu^4$ gives
$4\tau_{j+1}\le4\mu^{-1/4}\le\mu^{-1/8}\le\tau_j$. Thus
\begin{equation}\label{cc:uf:tower:eq:increment-tail}
 \tau_j\le\lambda_{j-1}^{-1/16},\qquad 4\tau_{j+1}\le\tau_j ,
\end{equation}
uniformly in $j$ for large $N$.

Finally the depth $P=6K_{j+1}$ in \eqref{cc:uf:tower:eq:step-orders}
requires estimates through order $P+1$. Since $P+1\le7K(q_{j+1})$,
$\mathcal B(q_{j+1})$ bounds all the fixed-order constants used in
Section~\ref{cc:sec:exact} and Lemma~\ref{cc:uf:trunc:lem:cutoff}
at these orders and this depth.
The choice of $q_{j+1}$ in \eqref{cc:uf:tower:eq:schedule} gives
$\mathcal B(q_{j+1})\le\mu$. With $\mathsf D=\mu^9$ and
$\lambda=\mu^{K_j}$, $K_j=320b(q_j+1)$, we obtain
\begin{equation}\label{cc:par:step}
 \mathcal B(q_{j+1})\mathsf D^b\le\mu^{10b}\le\lambda^{1/32}.
\end{equation}

\subsection{Proof of Theorem~\ref{cc:uf:tower:thm:finite}}

For $j=0$, the shear $U_0$ has only the Fourier modes $\pm e_2$,
$\|U_0\|_\infty\le2$, and its pressure and residual vanish. Its
deformation and inverse have norm at most $1+4\pi$ on $[0,S_0]$
by \eqref{cc:uf:base:eq:deformation}. These give \textup{(T1)}--\textup{(T3)}
for large $N$.

Suppose the construction has been carried out through step $j-1$.
For $j=1$, \eqref{cc:uf:tower:eq:increment-tail} gives
$S_1=2\tau_1\le S_0=1$. For $j\ge2$, the previous crossing occurred
by $t_{j-1}+\tau_{j-1}$. Since $2\tau_j\le\tau_{j-1}/2$,
\[
 S_j=t_j+2\tau_j\le t_{j-1}+\tfrac32\tau_{j-1}<S_{j-1}.
\]
Thus $U_{j-1}$ is available on the entire interval $[0,S_j]$ of
the endpoint problem.

\begin{enumerate}[leftmargin=*]
\item \emph{Estimates for the parent.}
On $[0,S_j]\subset[0,S_{j-1}]$, (T1) bounds the frequencies of $U=U_{j-1}$
by $2\Lambda_{j-1}=2\mu^5\le2\mathsf D$. By (T2) and
\eqref{cc:uf:tower:eq:derived-low}, $\sup|(D_aX)^{\pm1}|\le\mu$ and
$\|\nabla U\|_\infty\le\mu^9=\mathsf D$, as required in
\eqref{cc:uf:tower:eq:parent-bound}. The residual bound
$\|R_{j-1}\|_{C_tL^\infty}\le\lambda^{-3}$ in (T2) gives
\eqref{cc:uf:tower:eq:parent-small}. For the endpoint problem we also need
\eqref{cc:uf:tower:eq:parent-coercive}. The derivative bound preceding
\eqref{cc:uf:tower:eq:derived-low} gives $\|\nabla R_{j-1}\|_\infty
\le C\mu^5\lambda^{-3}$. In (T3) the sum is at most $2/N$, since $\lambda_0=N$ and
$\lambda_i\ge2\lambda_{i-1}$. Thus
$D^2\pi_{j-1}-\nabla R_{j-1}\le(2/N+C\mu^5\lambda^{-3})\Id$
as quadratic forms, which gives \eqref{cc:uf:tower:eq:parent-coercive} for large $N$.

\item \emph{Choosing wave $j$.}
For $j=1$ take the wave \eqref{cc:uf:base:eq:first-wave} of frequency
$\lambda_1$ on $[0,S_1]$. Since $4\pi S_1\le1/8$ for large $N$,
Lemma~\ref{cc:uf:base:lem:first} gives \eqref{cc:uf:tower:eq:wave-origin}
and the strength bounds in \eqref{cc:uf:tower:eq:late-values} at
$t_2=t_1=0$. The couplings there are $\nu_1=2\pi$ and
$\eta_1=2\pi\varepsilon_1$ by \eqref{cc:uf:base:eq:shear-couplings}.
The first wave has $\widetilde\chi_1\equiv1$.

For $j\ge2$, (T4) at $j-1$ gives the splitting
\eqref{cc:uf:tower:eq:split} on $[t_j,S_j]$ with $B=\nabla U_{j-2}(t,0)$,
the preceding pair $(\mathbf a_{j-1},\kappa_{j-1})$ and
$\|E_c\|_\infty\le\mu^{-1/4}$. The pair solves \eqref{cc:eq:wave-ode}
over $B$ because the flow of $U_{j-2}$ fixes the origin. Freeze its
couplings and strength at $t_j$. For $j=2$ these are the values of
wave~$1$ at time zero. For $j\ge3$ they are the values at the preceding crossing.
In both cases \eqref{cc:uf:tower:eq:late-values} gives
\[
 \varepsilon=\varepsilon_{j-1}=L_j^{-2},\qquad
 \pi\le\nu\le4\pi,\qquad
 c\mu^{1/4}\le\sigma\le C\mu^{1/4}L_j^2 .
\]
Choose the wave of Subsection~\ref{cc:uf:sec:later-wave}, whose wave
vector has length $\lambda$ at $t_j$. Since $r=\sqrt{2\sigma\nu}$,
$c\mu^{1/8}\le r\le C\mu^{1/8}L_j$. The normalized length
$r(S_j-t_j)=2r\mu^{-1/8}L_{j+1}^3$ is therefore between $2$ and
$CL_{j+1}^4$ for large $N$.

To check the comparison on this interval, \eqref{cc:uf:tower:eq:derived-low}
at $j-2$ gives $\|B\|_\infty\le\lambda_{j-2}^9\le\mu^{1/32}$ and
$\|B_t\|_\infty\le\mu^{1/16}$, since $\lambda_{j-2}=\mu^{1/K_{j-1}}$
and $K_{j-1}>288$. With $\eta=\nu\varepsilon$, the three terms in
\eqref{cc:uf:tower:eq:active-error} give
\[
 \delta\le C\bigl(\varepsilon^{-1}L_{j+1}^3\mu^{-1/16}
                   +\mu^{-3/32}+\mu^{-1/4}\bigr).
\]
By \eqref{cc:uf:tower:eq:uniform}, $L_j\le L_{j+1}\le(\log\mu)^2$.
Thus $C_0(2r\tau_j)^{C_0}\delta/\varepsilon$ is bounded by
$\mu^{-1/16}$ times a fixed power of $\log\mu$, and
$C_0\varepsilon\le C_0(\log N)^{-2}$. For large $N$ both conditions
of \eqref{cc:eq:finite-smallness} hold. In the resulting central estimates
\eqref{cc:uf:tower:eq:later-origin}, $CS=Cr/(2\nu)\le C\mu^{1/8}L_j$
and $2/(\nu\varepsilon)\le 2L_j^2/\pi$. Both are at most $\mu\le\mathsf D$
for large $N$. This proves \eqref{cc:uf:tower:eq:wave-origin} for wave $j$.

\item \emph{The next activation time.}
For $j\ge2$ we seek a crossing with $\eta_j/\nu_j=\varepsilon_j$ in
the first half of $[t_j,S_j]$. Then $t_{j+1}\le t_j+\tau_j$ leaves
room for the next interval, since $2\tau_{j+1}\le\tau_j/2$ by
\eqref{cc:uf:tower:eq:increment-tail}. The parameter choice
\eqref{cc:uf:tower:eq:schedule} gives $L_{j+1}/L_j=K_{j+2}\ge640$, so
\begin{equation}\label{cc:uf:tower:eq:ratio}
 \varepsilon_j=(L_j/L_{j+1})^2\varepsilon_{j-1}\le\varepsilon_{j-1}/32,
\end{equation}
and $(\varepsilon\varepsilon_j)^{-1/2}=L_jL_{j+1}\le L_{j+1}^2$.
In physical time the upper endpoint of the bracket in
Proposition~\ref{cc:prop:finite-return} lies at most $\tau_j$ after $t_j$. Indeed,
\[
 \frac{2(\varepsilon\varepsilon_j)^{-1/2}}{r}
 \le C\mu^{-1/8}L_{j+1}^2\le\tau_j ,
\]
using $r\ge c\mu^{1/8}$ and $L_{j+1}\ge\log N$. Take $t_{j+1}$ to be
the first crossing in that interval. Its normalized time is at least
$\tau_m\ge2/\varepsilon>1$, so $0<t_{j+1}-t_j\le\tau_j$ and
$\widetilde\chi=1$ on $[t_{j+1},S_j]$.

At this crossing \eqref{cc:eq:general-return} bounds $\sigma_j(t_{j+1})$
above and below by absolute multiples of $\lambda^{1/4}L_jL_{j+1}$,
since $\max_{[t_j,S_j]}|\mathbf a_j(t,0)|=\lambda^{-3/4}$.
As $1\le L_jL_{j+1}\le L_{j+1}^2$, this gives the strength bounds in
\eqref{cc:uf:tower:eq:late-values}. Equation~\eqref{cc:eq:actual-ratio}
applies at the same time and gives
\[
 \bigl|\log\bigl(\nu_j(t_{j+1})/\nu_{j-1}(t_j)\bigr)\bigr|
 \le2\varepsilon_{j-1}.
\]
Summing this inequality over the successive crossings from
$\nu_1(t_2)=2\pi$ and using
$\varepsilon_i\le32^{1-i}\varepsilon_1$ from \eqref{cc:uf:tower:eq:ratio}
gives $|\log(\nu_j(t_{j+1})/(2\pi))|\le4\varepsilon_1\le4(\log N)^{-2}$.
Thus $\pi\le\nu_j(t_{j+1})\le4\pi$ for large $N$, completing
\eqref{cc:uf:tower:eq:late-values}.

On the interval where the cutoff changes, \eqref{cc:uf:tower:eq:switch-small}
gives $\gamma\le e^{3-L_j}=e^3\lambda_{j+1}^{-4}$. Since
$|\widetilde\chi'|\le Cr$ and $Ce^3r\le\mu$ for large $N$ by item~2,
\begin{equation}\label{cc:uf:tower:eq:switch-induction}
 \|\widetilde\chi'\gamma\|_\infty\le\mu\lambda_{j+1}^{-4}.
\end{equation}
The same estimate holds for the first wave since $\widetilde\chi_1'\equiv0$.

\item \emph{The expansion and its truncation.}
Items 1 and 2 give the parent and wave assumptions of
Subsection~\ref{cc:uf:sec:parent-step}. By \eqref{cc:par:step},
$C_{\rm step}=\Gamma_P\mathsf D^b\le\mu^{10b}\le\lambda^{1/32}$,
so \eqref{cc:uf:step:eq:small} holds. We can therefore apply the
expansion and truncation estimates of Section~\ref{cc:sec:exact}.

The new residual must be small at the frequency $\lambda_{j+1}$ of the
next wave. For the expansion remainder, the depth
\eqref{cc:uf:tower:eq:step-orders} gives
$\lambda^{-2P/3}=\lambda^{-4K_{j+1}}=\lambda_{j+1}^{-4}$.
For the time cutoff, item 3 gives the bound
\eqref{cc:uf:tower:eq:switch-induction}, with the same decay up to a factor $\mu$. Thus
\eqref{cc:uf:step:eq:physical} and $C_{\rm step}\mu\le\lambda^{1/16}$ give
\begin{equation}\label{cc:uf:tower:eq:frequency}
 \norm{R_V}_{C_t\Holder{1}}\le2\lambda^{2+1/16}\lambda_{j+1}^{-4}
 \le\lambda_{j+1}^{-7/2}.
\end{equation}

To estimate the Fourier tail in \eqref{cc:uf:step:eq:truncation}, we
need $\norm V_{C_t\Holder{P}}$. For the parent $U=U_{j-1}$, (T1), (T2)
and the Bernstein bound \eqref{cc:uf:trunc:eq:lp} give
$\norm U_{C_t\Holder{P}}\le C_P\mu^{5(P+1)+1}\le\lambda^{P+2}$,
since $C_P\le\mu$ and $\lambda=\mu^{K_j}$ with $K_j\ge6$.
For the wave, \eqref{cc:par:step} bounds the factor
$C_{P,P}\mathsf D^{b(P+2)}$ in \eqref{cc:uf:step:eq:physical}
by $\lambda^{(P+2)/32}$. The larger of the two resulting exponents
for the oscillation and the mean is $3P/2+1/4+(P+2)/32$.
This gives $\norm{W^\natural}_{C_t\Holder{P}}\le\lambda^{2(P+1)}$
for large $N$. Hence $1+\norm V_{C_t\Holder{P}}\le3\lambda^{2(P+1)}$.
Inserting this with $\Lambda=\lambda^5$ and $9\mu\le\lambda$ gives
\begin{equation}\label{cc:uf:tower:eq:tail-small}
 \mu\varepsilon_\Lambda\le9\mu\lambda^{5(2-P)+4(P+1)}
 \le\lambda^{15-P}
 =\lambda^{15}\lambda_{j+1}^{-6}\le\lambda_{j+1}^{-5}.
\end{equation}
In particular $C_{2,P}\varepsilon_\Lambda\le\lambda_{j+1}^{-5}$.

\item \emph{The new velocity.}
Take $U_j=S_{\Lambda_j}V$ in Proposition~\ref{cc:uf:step:thm:step}.
The proposition gives (T1), and \eqref{cc:uf:step:eq:truncated} with
\eqref{cc:uf:tower:eq:frequency} and \eqref{cc:uf:tower:eq:tail-small} gives
\[
 \norm{R_j}_{C_tL^\infty}\le C\lambda_{j+1}^{-7/2}+\lambda_{j+1}^{-5}
 \le\lambda_{j+1}^{-3},\qquad
 \norm{v_\Lambda}_{C_tC^2}\le\lambda_{j+1}^{-5}.
\]
Since $U_j=V+v_\Lambda$ and $\lambda_{j+1}^{-5}\le\lambda^{-5/2}$,
the divergence-free field $v_\Lambda$ satisfies the hypothesis of
Lemma~\ref{cc:uf:step:lem:flow-bounds}. It follows that
$\sup|(D_aX_j)^{\pm1}|\le C_{\rm step}\lambda^{1/4}\le\lambda$.

By Proposition~\ref{cc:uf:hier:prop:identity}, $W_{\rm osc}$ is supported
in $|a|\le2\ell_*$ in label coordinates. The support of
$O_j=W_{\rm osc}^\natural$ lies in the image of this ball under the
parent flow. The flow fixes the origin, so
\eqref{cc:uf:tower:eq:parent-bound} and $\ell_*\le\mathsf D^{-3}$
place its Eulerian support in the ball of radius
$2\mathsf D\ell_*\le2\mathsf D^{-2}\le\mu^{-1}$.
For $j\ge2$, $O_j$ vanishes before $t_j$ by the same proposition.
It vanishes also at $t_j$, since $\widetilde\chi(t_j)=0$ and the higher
oscillatory profiles have zero data there by Subsection~\ref{cc:uf:sec:transverse}.
By \eqref{cc:uf:step:eq:physical-sup} and \eqref{cc:par:step},
$\norm{O_j}_{C_tL^\infty}\le\lambda^{1/32}\lambda^{-5/4}\le\lambda^{-1}$.

The rest of the velocity increment can extend across the torus. Since
$V=U_{j-1}+W^\natural$,
\[
 U_j-U_{j-1}-O_j=\bar W^\natural+v_\Lambda .
\]
For $0\le m\le q_j$, the coefficient $C_{m,P}\mathsf D^{b(m+2)}$
in \eqref{cc:uf:step:eq:physical} is at most $\mu^{10b(m+2)}\le
\lambda^{1/16}$, since $K_j=320b(q_j+1)$. Thus
$\|\bar W^\natural\|_{C_t\Holder{m}}\le\lambda^{1/16-5/2}
\le\lambda^{-3/2}$. For the Fourier tail, \eqref{cc:uf:trunc:eq:tail}
at order $P$ and \eqref{cc:uf:tower:eq:tail-small} give
\[
 \|v_\Lambda\|_{C_t\Holder{m}}
 \le\mu\Lambda_j^{m-P}\norm V_{C_t\Holder{P}}
 \le\Lambda_j^{m-2}\mu\varepsilon_\Lambda
 \le\lambda^{5q_j}\lambda_{j+1}^{-5}\le\lambda_{j+1}^{-4},
\]
where $K_{j+1}\ge K_j\ge5q_j$. The sum of the mean and tail bounds is
at most $\lambda^{-1}$, proving \eqref{cc:uf:tower:eq:increment}.
In particular $\norm{U_j}_{C_tL^\infty}\le\mu+2\lambda^{-1}\le\lambda$
by (T2) for the parent. This completes (T2) for $U_j$.

\item \emph{Pressure and the central velocity gradient.}
For the pressure increment, the leading term in \eqref{cc:uf:step:eq:pressure-error}
satisfies $C_{\rm step}\lambda^{-1/4}\le\lambda^{-7/32}\le\mu^{-2}$.
By \eqref{cc:uf:tower:eq:frequency} and \eqref{cc:uf:tower:eq:tail-small},
the other two terms are at most $C\lambda_{j+1}^{-7/2}+\lambda_{j+1}^{-5}$.
Hence $D^2\pi_j\le D^2\pi_{j-1}+\mu^{-1}\Id$ on $[0,S_j]$ for large $N$.
Since $\mu^{-1}=\lambda_{j-1}^{-1}$ and $S_j\le S_{j-1}$, adding this
increment to (T3) for the parent proves (T3) at step $j$.

At the origin truncation adds $\nabla v_\Lambda(t,0)$ to the gradient
in \eqref{cc:uf:tower:eq:central-return}. Thus the new error is
$E_{c,j}=E_c+\nabla v_\Lambda(t,0)$ with
$|E_{c,j}|\le\lambda^{-1/3}+\lambda_{j+1}^{-5}\le\lambda^{-1/4}$.
From $t_{j+1}$ onward the cutoff equals one, so the full rank-one
term is present for the next wave. With the coupling and strength bounds
from items 2 and 3, this proves (T4) and closes the induction.
\end{enumerate}

The activation times increase strictly from $t_2=0$ onward. By (T4)
their successive differences are bounded by $\tau_j$, so
\eqref{cc:uf:tower:eq:increment-tail} gives a finite limit
$T_\infty\ge t_3>0$. For $j\ge1$,
\[
 0<T_\infty-t_j\le\sum_{p\ge j}\tau_p\le\tfrac43\tau_j<S_j-t_j .
\]
Thus $T_\infty<S_j$ for $j\ge1$, and $T_\infty<S_0$ since $S_1\le S_0$.

\section{The limit and its terminal behavior}\label{cc:sec:limit}
On the common interval $[0,T_\infty]$ of Theorem~\ref{cc:uf:tower:thm:finite},
\eqref{cc:uf:tower:eq:increment} gives $\|U_j-U_{j-1}\|_{C_tL^\infty}
\le2\lambda_j^{-1}$. The sum converges, so the velocities converge
uniformly to a continuous field $w$ on $[0,T_\infty]\times\mathbb T^3$.
For convergence of derivatives we use that the oscillations vanish
before activation and have supports shrinking to the origin.

\subsection{Convergence before the accumulation time}

\begin{proposition}\label{cc:uf:terminal:agreement}
The data $U_j(0)$ converge in every $\Holder{m}$ to an odd divergence-free field $u_0$.
In particular, $u_0\in C^\infty(\mathbb T^3;\R^3)$.
Let $u$ be the maximal smooth
solution of the Euler equations from $u_0$, with maximal time $T_*>0$, supplied by
Section~\ref{cc:sec:problem}. Then $T_*\ge T_\infty$ and, for every
$T<T_\infty$ and integer $m\ge0$,
$U_j\to u$ in $C([0,T];\Holder{m})$. The solution is smooth and odd on
$[0,T_\infty)\times\mathbb T^3$.
\end{proposition}
\begin{proof}
Fix $T<T_\infty$ and $m\ge1$. Since $t_j\uparrow T_\infty$ and
$q_j\uparrow\infty$, we can choose $j_0$ with $T\le t_{j_0}$ and
$q_{j_0}\ge m$. For $j>j_0$ the oscillation $O_j$ vanishes on $[0,T]$.
Thus \eqref{cc:uf:tower:eq:increment} gives, for $\ell>n\ge j_0$,
\[
 \|U_\ell-U_n\|_{C([0,T];\Holder{m})}
 \le\sum_{j>n}\lambda_j^{-1}\le2\lambda_{n+1}^{-1}.
\]
The limit in each of these norms is $w$. In particular $u_0=w(0)$ is
smooth, and $w$ is odd and divergence free on $[0,T]$.

Integrating \eqref{cc:uf:tower:eq:leray} gives
\[
 U_j(t)=U_j(0)-\int_0^t\PP\diver(U_j\otimes U_j)\dd t'
        +\int_0^tR_j\dd t'.
\]
By \textup{(T2)}, the last integral has supremum norm at most
$T\lambda_{j+1}^{-3}$. Convergence in $C([0,T];\Holder{1})$, the product
estimate \eqref{cc:holder:product-growth} and boundedness of $\PP$ on
$\Holder{0}$ give uniform convergence of the projected nonlinear term.
Hence $\partial_t w=-\PP\diver(w\otimes w)$. The right-hand side is
continuous in every spatial H\"older norm, so differentiating the equation
gives joint smoothness. Uniqueness identifies $w$ with $u$ wherever both
exist. Since $w$ exists for every $T<T_\infty$, maximality gives
$T_*\ge T_\infty$.
\end{proof}

\subsection{Velocity, central vorticity and the maximal time}

\begin{proposition}\label{cc:uf:terminal:velocity-curl}
The solution has bounded velocity on $[0,T_\infty)$, and
$|\curl u(t,0)|\longrightarrow\infty$ as $t\uparrow T_\infty$.
\end{proposition}
\begin{proof}
The continuous field $w$ on $[0,T_\infty]\times\mathbb T^3$ is bounded,
and Proposition~\ref{cc:uf:terminal:agreement} identifies it with $u$
before $T_\infty$.

For the vorticity, fix a large $j$ and consider the whole interval
$[t_j,t_{j+1}]$. The preceding wave is fully active there. Applying
\textup{(T4)} of Theorem~\ref{cc:uf:tower:thm:finite} at steps $j-1$
and $j$ gives, at the origin,
\[
 \nabla U_j=\nabla U_{j-2}
   +\mathbf a_{j-1}\otimes\kappa_{j-1}
   +\widetilde\chi_j\,\mathbf a_j\otimes\kappa_j+E_{c,j-1}+E_{c,j}.
\]
A gradient $\mathbf a\otimes\kappa$ contributes $\kappa\times\mathbf a$
to vorticity. Projecting $\curl U_j$ along $\kappa_j$ therefore removes
the leading curl of wave $j$. The preceding wave contributes
$\kappa_{j-1}\times\mathbf a_{j-1}=\sigma_{j-1}E_3$ in its own frame.
By \eqref{cc:eq:scaled-pair} and Lemma~\ref{cc:uf:exact:lem:additive},
the $E_3$ component of $\kappa_j/\lambda_j$ is at least $1/2$ and its
length is at most $C(1+\tau^2)$ at normalized time $\tau$.
The normalized length of $[t_j,S_j]$ is at most $CL_{j+1}^4$ by item 2
of the proof of Theorem~\ref{cc:uf:tower:thm:finite}. Thus
\begin{equation}\label{cc:uf:terminal:eq:projection}
 \frac{\kappa_j}{|\kappa_j|}\cdot(\kappa_{j-1}\times\mathbf a_{j-1})
    \ge c\sigma_{j-1}L_{j+1}^{-8},\qquad
 \frac{\kappa_j}{|\kappa_j|}\cdot(\kappa_j\times\mathbf a_j)=0.
\end{equation}

The preceding wave retains its strength throughout $[t_j,t_{j+1}]$.
Put $\mu=\lambda_{j-1}$. At $t_j$, \eqref{cc:uf:tower:eq:late-values}
gives $\sigma_{j-1}\ge c\mu^{1/4}$. By \eqref{cc:eq:frame},
$|\partial_t\log\sigma_{j-1}|\le|\nabla U_{j-2}(t,0)|
\le\lambda_{j-2}^9\le\mu^{1/32}$, using
\eqref{cc:uf:tower:eq:derived-low} and item 2 above. Since the interval
has length at most $\tau_j\le\mu^{-1/16}$ by \textup{(T4)} and
\eqref{cc:uf:tower:eq:increment-tail}, $\log\sigma_{j-1}$ changes by
at most one. Hence $\sigma_{j-1}\ge c\mu^{1/4}$ on the whole interval.
Since $L_{j+1}^8\le(\log\mu)^{16}\le\mu^{1/8}$ for large $j$ by
\eqref{cc:uf:tower:eq:uniform}, the projected curl of the preceding wave
in \eqref{cc:uf:terminal:eq:projection} is at least $c\mu^{1/8}$.
The older gradient contributes at most $C\lambda_{j-2}^9$ to the curl
and the two central errors are at most $2\mu^{-1/4}$ in norm. Consequently
\begin{equation}\label{cc:uf:terminal:eq:curl-lower}
 |\curl U_j(t,0)|
 \ge c\lambda_{j-1}^{1/8}-C\lambda_{j-2}^9-C,
\end{equation}
which tends to infinity uniformly on $[t_j,t_{j+1}]$.

For $\ell>j$ the oscillation $O_\ell$ vanishes throughout $[t_j,t_{j+1}]$.
Applying \eqref{cc:uf:tower:eq:increment} to the remaining increments
at order one and taking the limit by Proposition~\ref{cc:uf:terminal:agreement}
gives
\[
 \sup_{t_j\le t\le t_{j+1}}|\curl(u-U_j)(t,0)|
 \le C\sum_{\ell>j}\lambda_\ell^{-1}.
\]
This tail tends to zero. The successive intervals cover
$[t_j,T_\infty)$, so \eqref{cc:uf:terminal:eq:curl-lower} gives
$|\curl u(t,0)|\to\infty$ as $t\uparrow T_\infty$.
\end{proof}

\begin{corollary}\label{cc:uf:terminal:time}
The maximal existence time is $T_*=T_\infty\in(0,\infty)$, and
both integrals in \eqref{cc:uf:eq:main-breakdown} diverge.
\end{corollary}
\begin{proof}
By Proposition~\ref{cc:uf:terminal:agreement}, $T_*\ge T_\infty$. If
$T_*>T_\infty$, then smoothness on the compact interval $[0,T_\infty]$
would bound $\curl u(t,0)$, contrary to
Proposition~\ref{cc:uf:terminal:velocity-curl}. Hence $T_*=T_\infty$,
which is positive and finite by Theorem~\ref{cc:uf:tower:thm:finite}.
The initial data are smooth, so the periodic Beale--Kato--Majda criterion
\cite{BKM}, in the formulation of \cite[Theorem~1]{GGPP}, implies
$\int_0^{T_*}\|\curl u(t)\|_\infty\dd t=\infty$.
The gradient integral also diverges, since
$\|\curl u\|_\infty\le C\|\nabla u\|_\infty$.
\end{proof}

\subsection{The velocity away from the origin at the singular time}

\begin{proposition}\label{cc:uf:terminal:off-origin}
For every compact $K\subset\mathbb T^3\setminus\{0\}$ and every integer
$m\ge0$, both $u(t)$ and $\partial_tu(t)$ converge in $C^m(K)$ as
$t\uparrow T_\infty$.
\end{proposition}
\begin{proof}
Fix $K$ and $m$ as in the statement. The support of $O_j$ lies in the
ball of radius $\lambda_{j-1}^{-1}$ about the origin throughout $[0,T_\infty]$.
For large $j$ the ball is disjoint from $K$ and $q_j\ge m$, so
\eqref{cc:uf:tower:eq:increment} gives
$\|U_j-U_{j-1}\|_{C([0,T_\infty];C^m(K))}\le\lambda_j^{-1}$.
The increments are summable in this norm, hence $u(t)\to w(T_\infty)$
in $C^m(K)$.

The time derivative also contains the pressure gradient, which depends
on the velocity on the whole torus. Choose a smooth cutoff $\chi$ equal
to one near $K$ and zero near the origin. On $[0,T_\infty)$ we have
\[
 \partial_tu=-\PP\diver(\chi u\otimes u)-\PP\diver\bigl((1-\chi)u\otimes u\bigr).
\]
Applying the spatial convergence just proved on $\supp\chi$ through
order $m+2$ gives convergence of $\chi u\otimes u$ in $\Holder{m+1}$
on the torus. Since $\PP\diver$ maps $\Holder{m+1}$ boundedly to
$\Holder{m}$, the first term converges in $C^m(K)$.

For the second term, $\supp(1-\chi)$ is a positive distance from $K$.
Only the pressure part $-\nabla\Delta^{-1}\diver\diver$ acts across
this distance. Its kernel is smooth off the diagonal, so for any
bounded tensor field $A$,
\[
 \|\PP\diver((1-\chi)A)\|_{C^m(K)}\le C_{m,K,\chi}\|A\|_{L^1}.
\]
Since $u(t)\otimes u(t)\to w(T_\infty)\otimes w(T_\infty)$ uniformly
on the torus, the difference tends to zero in $L^1$. The estimate
therefore gives convergence of the second term in $C^m(K)$.
\end{proof}

Propositions~\ref{cc:uf:terminal:agreement},
\ref{cc:uf:terminal:velocity-curl} and~\ref{cc:uf:terminal:off-origin}
and Corollary~\ref{cc:uf:terminal:time} prove
Theorem~\ref{cc:target:unforced}.

\appendix
\section{The H\"older inverse}\label{cc:app:mean}
Subsection~\ref{cc:uf:sec:mean} reduces the choice of initial velocity to
the endpoint problem \eqref{cc:uf:mean:integrated}. We solve it on the
normalized interval $I=[0,1]$ and estimate the velocity in H\"older norm. Write
$\norm f_k=\sup_{t\in I}\norm{f(t)}_{\Holder{k}}$.
Constants $C,p$ depend only on the dimension.

Let $K(t,a)$ be symmetric and $\Omega(t,a)$ skew-symmetric, continuous in time with
values in $\Holder{1}$. Set
\[
 (Jf)(t)=\int_0^t f(s)\dd s.
\]
Assume, for some $\kappa,\mu>0$,
\begin{align}
 K&\ge\kappa I,\label{cc:holder:eq:elliptic}\\
 \int_I\bigl(K(t,a)z_t\cdot z_t+\Omega(t,a)z_t\cdot z\bigr)\dd t
 &\ge\mu\int_I|z_t|^2\dd t
 \quad(a\in\mathbb T^3,\ z\in H^1_0(I;\R^3)).
 \label{cc:holder:eq:coercive}
\end{align}
Put
\[
 Q=2+\kappa^{-1}+\mu^{-1}+\norm K_1+\norm\Omega_1.
\]

\begin{theorem}\label{cc:holder:linear}
Under \eqref{cc:holder:eq:elliptic}--\eqref{cc:holder:eq:coercive},
let $L\in C(I;\Holder{0})$. There are unique
$z\in C^1(I;\Holder{0})$, $c\in\Holder{0}$ divergence free and
$\nabla\phi\in C(I;\Holder{0})$ satisfying
\[
 Kz_t-J(\Omega z_t)+L=c+\nabla\phi,
 \qquad \diver z=0,\qquad z(0)=z(1)=0.
\]
They obey
\begin{equation}\label{cc:holder:eq:base-linear}
 \norm{z_t}_0\le CQ^p\norm L_0.
\end{equation}
If $K$, $\Omega$ and $L$ are smooth on $I\times\mathbb T^3$, then so are
$z$, $c$ and $\nabla\phi$.
\end{theorem}

Write $v=z_t$ and $z=Jv$. The terminal condition requires $v$ to have
zero time mean at each label. Let $\mathcal Q$ subtract this mean.
The spatial Leray projection $\PP$ removes the pressure gradient, and
$\mathcal Q$ then removes the time-independent field $c$. The two
projections commute. With $\mathcal R=\PP\mathcal Q$, the equation for
the solenoidal velocity becomes
\[
 \mathcal R Bv=-\mathcal R L,\qquad v=\mathcal Rv,\qquad
 Bv=Kv-J(\Omega v).
\]
For a path of zero time mean $w\in H=L^2(I;\R^3)$, we have
$Jw\in H^1_0(I;\R^3)$ and $J^*w(t)=\int_t^1 w(s)\dd s=-Jw(t)$.
Hence at each label $a$,
\[
 \ip{B(a)w}{w}_H
 =\int_I\bigl(Kw\cdot w+\Omega w\cdot Jw\bigr)\dd t
 \ge\mu\norm w_H^2
\]
by \eqref{cc:holder:eq:coercive} applied to $Jw$.
This coercivity controls $v$ in $\Holder{0}_a(L^2_t)$. To obtain
continuity in time and \eqref{cc:holder:eq:base-linear}, we also invert
$v\mapsto\PP(K(t)v)$ at each time using \eqref{cc:holder:eq:elliptic}.

\subsection{Spatial freezing}
Both inversions have the form $\mathcal R(Bv)=f$ in the range of an
orthogonal projection. Freezing $B$ at a label makes this equation
diagonal in the spatial Fourier modes, even when $B$ acts on whole
time paths.
\begin{lemma}\label{cc:holder:freezing}
Let either $H=\R^3$, $\mathcal R=\PP$, or
$H=L^2(I;\R^3)$, $\mathcal R=\PP\mathcal Q$, where
$\mathcal Qv=v-\int_0^1v(t)\dd t$.
Let $B\in\Holder{1}(\mathbb T^3;\mathcal L(H))$ and $\mu>0$. Suppose
\begin{equation}\label{cc:holder:eq:abstract-coercive}
 \ip{B(a)w}{w}_H\ge\mu\norm w_H^2
\end{equation}
for every $a\in\mathbb T^3$ and every $w\in H$, with
$\mathcal Qw=w$ in the second case.
Then $\mathcal RB$ is an isomorphism on $\mathcal R\Holder{0}(\mathbb T^3;H)$, and
\begin{equation}\label{cc:holder:eq:freezing}
 \norm v_{\Holder{0}_a(H)}
 \le C(2+\mu^{-1}+\norm B_{C^1_a\mathcal L(H)})^{p}
            \norm{\mathcal RBv}_{\Holder{0}_a(H)}.
\end{equation}
If moreover $B\in\Holder{k}(\mathbb T^3;\mathcal L(H))$ for an
integer $k\ge1$, the same operator is an isomorphism on
$\mathcal R\Holder{k}(\mathbb T^3;H)$.
\end{lemma}

\begin{proof}
For $\mathcal Rv=v$, testing $f=\mathcal R(Bv)$ against $v$ and using
\eqref{cc:holder:eq:abstract-coercive} gives
\[
 \mu\norm v_{L^2_a(H)}^2
 \le \ip{Bv}{v}_{L^2_a(H)}=\ip{f}{v}_{L^2_a(H)} .
\]
This bound will control the lower-order term in the H\"older estimate.

Freeze $B_0=B(a_0)$ and let $r(n)$ be the symbol of $\mathcal R$.
For $n\ne0$ the Leray symbol projects onto $n^\perp$, and at zero it
is the identity. In the time-dependent case $r(n)$ also includes
$\mathcal Q$. Extend the frozen operator to all of $H$ by
\[
 A_0(n)=r(n)B_0r(n)+I-r(n).
\]
It acts as the identity on the orthogonal complement of $r(n)H$.
Applying \eqref{cc:holder:eq:abstract-coercive} to the real and
imaginary parts of $r(n)w$ gives
\[
 \operatorname{Re}\ip{A_0(n)w}{w}_H
 \ge\mu\norm{r(n)w}_H^2+\norm{(I-r(n))w}_H^2
 \ge\min(1,\mu)\norm w_H^2 .
\]
Lax--Milgram gives $\|A_0(n)^{-1}\|\le\min(1,\mu)^{-1}$
\cite[Section 6.2.1, Theorem 1]{Evans}.

The projection $r$ is smooth on the unit sphere. The formula
$D(A_0^{-1})=-A_0^{-1}(DA_0)A_0^{-1}$ bounds its first five angular
derivatives by a fixed polynomial in $2+\mu^{-1}+\norm B_\infty$. It
follows that
\begin{equation}\label{cc:holder:eq:frozen-bound}
 \norm{A_0(D)^{-1}f}_{\Holder{0}_a(H)}\le C_0\norm f_{\Holder{0}_a(H)},
 \qquad C_0\le C(2+\mu^{-1}+\norm B_\infty)^C.
\end{equation}
Indeed, localize either symbol $A_0^{-1}$ or $A_0^{-1}r$ to frequencies
of size $2^j$. Five integrations by parts on the rescaled annulus give
the Euclidean kernel estimate
\[
 \|\nabla^dG_j(y)\|_{\mathcal L(H)}
 \le CC_0\,2^{(3+d)j}(1+2^j|y|)^{-5},\qquad d=0,1 .
\]
For $d=0$ the uniform $L^1$ bound persists under periodization. Applying
the kernel to a slightly enlarged input block bounds the output block by
$CC_0\,2^{-j/2}\norm f_{\Holder{0}_a(H)}$. Summing differences with one
derivative on the low blocks and the supremum bound on the high blocks
gives \eqref{cc:holder:eq:frozen-bound}.
The zero mode acts on the mean with norm at most $\min(1,\mu)^{-1}$.

Only near $a_0$ is the coefficient difference $B_0-B$ small.
Choose a cutoff $\chi$ equal to one on the ball of radius $r$ about
$a_0$, supported in the ball of radius $2r$, with $0\le\chi\le1$ and
$|\nabla\chi|\le C/r$. Splitting the source $(B_0-B)v$ gives
\begin{align}\label{cc:holder:eq:local}
 v={}&A_0(D)^{-1}f+A_0(D)^{-1}\mathcal R\{\chi(B_0-B)v\}\notag\\
 &+A_0(D)^{-1}\mathcal R\{(1-\chi)(B_0-B)v\}.
\end{align}
The first term is controlled by \eqref{cc:holder:eq:frozen-bound}.
The coefficient $\chi(B_0-B)$ has supremum norm at most
$Cr\norm B_{C^1}$ and Lipschitz constant at most $C\norm B_{C^1}$.
Its $1/2$-H\"older norm is therefore at most $Cr^{1/2}\norm B_{C^1}$.
Put $Q_0=2+\mu^{-1}+\norm B_{C^1}$. Choosing $r=Q_0^{-M}$ with $M$ sufficiently
large bounds the second term in
$\Holder{0}_a(H)$ by $\frac14\norm v_{\Holder{0}_a(H)}$.

For the last term the source vanishes on the ball of radius $r$.
The periodized bounds for $d=1$ sum to $CC_0(|y|^{-4}+1)$ away from
zero. At points within $r/2$ of $a_0$, integrating over torus distances
$|y|\ge r/2$ bounds the gradient by
$CC_0r^{-1}\norm B_\infty\norm v_\infty$.
For pairs of points at distance at most $r/2$, center the cutoff at one
of them. The last term then contributes at most
$CC_0r^{-1/2}\norm B_\infty\norm v_\infty$ to the H\"older quotient.
For other pairs use $Cr^{-1/2}\norm v_\infty$ for the quotient of $v$.
Absorbing the second term of \eqref{cc:holder:eq:local} gives
\[
 \norm v_{\Holder{0}_a(H)}\le C_1
       \bigl(\norm f_{\Holder{0}_a(H)}+\norm v_{L^\infty_a(H)}\bigr),
 \qquad C_1\le CQ_0^C.
\]
Averaging on a ball of radius $\varepsilon^2$ gives
\[
 \norm v_\infty\le\varepsilon\norm v_{\Holder{0}_a(H)}
       +C\varepsilon^{-3}\norm v_{L^2_a(H)}.
\]
Taking $\varepsilon=(2C_1)^{-1}$ and using
$\norm v_{L^2_a(H)}\le\mu^{-1}\norm f_{L^2_a(H)}$ from coercivity
proves \eqref{cc:holder:eq:freezing}.

To obtain existence from this estimate, set $B_s=(1-s)I+sB$.
For $0\le s\le1$ its coercivity is at least $\min(1,\mu)$, so the
estimate is uniform in $s$. When $B\in\Holder{k}$, apply it to
$\partial^\nu v$ for $|\nu|\le k$. Spatial derivatives commute with
$\mathcal R$, so
\[
 \mathcal R B_s\partial^\nu v
 =\partial^\nu(\mathcal R B_s v)-\mathcal R[\partial^\nu,B_s]v .
\]
Only lower derivatives of $v$ occur in the commutator. Induction
therefore gives $\norm v_{\Holder{k}}\le M_k\norm{\mathcal R B_s v}_{\Holder{k}}$
uniformly in $s$. On $\mathcal R\Holder{k}$ the operator at $s=0$ is
the identity. The uniform estimate and the Neumann series extend its
inverse in finitely many steps to $s=1$, for $k=0$ and for every
admissible $k\ge1$.
\end{proof}

\subsection{Recovering the velocity at each time}
\begin{proof}[Proof of Theorem~\ref{cc:holder:linear}]
On $H=L^2(I;\R^3)$ the coefficient $B=K-J\Omega$ satisfies
$\norm B_{C^1_a\mathcal L(H)}\le C(\norm K_1+\norm\Omega_1)$.
With the coercivity obtained from \eqref{cc:holder:eq:coercive},
Lemma~\ref{cc:holder:freezing} gives a unique solution
$v=\PP\mathcal Qv$ of $\PP\mathcal Q(Bv+L)=0$, with
\begin{equation}\label{cc:holder:eq:mixed}
 \norm v_{\Holder{0}_a(L^2_t)}\le C Q^{p}\norm L_0.
\end{equation}
Taking the time mean recovers the field removed by $\mathcal Q$,
\[
 c=\PP\int_0^1(Bv+L)\dd t .
\]
Then $\PP(Bv+L)=c$ in $L^2_t L^2_a$. In particular, for almost every time,
\begin{equation}\label{cc:holder:eq:instant}
 \PP(Kv)=c+\PP J(\Omega v)-\PP L.
\end{equation}
The right-hand side, denoted by $g(t)$, is already continuous in
$\Holder{0}$. Indeed, \eqref{cc:holder:eq:mixed} bounds $c$ in
$\Holder{0}$, and Cauchy--Schwarz on $[s,t]$ bounds the time increment
of $J(\Omega v)$ in the same norm by
$|t-s|^{1/2}\norm{\Omega v}_{\Holder{0}_a(L^2_t)}$.
Together with the hypothesis on $L$, this gives $\norm g_0\le CQ^p\norm L_0$.

For every $t$ let $S(t)u=\PP(K(t)u)$ on $\PP\Holder{0}$.
By \eqref{cc:holder:eq:elliptic}, Lemma~\ref{cc:holder:freezing} with
$H=\R^3$ bounds $S(t)^{-1}$ uniformly by $CQ^p$. Since $K$ is continuous
in time, the resolvent identity gives continuity of $S(t)^{-1}$ in
operator norm. Thus $\widetilde v(t)=S(t)^{-1}g(t)$ is continuous in
$\Holder{0}$ and satisfies \eqref{cc:holder:eq:base-linear}.
For almost every $t$, the difference $u=v(t)-\widetilde v(t)$ is
solenoidal in $L^2_a$ and satisfies $\PP(K(t)u)=0$ by
\eqref{cc:holder:eq:instant}. Testing against $u$ gives
\[
 \kappa\norm u_{L^2_a}^2\le\ip{K(t)u}{u}_{L^2_a}
 =\ip{\PP(K(t)u)}{u}_{L^2_a}=0 .
\]
Thus $v$ has the continuous representative $\widetilde v$ with the same
zero time mean, so $z=Jv$ satisfies both endpoint conditions. At every
time $\PP(Bv+L-c)=0$, hence $Bv+L-c$ is a periodic gradient
$\nabla\phi\in C(I;\Holder{0})$. Uniqueness of the projected solution
$v$ gives uniqueness of $z$, $c$ and $\nabla\phi$.

For smooth data the higher-order assertion of Lemma~\ref{cc:holder:freezing}
first gives $v\in\Holder{k}_a(L^2_t)$, hence $c\in\Holder{k}$ and
$J(\Omega v)\in C(I;\Holder{k})$. The fixed-time inverse then gives
$v\in C(I;\Holder{k})$ for each $k$. Since $K$ is smooth, $S(t)^{-1}$
is smooth in operator norm on $\PP\Holder{k}$. In
\eqref{cc:holder:eq:instant}, $J(\Omega v)$ gains one time derivative
from $v$. Induction gives all time derivatives of $v$ in every $\Holder{k}$
and thus smoothness of $z$ and $\nabla\phi$.
\end{proof}

\section{The frequency cutoff}\label{cc:app:cutoff}
Fix a radial $\psi\in C_c^\infty(\R^3)$ with $0\le\psi\le1$, $\psi=1$ on
$|\zeta|\le1$ and $\supp\psi\subset\{|\zeta|\le2\}$. For $\Lambda\ge1$ let
$S_\Lambda$ multiply the spatial Fourier mode $e^{2\pi i\xi\cdot x}$,
$\xi\in\mathbb Z^3$, by $\psi(\xi/\Lambda)$.
Subsection~\ref{cc:uf:trunc:sec:truncation} takes $u=S_\Lambda V$ as the
next parent. The frequencies of $u$ are at most $2\Lambda$, and those of
$u\cdot\nabla u$ are at most $4\Lambda$. Thus \eqref{cc:uf:trunc:eq:lp}
bounds derivatives of both $u$ and its residual, while
\eqref{cc:uf:trunc:eq:tail} controls the error $u-V$.

\begin{lemma}\label{cc:uf:trunc:lem:cutoff}
$S_\Lambda$ commutes with derivatives, preserves parity and divergence
freeness, and maps into fields supported on frequencies at most
$2\Lambda$. Its operator norms on $L^\infty$ and on $\Holder{0}$ are bounded
in terms of $\psi$ alone. For every integer $m\ge0$ and every
normed-space-valued $f$ supported on frequencies at most $4\Lambda$,
\begin{equation}\label{cc:uf:trunc:eq:lp}
 \sum_{|\beta|=m}\norm{\partial^\beta f}_\infty\le C_m\Lambda^m\norm f_\infty .
\end{equation}
For all integers $0\le s\le M$ and every $f\in\Holder{M}$,
\begin{equation}\label{cc:uf:trunc:eq:tail}
 \norm{(\Id-S_\Lambda)f}_{\Holder{s}}
 \le C_M\Lambda^{s-M}\norm f_{\Holder{M}}.
\end{equation}
The constants $C_m$ and $C_M$ depend only on their index and on $\psi$.
\end{lemma}

\begin{proof}
The symbol is scalar, even and independent of time, so the cutoff
preserves parity and divergence freeness and commutes with differentiation.
Its support gives the frequency bound. For the norm estimates we use
the convolution kernel as in the Euclidean Bernstein inequalities
\cite[Appendix~A, (A.2)--(A.6)]{TaoDispersive}. Let $K=\check\psi$,
a Schwartz function, and set $K_\Lambda(y)=\Lambda^3K(\Lambda y)$.
To compute $S_\Lambda f$,
convolve the periodic extension of $f$ to $\R^3$ with $K_\Lambda$.
Since $\norm{K_\Lambda}_{L^1}=\norm K_{L^1}$, applying convolution
to $f$ and to its spatial increments gives the two operator bounds.
If $f$ has frequencies at most $4\Lambda$, then $f=S_{4\Lambda}f$.
The identity $\norm{\partial^\beta K_{4\Lambda}}_{L^1}
=(4\Lambda)^{|\beta|}\norm{\partial^\beta K}_{L^1}$ gives
\eqref{cc:uf:trunc:eq:lp} also for normed-space-valued $f$.

For \eqref{cc:uf:trunc:eq:tail} put $k=M-s$. When $k=0$, apply the
$\Holder{0}$ bound to the derivatives of $f$ through order $s$.
For $k\ge1$, the kernel has integral one and all its positive-order
moments vanish because $\psi$ is constant near zero. Thus in
$f-S_\Lambda f$ only the remainder after expanding $f(x-y)$ at $y=0$
through degree $k-1$ contributes. Its $\Holder{s}$ norm in $x$ is bounded
by $C_M|y|^k\norm f_{\Holder{M}}$, by Taylor's formula for each derivative
through order $s$ and its spatial increments. Integrating this bound
against $|K_\Lambda(y)|$ gives \eqref{cc:uf:trunc:eq:tail}, since
$\int_{\R^3}|y|^k|K_\Lambda(y)|\dd y
=\Lambda^{-k}\int_{\R^3}|y|^k|K(y)|\dd y$.
\end{proof}

\section{The parameter inequalities}\label{cc:app:parameters}
\begin{proof}[Proof of Lemma~\ref{cc:uf:tower:lem:slow-orders}]
The current order $q_j$ is always admissible in the maximum
\eqref{cc:uf:tower:eq:schedule}. Indeed, $\mathcal B(q_1)=\mathcal B(1)
\le N=\lambda_0$, and for $j\ge2$ the preceding choice gives
$\mathcal B(q_j)\le\lambda_{j-2}\le\lambda_{j-1}$. Thus $q_{j+1}=q_j$
until $\lambda_{j-1}$ reaches $\mathcal B(q_j+1)$, at which point the
order increases by one. Since $N>1$ and $K_j\ge2$, the frequencies tend
to infinity. They eventually exceed every fixed threshold, so
$q_j\to\infty$.

The gain $L_{j+1}$ involves three successive frequency steps. By
\eqref{cc:uf:tower:eq:schedule},
\[
 L_{j+1}=4\log\lambda_{j+2}
        =4K_jK_{j+1}K_{j+2}\log\lambda_{j-1}.
\]
Each exponent is at most $K(q_{j+1}+1)$, since the orders increase
by at most one. Thus \eqref{cc:uf:tower:eq:uniform} follows from
\eqref{cc:par:floor} and the choice of $q_{j+1}$, which give
\[
 4K_jK_{j+1}K_{j+2}\le4K(q_{j+1}+1)^3
 \le\log\mathcal B(q_{j+1})\le\log\lambda_{j-1}. \qedhere
\]
\end{proof}

\section{Finite norms and coefficient estimates}\label{cc:uf:finite:sec:calculus}
In Section~\ref{cc:sec:exact} we estimate the mean and the oscillatory
corrections at every derivative order. We measure label derivatives on
the envelope radius $\ell_*$ and phase derivatives on the profile width
$h$, as chosen in Subsection~\ref{cc:uf:sec:parent-step}.
With these weights the power of the parent bound $\mathsf D$ in the
inverse estimates is independent of the derivative order.

Write $|f|_m=\sum_{|\beta|=m}\sup|\partial^\beta f|$ for the spatial
derivatives of order $m$, taken in the label $a$ or the position $x$
as appropriate. The supremum includes all variables. For a field
$f(t,a,s)$ of period one in $s$ and an integer $n\ge0$, set
\begin{equation}\label{cc:uf:trunc:eq:field-norm}
 \|f\|_{[n]}=
 \sum_{m+j\le n}\ell_*^mh^j|\partial_s^jf|_m,
 \qquad \|f\|_{[n],\gamma}=\|f/\gamma\|_{[n]}.
\end{equation}
The weighted norm is taken on $[t_0,T]$, where $\gamma(t)>0$ measures
the amplitude of the seed at the origin. All norms may be restricted to a closed
subinterval of their time domain.

For a phase-independent field $c(t,a)$, only the terms with $j=0$
occur in \eqref{cc:uf:trunc:eq:field-norm}. The mean inverse of
Theorem~\ref{cc:uf:mean:inverse} uses spatial H\"older norms, so we put
\[
 |c|_{m,1/2}=\sum_{|\beta|=m}\sup_t\norm{\partial_a^\beta c(t)}_{\Holder{0}},
 \qquad
 \|c\|_{[n],1/2}=\sum_{m\le n}\ell_*^m|c|_{m,1/2}.
\]
To pass from supremum bounds to this norm costs one derivative and a
factor $\ell_*^{-1/2}$. Indeed, bounding $|v(x)-v(y)|$ by both
$2\|v\|_\infty$ and $\|\nabla v\|_\infty|x-y|$ gives
$[v]_{\Holder{0}}\le(2\|v\|_\infty\|\nabla v\|_\infty)^{1/2}$.
Applying this to $\partial_a^\beta c$ and using $\ell_*\le1$ gives
\begin{equation}\label{cc:uf:trunc:eq:holder-sup}
 \ell_*^m|c|_{m,1/2}
 \le C\ell_*^{-1/2}\bigl(\ell_*^m|c|_m+\ell_*^{m+1}|c|_{m+1}\bigr),
 \qquad\hbox{hence}\quad
 \|c\|_{[n]}\le\|c\|_{[n],1/2}\le C\ell_*^{-1/2}\|c\|_{[n+1]} .
\end{equation}

\begin{lemma}\label{cc:uf:trunc:lem:fields}
For every $n\ge0$ and fields $f,g$,
\begin{align}
 \sum_i\|\partial_{a_i}f\|_{[n],\gamma}
 &\le C\ell_*^{-1}\|f\|_{[n+1],\gamma},\qquad
 \|\partial_sf\|_{[n],\gamma}
 \le h^{-1}\|f\|_{[n+1],\gamma},
 \label{cc:uf:trunc:eq:deriv}\\
 \|fg\|_{[n],\gamma}
 &\le C_n\|f\|_{[n]}\|g\|_{[n],\gamma}.
 \label{cc:uf:trunc:eq:product}
\end{align}
With $\av f$ the phase average,
\begin{equation}\label{cc:uf:trunc:eq:phase-ops}
 \|\av f\|_{[n],\gamma}\le\|f\|_{[n],\gamma},\qquad
 \|f-\av f\|_{[n],\gamma}\le2\|f\|_{[n],\gamma},\qquad
 \|\partial_s^{-1}f\|_{[n],\gamma}\le C\|f\|_{[n],\gamma}
 \quad(\av f=0).
\end{equation}
For phase-independent fields,
\eqref{cc:uf:trunc:eq:product} also holds with
$\|\cdot\|_{[n],1/2}$ in place of $\|\cdot\|_{[n]}$ and
$\gamma\equiv1$.
\end{lemma}
\begin{proof}
Each label derivative of order $m+1$ occurs at most three times in the
sum over $i$ in \eqref{cc:uf:trunc:eq:deriv}. The phase derivative
shifts $j$ by one in \eqref{cc:uf:trunc:eq:field-norm}. For
\eqref{cc:uf:trunc:eq:product}, apply Leibniz' rule to $f(g/\gamma)$,
using that $\gamma$ depends only on time. For the H\"older version use
the product inequality in $\Holder{0}$.

Phase averaging has supremum norm at most one. For a zero-mean $f$,
subtracting the average of $\int_0^s f(t,a,\sigma)\dd\sigma$ gives its
zero-mean periodic primitive, of supremum norm at most $2\sup|f|$.
The phase operators commute with label and phase derivatives, so these
bounds give \eqref{cc:uf:trunc:eq:phase-ops} at every order.
\end{proof}

The coefficients in Subsections~\ref{cc:uf:mean:sec:finite}
and~\ref{cc:uf:sec:transverse} are independent of the phase. Their
label derivatives satisfy bounds of the form
\begin{equation}\label{cc:uf:trunc:eq:power}
 |c|_m\le C_m\mathsf D^{C(m+1)}\qquad(m\ge0).
\end{equation}
The absolute exponent $C$ may depend on the coefficient.
The next lemma allows us to multiply and invert such coefficients,
and shows how much of the factor $\ell_*^m$ is available when at
least one derivative falls on a coefficient.

\Needspace{4\baselineskip}
\begin{lemma}\label{cc:uf:trunc:lem:coefficients}
Let $c,d$ be smooth phase-independent matrix fields that satisfy
\eqref{cc:uf:trunc:eq:power}.
\begin{enumerate}[label=(\alph*),leftmargin=*]
\item The product $cd$ satisfies \eqref{cc:uf:trunc:eq:power}.
\item If $c$ is invertible and $\sup_{t,a}|c^{-1}|$ is at most a power
of $\mathsf D$ with an absolute exponent, then $c^{-1}$ satisfies
\eqref{cc:uf:trunc:eq:power}.
\item If $\ell_*\le\mathsf D^{-C}$ with the exponent $C$ of
\eqref{cc:uf:trunc:eq:power} for $c$, then
$\|c\|_{[n]}\le C_n\mathsf D^C$ for every $n\ge0$, and
$\ell_*^m|c|_m\le C_m\mathsf D^{2C}\ell_*$ for every $m\ge1$.
\end{enumerate}
\end{lemma}
\begin{proof}
For (a) and (b) take $C$ large enough for the bounds on $c,d$ and,
in (b), on $\sup|c^{-1}|$. Leibniz' rule gives
$|cd|_m\le C_m\mathsf D^{C(m+2)}$, proving (a) with exponent $2C$.
For $e=c^{-1}$, differentiating $ce=I$ gives
\[
 |e|_m\le
 C_m\mathsf D^C\sum_{i=1}^m|c|_i|e|_{m-i}\qquad(m\ge1).
\]
Starting from $|e|_0\le\mathsf D^C$, induction gives
$|e|_m\le C_m\mathsf D^{C(3m+1)}$, proving (b) with exponent $3C$.
For (c) use the exponent $C$ for $c$ in \eqref{cc:uf:trunc:eq:power}.
Then $\ell_*^m|c|_m\le C_m\mathsf D^C(\ell_*\mathsf D^C)^m$.
Since $\ell_*\mathsf D^C\le1$, summing gives the first bound, and
retaining one factor $\ell_*\mathsf D^C$ for $m\ge1$ gives the second.
\end{proof}

For the mean we differentiate the integrated equation
\eqref{cc:uf:mean:integrated} in the unweighted H\"older norm. For
the oscillations we differentiate \eqref{cc:uf:tower:eq:coordinate}
in the $\gamma$-weighted norm. Its unknown takes values in a fixed
$\R^2$, so label differentiation preserves the space of unknowns.
In both equations a derivative falling on a coefficient leaves fewer
derivatives on the solution. The following estimate uses this fact to
retain just one factor of the zero-order inverse bound.

\begin{lemma}\label{cc:uf:trunc:lem:lifting}
Consider a linear equation of the form
\[
 \mathcal Lv=\mathcal L_0v+\sum_\nu T_\nu(c_\nu v)=F
\]
with finitely many phase-independent coefficients $c_\nu(t,a)$.
Assume that $\mathcal L_0$ and the linear operators $T_\nu$ commute
with label and phase derivatives. The equation may hold modulo a fixed
subspace, such as gradients. Assume that differentiation preserves this
subspace and the linear constraints on $v$, and commutes with the maps
giving any prescribed data. Write $f$ for the forcing $F$ together with
these data, with the sum of their norms. Suppose that for every
admissible source one has the estimate
\[
 \|v\|_{[0],\gamma}\le C_{\rm low}\|f\|_{[0],\gamma}.
\]
Let $b_\nu$ be the operator norm of $T_\nu$ from
$\|\cdot\|_{[0],\gamma}$ to the forcing norm. If
\begin{equation}\label{cc:uf:trunc:eq:slow}
 C_{\rm low}\sum_\nu b_\nu\ell_*^m|c_\nu|_m\le C_m\qquad(m\ge1),
\end{equation}
then every smooth solution satisfies, for $n\ge0$,
\[
 \|v\|_{[n],\gamma}\le C_n C_{\rm low}\|f\|_{[n],\gamma}.
\]
Here $C_n$ depends only on $n$ and the bounds $C_m$ in
\eqref{cc:uf:trunc:eq:slow} for $m\le n$. For phase-independent fields
with $\gamma\equiv1$, the same assertion holds using
$\|\cdot\|_{[n],1/2}$ and $|c_\nu|_{m,1/2}$ throughout, with $b_\nu$
taken in the corresponding zero-order norms.
\end{lemma}
\begin{proof}
Phase derivatives commute with the equation. For a label multi-index
$\beta$ and phase order $j$, the differentiated forcing is
\[
 \mathcal L(\partial_a^\beta\partial_s^j v)
 =\partial_a^\beta\partial_s^j F-
 \sum_{\nu}\sum_{0<\beta'\le\beta}\binom{\beta}{\beta'}
 T_\nu\bigl((\partial_a^{\beta'}c_\nu)
          \partial_a^{\beta-\beta'}\partial_s^j v\bigr).
\]
The zero-order estimate applies with the differentiated data. Multiply
by $\ell_*^{|\beta|}h^j$ and sum over $|\beta|+j\le n$. A term with
$i=|\beta'|\ge1$ assigns $\ell_*^i$ to the coefficient and leaves
a norm of order at most $n-i$ on $v$. Thus for $n\ge1$,
\[
 \|v\|_{[n],\gamma}\le C_{\rm low}\|f\|_{[n],\gamma}
 +C_n\sum_{i=1}^n
 \Bigl(C_{\rm low}\sum_\nu b_\nu\ell_*^i|c_\nu|_i\Bigr)
 \|v\|_{[n-i],\gamma}.
\]
By \eqref{cc:uf:trunc:eq:slow}, the sum is bounded by
$C_n\|v\|_{[n-1],\gamma}$. Induction from order zero gives the claim.
For the H\"older norms the same argument uses the product inequality
in $\Holder{0}$.
\end{proof}

In Section~\ref{cc:sec:exact} the zero-order bounds come from
Theorem~\ref{cc:uf:mean:inverse} and \eqref{cc:uf:tower:eq:weighted-forward}.
For $m\ge1$, Lemma~\ref{cc:uf:trunc:lem:coefficients}(c) leaves a
factor $\ell_*$ in $\ell_*^m|c_\nu|_m$. In the H\"older norm,
\eqref{cc:uf:trunc:eq:holder-sup} leaves a factor $\ell_*^{1/2}$,
as in \eqref{cc:uf:mean:thresholds}. In each case $C_{\rm low}\sum_\nu b_\nu$
is bounded by a fixed power of $\mathsf D$. One choice
$\ell_*=\mathsf D^{-A}$ with $A$ sufficiently large therefore gives
\eqref{cc:uf:trunc:eq:slow} at every order.

\section*{Acknowledgments}
T.~B. acknowledges support from the National Science Foundation under
grants DMS-2243205 and DMS-2244879 and from the Simons Foundation through
its Wave Turbulence collaboration. He thanks Diego C\'ordoba and Luis
Mart\'inez-Zoroa for many discussions of their work, which laid the
foundation for this construction.
M.~P.~C.\ acknowledges support at New York University through a Simons
Society of Fellows Junior Fellowship. Simons Foundation International
funds the fellowship, and the Simons Foundation administers it.

We thank Anthropic for supporting the collaboration from which this paper
grew.

\end{document}